\documentclass[a4paper,11pt]{article}
\usepackage{amscd,amssymb,amsmath,amsthm}
\usepackage{rotating}
\usepackage{hyperref}
\usepackage{graphicx,caption,dsfont}
\usepackage{mathtools}

\usepackage{subcaption}
\usepackage{amsfonts,enumerate}
\usepackage{fullpage}
\usepackage[numbers,sort & compress]{natbib}
\usepackage{soul}
\usepackage{xfrac}
\usepackage{xcolor}
\usepackage{epstopdf}
\usepackage{booktabs}
\numberwithin{equation}{section}
\usepackage{xcolor}

\newtheorem{theorem}{Theorem}[section]
\newtheorem{example}{Example}[section]
\newtheorem{proposition}{Proposition}[section]

\newtheorem{note}{ Note}[section]

\newtheorem{remark}{Remark}
\newtheoremstyle{case}{}{}{}{}{}{:}{ }{}
\theoremstyle{case}

\allowdisplaybreaks

\makeatletter
\def\and{%
\end{tabular}%
\begin{tabular}[t]{c}}%
\def\@fnsymbol#1{\ensuremath{\ifcase#1\or 1\or 2\or 3\or
		4\or 5\or 6\or 7\or 8\or 9\else\@ctrerr\fi}}
\makeatother
\vspace{-8ex}
\date{}
\begin{document}
	\baselineskip=0.70cm
	\begin{center}\large{\textbf{ Analytical solutions of higher-dimensional coupled system of nonlinear
				time-fractional diffusion-convection-wave equations}}
	\end{center}
	\vspace{-0.5cm}
	\begin{center}
		  K.S. Priyendhu$^{1}$, P. Prakash$^{1*}$,  M. Lakshmanan$^{2}$\\
		$^{1}$Department of Mathematics,\\
		Amrita School of Engineering, Coimbatore-641112,\\
		Amrita Vishwa Vidyapeetham, INDIA.\\
		$^{2}$Department of Nonlinear Dynamics,\\
		Bharathidasan University,\\
		Tiruchirappalli-620 024, INDIA.\\
		E-mail ID's: vishnuindia89@gmail.com \&\ p$_{-}$prakash@cb.amrita.edu (P Prakash)\\
		\hspace{1.6cm} : priyendhusasikumar@gmail.com (K.S. Priyendhu)\\
		\hspace{1.6cm} : lakshman.cnld@gmail.com (M. Lakshmanan)\\
		$*$-Corresponding Author
	\end{center}
	\begin{abstract}This article develops how to generalize the invariant subspace method for deriving the {analytical solutions} of the {multi-component} $(N+1)$-dimensional coupled nonlinear time-fractional PDEs (NTFPDEs) in the sense of Caputo fractional-order derivative for the first time.
		{Specifically,}  we describe how to systematically find different invariant  product linear spaces with various dimensions for the considered system. Also, we observe that the obtained invariant  product linear spaces help to reduce the {multi-component} $(N+1)$-dimensional coupled  NTFPDEs into a system of fractional-order ODEs, which can then be solved using the well-known analytical methods.
		More precisely, we illustrate the effectiveness and importance of this developed method for obtaining a long list of invariant  product linear spaces for the {multi-component} $(2+1)$-dimensional coupled  nonlinear time-fractional diffusion-convection-wave equations.
		In addition, we have shown how to find different {kinds} of generalized separable {analytical} solutions for a {multi-component} $(2+1)$-dimensional coupled  nonlinear time-fractional diffusion-convection-wave equations  along with the initial and boundary conditions using invariant  product linear spaces obtained. Finally, we provide appropriate graphical representations of some of the {derived} generalized separable {analytical} solutions with various fractional-order values.
	\end{abstract}
\section{Introduction}
The theory of fractional calculus is the generalization of integer-order (classical) calculus, and it  deals with the theoretical development and applications of arbitrary-order (fractional-order) derivatives as well as  arbitrary-order integrals, commonly known as { differ-integrals \cite{Tarasov2013a,Tarasov2011,ionescu2017,bagley1984,Tarasov2013b,Podlubny1999,Diethelm2010,Mainardi1997}}.
{In the literature, there are different approaches of generalization to the  differ-integral definitions.
The most prominent among them is the one based on Cauchy $n$-fold integration, and the widely used definitions of such arbitrary-order derivatives and integrals were proposed by Riemann, Liouville,  Caputo and many others  \cite{Podlubny1999,Diethelm2010}.
Each of the generalizations has its importance in the sense of its advantages in different applications, especially in the fields of science and engineering  \cite{Tarasov2013a,Tarasov2011,ionescu2017,bagley1984,Tarasov2013b,Podlubny1999,Diethelm2010,Mainardi1997}.
Recent studies proved that  arbitrary-order derivatives are more efficient and powerful tools for formulating the equations that govern physical processes, which  {can not only be described by their present state but also by their previous states}  in the history of the process.
This power-law long-term memory property is the central essence of  arbitrary-order derivatives, which helps them to describe the physical processes more accurately and precisely than the integer-order {derivatives \cite{Tarasov2013a,Tarasov2011,ionescu2017,bagley1984,Tarasov2013b,Podlubny1999,Diethelm2010,Mainardi1997,povstenko2013}.}
It is well-known that the Caputo fractional derivative considers the initial conditions that can be interpreted physically in science and engineering applications. Also, we know that the Caputo fractional derivative is  {being utilized} to describe many real-world phenomena with a nonlocal nature that depend   on their past interactions.
{In the literature,}  the Caputo  fractional-order derivative is called  as singular kernel derivative.
Additionally, we note that  the singular kernel Caputo fractional-order derivative  violates some standard properties,    {which include the chain rule, Leibniz rule,}  and  semi-group property. Consequently, the Caputo fractional-order derivative is said to possess nonstandard properties. These nonstandard properties allow singular kernel  Caputo fractional-order derivative to possess  long-term power-law memory property, as well as other fascinating properties such as nonlocal behavior, irreversible nature, and fractal-type properties  \cite{Tarasov2013a,Tarasov2011}. 
Consequently,  {fractional-order differential equations} (FDEs) are used to  formulate many real-world processes, primarily in the fields of physics and {mechanics \cite{Tarasov2013a,Tarasov2011,ionescu2017,bagley1984,Tarasov2013b,Podlubny1999,Diethelm2010,Mainardi1997,povstenko2013}}.
For example, one of the most prominent applications of the FDEs in the sense of Caputo fractional-order derivative can be seen in the mathematical formulation of  the processes that govern anomalous diffusion problems in transport media like heat transfer in solids, thermo-elasticity without energy dissipation, chemical kinetics, biological population dynamics, vorticity in fluids, electro-dynamics,  contamination transport in groundwater, collision in radioactive elements, sorption of organic compounds, ion exchange, geophysical-environmental mechanics, and many more  \cite{ionescu2017,bagley1984,Tarasov2013b,Podlubny1999,Diethelm2010,Mainardi1997,povstenko2013,Povstenko2015}. So, in this work, we consider the Caputo fractional-order derivative for our investigation. }
\subsection{Motivation}
The  non-Fickian diffusion-convection-wave processes are governed by the two relations as follows:
\begin{eqnarray*}
	&	(i) \text{ Conservation law of mass: }&
	\dfrac{\partial u}{ \partial t}=-\bigtriangledown\cdot \mathbf{Q} ,\, 	\&
	\\
\, &(ii) \begin{array}{ll}
		\text{ Constitutive equation \cite{Povstenko2015}: }
	\end{array}
	& 	
	\mathbf{Q}=\left\{\begin{array}{ll}
		\dfrac{1}{\Gamma (\alpha)} \dfrac{\partial}{\partial t} \int_{0}^t(t-y)^{\alpha-1} \mathbf{q} dy ,  \text{ if }\alpha\in(0,1],
		\\
		\dfrac{1}{\Gamma (\alpha-1)}  {\int_{0}^t}(t-y)^{\alpha-2}  \mathbf{q} dy, \text{ if } \alpha\in(1,2],
	\end{array}\right.
\end{eqnarray*}
where $\mathbf{Q}$ is the matter flux such that  $  \mathbf{q}=-\kappa\bigtriangledown u+\mathbf{v}u $.
Here $\kappa$  and  $\mathbf{v}$ are the diffusion coefficient and the velocity of the particle with concentration $u=u(x_1,\dots,x_N,t),$ respectively, and $ \bigtriangledown =\left( \dfrac{\partial }{\partial x_1},\dots,\dfrac{\partial }{\partial x_N}\right) $ is the gradient operator.
Thus, the nonlocal-time dependence of the matter flux described by a power-law long-tail memory kernel gives the time-fractional  diffusion-convection-wave equation  {(TFDCWE)} \cite{Povstenko2015} in the form
\begin{eqnarray}
	\label{CDEQN}
	\dfrac{\partial^\alpha u}{ \partial t^\alpha}= \kappa\bigtriangleup u-\mathbf{v}\cdot \bigtriangledown u,\alpha\in(0,2],
\end{eqnarray}
where $ \bigtriangleup=\sum\limits_{r=1}^N \dfrac{\partial^2 }{\partial x_r^2}$ denotes the Laplacian and $	\dfrac{\partial^{\alpha} }{\partial t^{\alpha}}(\cdot)$ is the  {Caputo fractional-order derivative on the time variable}, $t>0$ of order $ \alpha,\alpha \in(0,2],$ which is defined  \cite{Podlubny1999,Diethelm2010} as
\begin{equation}
	\label{caputo}
	\dfrac{\partial^{\alpha} u}{\partial t^{\alpha}}=\left\{\begin{array}{ll}
		\dfrac{1}{\Gamma(k-{\alpha})}
		\int_{0}^{t}(t-y)^{k-(\alpha+1)}
		{\left( \dfrac{\partial^k  u(x_1,\dots,x_N,y)}{\partial y^k}\right) }dy,
		\\ \text{ if } \alpha\in(k-1,k], k=1,2,
		\\
		\dfrac{\partial^k  u}{\partial t^k}, \quad \text{ if } k=\alpha\in\{1,2\}.
	\end{array}
	\right.
\end{equation}
It is interesting to observe that for $\alpha=1,$ the equation \eqref{CDEQN} is equivalent to the classical integer-order diffusion-convection equation. For $v=0$ (neglecting convective terms) and $\alpha=2,$ the equation \eqref{CDEQN} represents the classical wave equation. Additionally, when $v=0$ in \eqref{CDEQN}, it denotes {sub-diffusion \cite{Tarasov2013b,Povstenko2015,bakkyaraj2015,yu2015}} if $\alpha\in(0,1)$ and {super-diffusion \cite{Tarasov2013b,Povstenko2015,bakkyaraj2015,yu2015}} if $\alpha\in(1,2)$.

{Motivated by the above studies, we  formulate a {multi-component} $(2+1)$-dimensional  {TFDCWEs}   using the following governing equations:
		$$\begin{aligned}
		(i) & \text{ Conservation law of mass: }
		\dfrac{\partial u_s}{ \partial t}=-\bigtriangledown\cdot \mathbf{Q_s}
		,\\
		(ii) &
		\text{ Constitutive equation: }
	\\&	\mathbf{Q_s}=\left\{\begin{array}{ll}
			\dfrac{1}{\Gamma (\alpha_s)} \dfrac{\partial}{\partial t} \int_{0}^t(t-y)^{\alpha_s-1} (-\kappa_s \bigtriangledown u_s(x_1,x_2,y)+{\mathbf{v}_s}\ u_s(x_1,x_2,y) ) dy , 
			\\ \text{ if }\alpha_s\in(0,1],
			\\
			\dfrac{1}{\Gamma (\alpha_s-1)}  {\int_{0}^t}(t-y)^{\alpha_s-2}  (-\kappa_s\bigtriangledown u_s(x_1,x_2,y)+{\mathbf{v}_s}\ u_s(x_1,x_2,y) ) dy, 
			\\\text{ if } \alpha_s\in(1,2],
		\end{array}\right.
	\end{aligned}	$$
	where $\mathbf{Q}_s$ is the matter flux,  $\kappa_s$  and  $\mathbf{v}_s$ are the diffusion coefficient and the velocity of the $s$-th particle with concentration $u_s=u_s(x_1,x_2,t),s=1,2,$ respectively, and $ \bigtriangledown =\left( \dfrac{\partial }{\partial x_1},\dfrac{\partial }{\partial x_2}\right) $ is the gradient operator.
	Substituting  $(ii)$  in $(i),$ we obtain the following {multi-component} $(2+1)$-dimensional  {TFDCWEs}:
	\begin{eqnarray}
		\dfrac{\partial^{\alpha_s} u_s}{ \partial t^{\alpha_s}}=-\kappa_s\bigtriangleup u_s(x_1,x_2,t)+\bigtriangledown \cdot {\mathbf{v}_s}\ u_s(x_1,x_2,t) ,\alpha_s\in(0,2],s=1,2
		,\label{2+1modelcd}
	\end{eqnarray}
	where $\bigtriangleup =\dfrac{\partial^2 }{\partial x_1^2}+\dfrac{\partial^2 }{\partial x_2^2} $ denotes the Laplacian operator and  $	\dfrac{\partial^{\alpha_s} }{ \partial t^{\alpha_s}}(\cdot)$ is the Caputo fractional derivative \eqref{caputo}.}
{Recent studies show that the significant advantage of diffusion models  with time-fractional derivatives is that, in contrast to the classical integer-order diffusion models, they can accurately describe the nonlinear time-dependency of anomalous diffusion processes.
Multi-component $(2+1)$-dimensional time-fractional diffusion equations, widely used in the literature, describe many complex physical models involving more than one diffusive particle or when both heat conduction and diffusion co-occur. For example, Zhang et al. \cite{zhang2019} have studied  the coupled heat conduction and moisture diffusion in a hollow cylindrical pipe in hygrothermoelastic environments using the following {multi-component} $(2+1)$-dimensional  time-fractional PDEs
	\begin{eqnarray}
		\begin{aligned}
			\label{2+1diffzhang}
			\dfrac{\partial^{\alpha_1} T}{ \partial t^{\alpha_1}}-\beta_1	\dfrac{\partial^{\alpha_1} C}{ \partial t^{\alpha_1}}=&\kappa_1\bigtriangleup T,
			\\
			\dfrac{\partial^{\alpha_2} C}{ \partial t^{\alpha_2}}-\beta_2	\dfrac{\partial^{\alpha_2} T}{ \partial t^{\alpha_2}}=&\kappa_2\bigtriangleup C,\alpha_s\in(0,2],s=1,2,
		\end{aligned}
	\end{eqnarray}
	where
	$T=T(x_1,x_2,t)$ represents the temperature, and $C=C(x_1,x_2,t)$ is the concentration of water vapor in a porous media.
	In the above system \eqref{2+1diffzhang},
	$\beta_s=\dfrac{\gamma_s a_s}{\gamma_s a_j+c_s},j\neq s,j,s=1,2,\kappa_1=\dfrac{\lambda_1}{\gamma_2(c_1+\gamma_1 a_2)},$ and $\kappa_2=\dfrac{\lambda_2c_2}{c_2+\gamma_2 a_1},$  where  the arbitrary constants $\gamma_1,\gamma_2,c_1,c_2,$ $\lambda_1$ and $\lambda_2$ represent the
	amount of heat released,
	the density of the material,
	the specific heat capacity at constant pressure,
	the volume fraction of moisture,
	heat conduction and diffusion coefficients, respectively, and $a_s,s=1,2,$ are material constants.
	We note that under the substitution $T-\beta_1 C =u_1$ and $C-\beta_2 T =u_2,$ the above system \eqref{2+1diffzhang} gets reduced into the following {multi-component} $(2+1)$-dimensional coupled  time-fractional diffusion-wave equations:
	\begin{eqnarray}
		\begin{aligned}
			\label{2+1diffzhang-reduced}
			\dfrac{\partial^{\alpha_1} u_1}{ \partial t^{\alpha_1}}=&\kappa_1\bigtriangleup (\tau_{0}u_1+\tau_{1}u_2),
			\\
			\dfrac{\partial^{\alpha_2} u_2}{ \partial t^{\alpha_2}}=&\kappa_2\bigtriangleup (\tau_{2}u_1+\tau_{0}u_2),\alpha_s\in(0,2],s=1,2,
		\end{aligned}
	\end{eqnarray}
	where $\tau_{0}=\dfrac{1}{1-\beta_1\beta_2},$ and $\tau_{s}=\dfrac{\beta_s}{1-\beta_1\beta_2},s=1,2.$
	Another physical application of the above type of system \eqref{2+1modelcd} was discussed by   Yang et al.  \cite{yang2019} to describe the anomalous diffusion of  sodium ( Na$+$) and  potassium (K$^+$) ions in a neuronal phenomenon  as follows:
	\begin{eqnarray}
		\begin{aligned}
			\label{2+1diffyang}
			\dfrac{\partial^{\alpha} u_s}{ \partial t^{\alpha}}=&\sum\limits_{i=1}^2	
			{\frac {\partial }{\partial x_i}}
			\left[ (\tau_1 u_s+\tau_0 ) {\frac {\partial u_{{1}} }{	\partial x_i}}
			+(\tau_1 u_s+\tau_0 ){\frac {\partial u_{{2}} }{
					\partial x_i}} \right]
			-I_s,\alpha\in(0,1],s=1,2,
		\end{aligned}
	\end{eqnarray}
	where $ \dfrac{\partial^{\alpha} }{ \partial t^{\alpha}}(\cdot)$ is the Caputo fractional derivative \eqref{caputo}, $\tau_0,\tau_1$ are known constants, $u_s=u_s(x_1,x_2,t),$ $s=1,2, $ represent the concentration of  Na$+$ and K$^+$ ions with the current density $I_s,s=1,2.$
They have developed the Jacobi special collocation method for the above system \eqref{2+1diffyang}, also known as the coupled time-fractional Nernst-Plank equations \cite{yang2019}, to find an approximate solution of the above system \eqref{2+1diffyang}  along with the reduction$r=\sqrt{x_1^2+x_2^2}.$

It is well-known that the essential characteristics of a physical process can be analyzed using its long-term asymptotic behavior and  qualitative-quantitative properties.
These properties are described exactly by the analytical solutions of the corresponding governing FDEs.
However, due to the  above-mentioned mysterious properties of fractional derivatives, the task of finding {analytical (exact)}  solutions for  nonlinear FDEs is  difficult and challenging, particularly  for nonlinear fractional PDEs.
Further, there is still no well-defined analytical method available to find {analytical} solutions to FDEs in general.
Due to the rigorous investigations in  recent years,  some  numerical and analytical methods have been developed to find solutions for different classes of FDEs which include the Adomain decomposition method \cite{gejji2005,jafari2006,gafiychuk2009},  {the exp-function method \cite{bekir2013,bekir2017}}, the sub-equation method \cite{bekirmma}, the generalized differential transform method \cite{odibat2008}, 
the hybrid method \cite{ma2020}, the separation of variable methods \cite{wu2018,ru2020,rui2022,ruinonlineardyn2022}, the  homotopy analysis method \cite{jafari2009b}, the homotopy perturbation method \cite{jafari2007}, the Lie symmetry analysis {method \cite{bakkyaraj2015,yu2015,gazizov2009,sahadevan2012,sahadevan2017b,prakash2017,nass2019,bakkyaraj2020,sahadevan2019,prakash2021,sethu2021}},
and the invariant subspace {method \cite{gazizov2013,harris2013,sahadevan2015,sahadevan2016,sahadevan2017a,choudary2017,harris2017,harris2014,rui2018,choudary2019b,prakash2020a,prakash2024,prakash2020b,gara2021,choudary2019a,kader2021,prakash2022a,prakash2022b,priyendhu2023,prakash2023,Polyanin2022,chu2022,weng2018,liu2020},
	to name a few.
	
	Among these  analytical methods,  the invariant subspace method has proven to be an effective and powerful analytic tool for finding generalized separable {analytical} solutions for many classes of nonlinear fractional PDEs.
	This method was initially proposed by Galaktionov and Svirshchevskii \cite{Galaktionov2007}, which is also addressed as the generalized separation of the variables method or the Titov-Galaktionov method \cite{Polyanin2022}, for finding generalized separable  {analytical} solutions of  {nonlinear integer-order evolution PDEs.}
	Later, the potential of this method was further investigated by Ma \cite{ma2012a}, Zhu and Qu \cite{zhu2016}, and many others \cite{ma2012b,ma2012c,ye2014,qu2009,song2013,prakash2019,polyanin2014,Polyanin2022} to solve many important classes of  nonlinear coupled  and scalar  integer-order PDEs.
	In 2013, Gazizov and Kasatkin \cite{gazizov2013} extended the invariant subspace method to solve NTFPDEs.
	{Recently,} this method was extended to solve many classes of  {nonlinear scalar} fractional PDEs in {$(1+1)$-dimensions \cite{harris2013,sahadevan2015,sahadevan2016,sahadevan2017a,choudary2017,harris2017,harris2014,rui2018,choudary2019b,prakash2020a,prakash2020b,Polyanin2022,chu2022,gara2021,weng2018,liu2020}.}
	Moreover, this method has been successfully extended and applied to find generalized separable  {analytical} solutions for a $(1+1)$-dimensional coupled  nonlinear fractional PDEs  \cite{sahadevan2017b,sahadevan2017a,choudary2019b,prakash2020a,prakash2024}. In 2019, Sangita et al. \cite{choudary2019a}  applied this method  to find generalized separable {analytical} solutions for the $(2+1)$-dimensional coupled  NTFPDEs.
	\subsection{Novelty of our work}
	We believe that no one has developed the invariant subspace method to find generalized separable analytical solutions for multi-component $(N + 1)$-dimensional coupled NTFPDEs as far as we know.
So, we describe for the first time the generalization of the invariant subspace method for deriving different kinds of analytical solutions for multi-component (N + 1)-dimensional coupled NTFPDEs.
	The main aim of our work is to systematically develop the invariant subspace method for {obtaining}  different types of invariant  product linear spaces and  {analytical} solutions for  {multi-component} coupled NTFPDEs mentioned above.
	Here, it is noteworthy  that the obtained finite-dimensional invariant  product linear spaces help to reduce the {multi-component} $(N+1)$-dimensional coupled  NTFPDEs
	{into a system} of fractional ODEs (FODEs). More precisely, the importance of the developed invariant subspace method has been  {discussed through finding the  {analytical} solutions in the generalized separable form} for the following  initial value problems of the {multi-component} $(2+1)$-dimensional coupled  nonlinear TFDCWEs
	\begin{eqnarray}
		\label{cdsystemeqn}
		\begin{aligned}
			\left(
			\dfrac{\partial^{\alpha_1} u_1}{\partial t^{\alpha_1}}, \dfrac{\partial^{\alpha_2} u_2}{\partial t^{\alpha_2}}
			\right)
			& =
			\mathbf{F}\left(\mathbf{U}\right)\equiv(	F_1(u_1,u_2),	F_2(u_1,u_2)),
			{{\alpha_1,\alpha_2}\in(0,2],}
		\end{aligned}
	\end{eqnarray}
	 {together with  initial conditions:}
	\begin{eqnarray}
		\begin{aligned}
			\label{ic}
			&(i) \, u_s(x_1,x_2,0)= \omega_s(x_1,x_2)
			\,\text{ if }\  {\alpha_s}\in(0,1],s=1,2,\\
			&(ii)\,	u_s(x_1,x_2,0)= \omega_s(x_1,x_2),\, \,
			{\dfrac{\partial u_s}{\partial t}}\big{|}_{t=0}=\vartheta_s(x_1,x_2)
			\,	\text{if}\  {\alpha_s}\in(1,2],s=1,2,
		\end{aligned}
	\end{eqnarray}
	where the components  	
	$F_s(u_1,u_2),s=1,2,$ are considered as follows:
	\begin{eqnarray}
		\label{componentscdsystem}
		\begin{aligned}
			F_s(u_1,u_2)=&
			\sum\limits_{i=1}^2	
			{\frac {\partial }{\partial x_i}}
			\left(  \beta_{si}(u_1,u_2) {\frac {\partial u_{{1}} }{	\partial x_i}}
			+ \rho_{si} (u_1,u_2){\frac {\partial u_{{2}} }{
					\partial x_i}}  \right)
		\\&	+\sum\limits_{i=1}^2 \left( \zeta_{si}(u_1,u_2) {
				\frac {\partial u_{{1}} }{\partial x_i}}
			+ \upsilon_{si} (u_1,u_2){\frac {\partial u_{{2}} }{\partial x_i}}\right)
			.
		\end{aligned}
	\end{eqnarray}
	Here $\mathbf{U}=(u_1,u_2),u_s=u_s(x_1,x_2,t),x_i\in\mathbb{R},t>0,$  the functions $\beta_{si}(u_1,u_2),\rho_{si} (u_1,u_2)$ denote the diffusion coefficients and  $\zeta_{si}(u_1,u_2),$ $\upsilon_{si} (u_1,u_2),s,i=1,2, $ represent the convective terms of the physical system.
	We observe that the obtained system \eqref{2+1diffzhang-reduced} is a particular case of the considered coupled NTFPDEs \eqref{cdsystemeqn}.	Additionally, we would like to point out that the above system \eqref{2+1diffyang} with negligible current density (that is, $I_s=0$) is a particular case of the discussed coupled NTFPDEs \eqref{cdsystemeqn}.
	In the literature, finding the {analytical} solutions of the coupled  diffusion-convection-wave equations has gained considerable attention because of its widespread application as a mathematical formulation for the dynamics of long-range interactions between more than one element or species in  physical processes like  convection-dominated transport media problems.
	Some important classes of coupled  diffusion-convection equations include the coupled Burgers equations, the hydrodynamic model for semiconductor device simulation, magneto-hydrodynamics, shallow water equations, and the Navier-stokes equations \cite{espinov,wu2019,evangelista2018}. In \cite{cockburn1998}, LDG (local discontinuous Galerkin) method was used efficiently to study {analytical} solutions for the integer-order coupled  diffusion-convection equations \cite{cockburn1998}  in the form
	\begin{equation}\label{DGM}
		\dfrac{\partial \mathbf{U}}{\partial t}=\bigtriangledown \cdot \mathbf{F}(\mathbf{U},D\mathbf{U}),
	\end{equation}
	where $\mathbf{U}=(u_1,u_2,\dots,u_m),
	u_s=u_s(x_1,x_2,\dots,x_N,t),x_r\in\mathbb{R},t>0,r=1,2,\dots,N, s=1,2,\dots,m,N,m\in\mathbb{N},$ $ \bigtriangledown$ is the gradient and $D$ is the total derivative of $\mathbf{U}$ of order one.
	Note that for $\alpha=1,$ the above-coupled equations \eqref{cdsystemeqn} belongs to the class of  given coupled  diffusion-convection equations \eqref{DGM} with $m=N=2.$
	\subsection{Structure of the article}
	 {This article  is structured as follows:}
	In  section \ref{sec4},   we develop the invariant subspace method  {to identify} different types of invariant product linear spaces and generalized separable analytical solutions  {for} {multi-component} $(N+1)$-dimensional coupled  NTFPDEs. In addition,  we systematically  {discuss} how to construct $N$-different types of invariant   product linear spaces for the system mentioned above with their invariance criteria.
	Section \ref{sec5}  briefly describes how to systematically  {obtain} various kinds of  invariant  product linear spaces to the {multi-component} $(2+1)$-dimensional coupled  nonlinear  {TFDCWEs} \eqref{cdsystemeqn}  with cubic power-law nonlinearity using the developed invariant subspace method. Also, we present a  {comprehensive} explanation of how to find generalized separable  {analytical} solutions for the considered coupled equations \eqref{cdsystemeqn} equipped {with initial and boundary conditions.} Additionally, we show the graphical representations of some of the obtained generalized separable {analytical} solutions with various $\alpha$-values. {{Finally, section \ref{conclusion} discusses} a brief review and concluding  remarks on the developed method with their obtained results.}
	\section{{Multi-component} $(N+1)$-dimensional coupled  NTFPDEs: Invariant subspace method }\label{sec4}
	This section presents the systematic procedure for deriving different kinds of {analytical} solutions of the {multi-component} $(N+1)$-dimensional coupled  NTFPDEs using the invariant subspace method.
	Let us consider the generalized {multi-component} $(N+1)$-dimensional coupled  NTFPDEs  as follows:
	\begin{eqnarray}
		\begin{aligned}
			\label{Nd1.1}
			\dfrac{\partial^{\alpha_s} u_s}{\partial t^{\alpha_s}}=&{{{A}}_s}(\mathbf{U})
		\\	=&
			{{{A}}_s}\left(x_1,x_2,\dots,x_N,u_1,\dots, u_m,\dfrac{\partial u_1}{\partial x_1},\dots,\dfrac{\partial^{q_s} u_1}{\partial x_1^{i_1}\partial x_2^{i_2}\cdots\partial x_N^{i_{N}}},\dots,
				\right.	\\ & 
			\left. \dfrac{\partial u_m}{\partial x_{1}},
	\dots,\dfrac{\partial^{q_s} u_m}{\partial x_1^{i_1}\partial x_2^{i_2}\cdots\partial x_N^{i_{N}}}\right),\alpha_s>0, s=1,2,\dots,m,
		\end{aligned}
	\end{eqnarray}
	where    ${{{A}}_s}(\mathbf{U}) $ is the $q_s$-th order   sufficiently smooth differential  operator  of the $m$-component vector function
	$\mathbf{U}=(u_1,\dots,u_m),$   $u_s=u_s(x_1,\dots,x_N,t),$ $\sum\limits_{j=1}^N{i_j}=q_s,$  $i_j\in\{0,1,\dots,q_s\},$ $q_s\in\mathbb{N},$ $j=1,2,\dots,N,$  and $ \dfrac{\partial^{\alpha_s} }{\partial t^{\alpha_s}}(\cdot)$ denotes the  Caputo derivative \eqref{caputo} of order $ \alpha_s>0, s=1,2,\dots,m.$
	
	In this work, we develop the invariant subspace method  for the above-coupled NTFPDEs \eqref{Nd1.1}
	to derive possible types of {analytical} solutions in the  generalized separable form  based on ${w_{s,i}}(x_1,x_2,\dots,x_N)$ as	$\mathbf{U}=(u_1,u_2,\dots,u_m) $ with $s$-th component function
	\begin{eqnarray}\label{expected-solutionform}
		\begin{aligned}
			u_s:= u_s(x_1, \dots,x_N,t)=
			\sum\limits_{i=1}^{k_s}	
			{C_{s,i}}(t)
			{w_{s,i}}(x_1, \dots,x_N) ,k_s\in\mathbb{N},
			s=1,2,\dots,m,
		\end{aligned}
	\end{eqnarray}
	where the unknown functions ${C _{s,i}}(t)$ and
	$  {w_{s,i}}(x_1,x_2,\dots,x_N) $
	are to be obtained that depend on the considered coupled NTFPDEs  \eqref{Nd1.1} and its invariant  product linear spaces.
	
	Let us first  define the  product linear space $\mathbf{Y}_n$ as
	\begin{eqnarray}\begin{aligned}
			\label{product}
			\mathbf{Y}_n&
			 	= \mathit{Y}_{1,n_1}\times\mathit{Y}_{2,n_2}\times \cdots\times \mathit{Y}_{m,n_m}
		\\&	=\{\mathbf{y}=(y_1,y_2,\dots,y_m) \ {|} \, y_s\in \mathit{Y}_{s,n_s} ,s=1,2,\dots,m \},
		\end{aligned}
	\end{eqnarray}
	along with  linear component spaces
	$ \mathit{Y}_{s,n_s} $ of dimension $n_s (=\text{dim}( \mathit{Y}_{s,n_s}) <\infty),$  for each $s=1,2,\dots,m ,$ under the component-wise addition  `{$ +$}', and  scalar multiplication  `{$ \cdot$}' as
	\begin{eqnarray*}
		\begin{aligned}
			&	\mathbf{y}+\mathbf{z}=(y_1 + z_1,\dots,y_m +z_m),
			\text{ and }	\mu \cdot 	\mathbf{y} =(\mu \cdot y_1,\ldots,\mu\cdot y_m),
		\end{aligned}
	\end{eqnarray*}
	respectively, where $\mathbf{y}=(y_1,\dots,y_m),\mathbf{z}=(z_1,\dots,z_m)\in \mathbf{Y}_n,$ and $ \mu\in\mathbb{R}.$
	The following proposition explains how to determine the dimension of  {the above-mentioned} product linear space  \eqref{product}.
	\begin{theorem}
		\cite{Axler20142014}\label{axlerprop}
		Let $  \mathbf{Y}_n= \mathit{Y}_{1,n_1}\times\mathit{Y}_{2,n_2}\times \cdots\times \mathit{Y}_{m,n_m}$ given in  \eqref{product} be the  product linear space with the $n_s$-dimensional linear component  spaces  $\mathit{Y}_{s,n_s}, $  where  $n_s (=\text{dim}( \mathit{Y}_{s,n_s})) <\infty,s=1,2,\dots,m. $ Then, the dimension of the   product linear space $\mathbf{Y}_n$ is $ n=n_1+n_2+\dots+n_m.$ That is,
		$\text{dim}(\mathbf{Y}_n)= \sum\limits_{s=1}^m\text{dim}( \mathit{Y}_{s,n_s}).$
	\end{theorem}
	Next, we discuss the invariant  product linear spaces for the given coupled NTFPDEs \eqref{Nd1.1} with their invariant criteria. Also,  we present the systematic algorithmic study  to derive generalized separable  {analytical} solutions of the form \eqref{expected-solutionform} using the invariant  product linear spaces.
\subsection{Finding invariant  product linear spaces and  generalized separable  {analytical} solutions of the coupled NTFPDEs \eqref{Nd1.1}}
	Let us consider the vector differential operator $\mathbf{\hat{A}}[\mathbf{U}]=\left( {{A}}_1(\mathbf{U}),{{A}}_2(\mathbf{U}),\dots {{A}}_m(\mathbf{U})\right) ,$ where ${{A}}_s(\mathbf{U})$ is the nonlinear  differential operator   given in \eqref{Nd1.1} for $s=1,2,\dots,m$.
	Now, we introduce the   product linear space
	${\mathbf{\hat{W}}_k}$ as
	\begin{eqnarray}
		\begin{aligned}
			\label{Nd2}
			{\mathbf{\hat{W}}_k }=& {{W}_{1,k_1}}\times {{W}_{2,k_2}}\times {{W}_{3,k_3}}\times \cdots \times {{W}_{m,k_m}}
			\\	=&\Big\{{\mathbf{w}(x_1,\dots,x_N)}=({w_1}(x_1,\dots,x_N),\dots,{w_m}(x_1,\dots,x_N))
			\big{|}
			\\& {w_s(x_1,\dots,x_N)}\in  {{W}_{s,k_s}},  s=1,2,\dots,m
			\Big\},
		\end{aligned}
	\end{eqnarray}
	where  ${{W}_{s,k_s}}$ are $k_s$-dimensional linear component  spaces for $s=1,2,\dots,m.$ Additionally, from the Proposition \ref{axlerprop}, we note that   dim$ ({\mathbf{\hat{W}}_k})=k=\sum\limits_{s=1}^{m}k_s(<\infty).$
	Thus, we  say that  ${\mathbf{\hat{W}}_k}$  is an {invariant finite-dimensional product linear space} for the nonlinear vector differential operator
	$\mathbf{\hat{A}}(\mathbf{U}) $	if either one of the following invariant conditions holds:
	\begin{itemize}
		\item[(a)]	For every $ \mathbf{U}=(u_1,\dots,u_m) \in {\mathbf{\hat{W}}_k} , $  $ \mathbf{\hat{A}}[\mathbf{U}]=\left( {{A}}_1(\mathbf{U}),\dots {{A}}_m(\mathbf{U})\right)\in {\mathbf{\hat{W}}_k} $ or $ \mathbf{\hat{A}}[ {\mathbf{\hat{W}}_k} ] \subseteq  {\mathbf{\hat{W}}_k} . $
		\item[(b)] If $  \mathbf{U} \in {\mathbf{\hat{W}}_k}, $ then $ {{{A}}_s}(\mathbf{U})\in {{W}_{s,k_s}} , $ or
		${{{A}}_s}({\mathbf{\hat{W}}_k})= {{{A}}_s}( {{W}_{1,k_1}}\times {{W}_{2,k_2}}\times  \cdots \times {{W}_{m,k_m}} )\subseteq {{W}_{s,k_s}},\, \forall u_s\in  {{W}_{s,k_s}},s=1,2,\dots,m.$		
	\end{itemize}
	Now, let  $ B_{s}=\left\{{w_{s,1}(x_1,x_2,\dots,x_N)},{w_{s,2}(x_1,x_2,\dots,x_N)},\dots, {w_{s,k_s}(x_1,x_2,\dots,x_N)}	\right\} $ be a basis for  the linear component spaces $ {{W}_{s,k_s}},$
	then $ {{W}_{s,k_s}},s=1,2,\dots,m,$ can be written  as follows:
	\begin{eqnarray}
		\begin{aligned}
			\label{sNvs}
			{{W}_{s,k_s}}
			&	=\text{Span}\left(B_{s}\right)
		\\&	=\Big\{{w_s(x_1, \dots,x_N)} =\sum\limits_{i=1}^{k_s}
			{\mu_{s,i}} \ {w_{s,i}(x_1, \dots,x_N)}
			\Big{|} {\mu_{s,i}}\in \mathbb{R},i=1,2,\dots,k_s\Big \}.
		\end{aligned}
	\end{eqnarray}
	Thus, using the basis  $ B_{s} $  of the linear component spaces ${{W}_{s,k_s}},s=1,2,\dots,m,$ we define the basis $ \mathbf{\hat{B}} $ for the given finite-dimensional  product linear space ${\mathbf{\hat{W}}_k}$  as follows:
	$$
	\begin{aligned}
		\mathbf{\hat{B}}=& \Big\{({{w_{1,1}}},0,\dots,0),\dots,({{w_{1,k_1}}},0,\dots,0),\dots,(0,\dots,0, {w_{m,1}}),\dots,(0,\dots,0, {w_{m,k_m}})
		\Big{|}
	\\&	{w_{s,i}}={{w_{s,i}}(x_1,x_2,\dots,x_N)}\in B_{s},
		i=1,2,\dots,k_s,s=1,2,\dots, m\Big\}.
	\end{aligned}
	$$
	Hence, the invariant criteria  for the invariant  product linear  space ${\mathbf{\hat{W}}_k}$ of $ \mathbf{\hat{A}}(\mathbf{U})$  can be obtained as follows:
	\begin{equation}\label{Nd1.4}
		\begin{aligned} &\mathbf{\hat{A}}\big(\mathbf{U}\big)=\Big({{{A}}_1}(\mathbf{U}),\dots,{{{A}}_m}(\mathbf{U})\Big)
		\text{ with }
		\\&		{{{A}}_s}(\mathbf{U})=\sum\limits_{i=1}^{k_s}	{{{\Lambda}}_{s,i}({\lambda_{1,1}},\dots,{\lambda_{m,k_m}})}\ {{w_{s,i}}(x_1,x_2,\dots,x_N)} ,
				\end{aligned}
		\end{equation}
	 { whenever }
			$\mathbf{U}=(u_1,u_2,\dots,u_m) $   along with $	u_s=\sum\limits_{i=1}^{k_s}	{\lambda_{s,i}} \ {{w_{s,i}}(x_1,x_2,\dots,x_N)} ,	{\lambda_{s,i}} \in\mathbb{R},$	
	where  the coefficients $ {{{\Lambda}}_{s,i}}(\cdot),$ $i=1,2,\dots,{k_s},$ $s=1,2,\dots,m, $ are obtained from the expansion corresponding  to the  basis $ \mathbf{\hat{B}}. $
	Thus, we can obtain the analytical solution \eqref{expected-solutionform}  for the considered  coupled NTFPDEs \eqref{Nd1.1} whenever the invariant criteria \eqref{Nd1.4} hold.

	Now, we state  how  to reduce the given {multi-component} $(N+1)$-dimensional coupled NTFPDEs \eqref{Nd1.1}  into a system  of time-fractional ODEs under  the  invariant  product linear space \eqref{Nd2}.
	\begin{theorem}
		\label{thm3}
		Let  $
		\mathbf{\hat{A}}[\mathbf{U}]=\left( {{A}}_1(\mathbf{U}),{{A}}_2(\mathbf{U}),\dots {{A}}_m(\mathbf{U})\right)
		$
		be  the vector differential operator, which admits the  {invariant finite-dimensional product linear space} ${\mathbf{\hat{W}}_k }= {{W}_{1,k_1}}\times {{W}_{2,k_2}}\times \cdots \times {{W}_{m,k_m}},$ where  $ {{W}_{s,k_s}} $ are the $k_s$-dimensional  linear component spaces defined in \eqref{sNvs} for $s=1,2,\dots,m$.
		Then, there exists an {analytical} solution for the given {multi-component} $(N+1)$-dimensional coupled NTFPDEs \eqref{Nd1.1}   in the  generalized separable  form  $\mathbf{U}=(u_1, \dots,u_m)$ with $u_s=u_s(x_1, \dots,x_N,t)$ as
		\begin{equation}\label{Nsoln}
			u_s=\sum\limits_{i=1}^{k_s}	{{{C}}_{s,i}(t)}{{w_{s,i}}(x_1,x_2,\dots,x_N)} ,	\end{equation}
		and the  functions ${{{C}}_{s,i}(t)} $ satisfying the following system of time-fractional ODEs
		\begin{eqnarray}
			\begin{aligned}
			\label{Nsysfodes}
			\dfrac{d^{\alpha_s}{C_{s,i}(t)}}{dt^{\alpha_s}} =&{\Phi_{s,i}}\Big( {C_{1,1}(t)},\dots,{C_{1,k_1}(t)},
			\dots,
		{C_{m,1}(t)},\dots,{C_{m,k_m}(t)}\Big) , 
			\end{aligned}
		\end{eqnarray}
for all $i=1,2,\dots,{k_s},s=1,2,\dots,m,$		where  $ {\Phi_{s,i}}(\cdot) $ are the functions based on the expansion corresponding  to the  basis $ \mathbf{\hat{B}}. $
	\end{theorem}
	\begin{proof}
		Let us assume  that ${\mathbf{\hat{W}}_k}$  is an {invariant finite-dimensional product linear space} of the considered nonlinear vector differential operator
		$ \mathbf{\hat{A}}\big(\mathbf{U}\big).$ Then, we obtain the invariant criteria given in \eqref{Nd1.4}.
	From this invariant criteria \eqref{Nd1.4}, we get	   $ \mathbf{U}=(u_1,u_2,\dots,u_m)  $ with $u_s$   as given in \eqref{Nsoln}, and hence we obtain
		\begin{eqnarray}
			\begin{aligned}\label{prf1.2}
				{{{A}}_s}(\mathbf{U})=&{{{A}}_s}
				\Big( \sum\limits_{i=1}^{k_1}	{{{C}}_{1,i}(t)}{w_{1,i}}(x_1,\dots,x_N), 
				\dots,
				\sum\limits_{i=1}^{k_m}	{{{C}}_{m,i}(t)}{w_{m,i}}(x_1,\dots,x_N) \Big)
				\\	=&\sum\limits_{i=1}^{k_s}	{\Phi_{s,i}}\Big( {C_{1,1}(t)},\dots,{C_{1,k_1}(t)},
				\dots,
				 {C_{m,1}(t)},\dots,{C_{m,k_m}(t)}
				\Big){w_{s,i}}(x_1,\dots,x_N),
			\end{aligned}
		\end{eqnarray}
		where $ {\Phi_{s,i}}(\cdot),i=1,2,\dots,k_s,
		s=1,2,\dots,m,$ are the functions based on  expansion corresponding  to the basis $ \mathbf{\hat{B}}. $

		Now, we consider the given coupled NTFPDEs \eqref{Nd1.1} having the form
		\begin{eqnarray}\label{thm:eqn}	\dfrac{\partial^{\alpha_s} u_s}{\partial t^{\alpha_s}}={{{A}}_s}(\mathbf{U}),\alpha_s>0,s=1,2,\dots,m.\end{eqnarray}
		Then, by the virtue of the linear property of the Caputo fractional-order derivative and, for every $ \mathbf{U}=(u_1,u_2,\dots,u_m)  $  with $u_s$  as given in \eqref{Nsoln}, we obtain the left-hand side of  the given coupled NTFPDEs \eqref{thm:eqn}  as
		\begin{eqnarray}
			\begin{aligned}
				\label{prf1.1}
				\dfrac{\partial^{\alpha_s} u_s}{\partial t^{\alpha_s}}
				&	=\dfrac{\partial^{\alpha_s} }{\partial t^{\alpha_s}}\left(	\sum\limits_{i=1}^{k_s}	{{{C}}_{s,i}(t)}{w_{s,i}}(x_1,x_2,\dots,x_N) \right)
				\\
				&	=\sum\limits_{i=1}^{k_s}\left(\dfrac{\partial^{\alpha_s} {{{C}}_{s,i}(t)}}{\partial t^{\alpha_s} 	 }	 \right){w_{s,i}}(x_1,x_2,\dots,x_N),s=1,2,\dots,m,
			\end{aligned}
		\end{eqnarray}
		and the right-hand side of the  given coupled NTFPDEs  \eqref{thm:eqn} reduces  into the form  \eqref{prf1.2}.
		Thus,  substituting the obtained equations \eqref{prf1.2} and \eqref{prf1.1} into the given NTFPDEs  \eqref{thm:eqn},  we get
		$$
		\begin{aligned}
	&	\sum_{i=1}^{k_s}	\Big[ \dfrac{d^{\alpha_s}{C_{s,i}(t)}}{dt^{\alpha_s}} - {\Phi_{s,i}}( {C_{1,1}(t)},\dots,{C_{1,k_1}(t)},\dots,
		 {C_{m,1}(t)},\dots,{C_{m,k_m}(t)})\Big]  {w_{s,i}}=0 ,
\end{aligned}	$$
where ${w_{s,i}}={w_{s,i}}(x_1,\dots,x_N)$ for   $s=1,2,\dots,m.$		Since  the  functions ${w_{s,i}}\in B_s,$ $i=1,2,\dots,k_s$ are  linearly independent for $s=1,2,\dots,m,$ we find  that the  given coupled NTFPDEs \eqref{thm:eqn} reduces into the system of time-fractional ODEs as given in \eqref{Nsysfodes} whenever $ \mathbf{U}=(u_1,\dots,u_m)  $ with $u_s$  as given in \eqref{Nsoln}, which completes the proof.
	\end{proof}
	Note that  Theorem \ref{thm3} provides  a necessary condition for  {the existence of the generalized separable  {analytical} solution} \eqref{expected-solutionform} for the given coupled NTFPDEs \eqref{Nd1.1}.
	Thus, we can conclude that whenever  ${\mathbf{\hat{W}}_k}$  is an invariant  product linear space of the considered nonlinear vector differential operator
	$ \mathbf{\hat{A}}\big(\mathbf{U}\big)$,  then we can find the generalized separable {analytical} solution $ \mathbf{U}=(u_1,\dots,u_m)  $ with $u_s=u_s(x_1,\dots,x_N,t)$  given in \eqref{expected-solutionform}  for the  given  {multi-component} $(N+1)$-dimensional coupled  NTFPDEs \eqref{Nd1.1}.
	
	It is interesting to observe that we can derive the following $N$-different types of generalized separable  {analytical} solutions
	for the {multi-component} $(N+1)$-dimensional coupled  NTFPDEs \eqref{Nd1.1}:
	\begin{itemize}
		\item[$ \bullet $]
		The most generalized separable  {analytical} solution for the coupled NTFPDEs \eqref{Nd1.1} is obtained as $ \mathbf{U}=(u_1,u_2,\dots,u_m) $ with components
		\begin{eqnarray}
			\begin{aligned}\label{Ndsolnt1}
			&	u_s	= 
				\sum\limits_{r_1,\dots,r_N}	{C_{s,r_1,\dots,r_N}}(t) \ \prod_{i=1}^{N}{\xi^{i}_{s,r_i}}(x_i),
									\end{aligned}
		\end{eqnarray}
	where $	r_i=1,2,\dots,{k_{s,i}},$ $i=1,2,\dots,N,$ $ s=1,2,\dots,m.$	
	\newpage 
	\item[$ \bullet $]
		Other $ (N-1) $ types of generalized separable
		{analytical} solutions for the coupled NTFPDEs \eqref{Nd1.1} are obtained as
		$ \mathbf{U}=(u_1,u_2,\dots,u_m) $ with components
		\begin{eqnarray}
			\begin{aligned}\label{Ndsolntp}
					u_s	=&
				{C_{s,0}}(t)+\sum\limits_{i=1}^N \sum\limits_{r_i=1}^{k_{s,i}-1}
				C^i_{s,r_i}(t)\
				{\xi^{i}_{s,r_i}(x_i)}+\ldots
			\\&	+
				\sum\limits_{{{i_1}},\dots,{{i_{N-p}}}}	\left( \sum\limits_{{r_{i_1}},\dots,{r_{i_{N-p}}}}{C^{i_1\dots,{i_{N-p}}}_{s,{r_{i_1},\dots, r_{i_{N-p}}}}}(t)\prod\limits_{j=1}^{N-p}{\xi^{i_j}_{s,{r_{i_j}}}(x_{i_j})}\right)  ,
						\end{aligned}
		\end{eqnarray}
		where  	$r_{i_j}=1,2,\dots,{k_{s,i_j}-1},i_j=1,2,\dots,N,j=1,2,\dots,N-p,$ $i_1<i_2<\dots<i_{N-p},$ $
		s=1,2,\dots,m,$ $p=1,2,\dots,N-1.$
	\end{itemize}
	Note that the unknown functions   $ {\xi^{i}_{s,r_i}}(x_i),$   ${\xi^{i_j}_{s,{r_{i_j}}}(x_{i_j})}, $ $	{C_{s,0}}(t),$ $ C^i_{s,r_i}(t),$	 $\dots,$ ${C^{i_1\dots,{i_{N-p}}}_{s,{r_{i_1},\dots, r_{i_{N-p}}}}}(t),$  
	 and  ${C_{s,r_1,\dots,r_N}}(t) $   are to be determined  using  the considered coupled NTFPDEs \eqref{Nd1.1} and  corresponding invariant finite-dimensional product linear spaces.
	The details of finding the above
	{analytical} solutions \eqref{Ndsolnt1}-\eqref{Ndsolntp} using the various kinds  of invariant  product linear spaces {that} are discussed in the following subsection.
	
	\subsection{Determining different types of linear component spaces $ {{W}_{s,k_s}} $ with their invariant  product linear space $ {\mathbf{\hat{W}}_k} $ and
		{analytical} solutions for the  given  {multi-component}  NTFPDEs  \eqref{Nd1.1} }
	This subsection presents the details of  determining the  $N$-different types of invariant  product linear spaces ${\mathbf{\hat{W}}_k}$  and generalized separable
	{analytical} solutions for the  given  {multi-component} coupled NTFPDEs \eqref{Nd1.1}.
	Let us first explain the details of the systematic procedure to  determine  the different types of invariant  product linear  spaces ${\mathbf{\hat{W}}_k}$  for the vector differential operator $\mathbf{\hat{A}}\left(\mathbf{U} \right) .$
	Before proceeding further, we explain what kind of notations are used for linear spaces in the remaining sections, which are discussed in the following for easy understanding.
	\begin{itemize}
		\item[(I)] ${{W}_{s,k_s}^p}$ represents the  $s$-th component  type-$p$  linear component space with dimension $k_s$ for $s=1,2,\dots,m.$
		\item[(II)]$ {\mathbf{\hat{W}}_k^p} $  represents the  type-$p$  product linear vector space with dimension $k.$
	\end{itemize}
	Now, we aim to construct the  linear component spaces ${{W}_{s,k_s}^p}$  of the type-$p$   product linear space $ {\mathbf{\hat{W}}_k^p},p=1,2,\dots,N, $   {using the}
	 	system of ${k_{s,i}}$-th order linear  homogeneous ODEs:
	\begin{eqnarray}
		\begin{aligned}
		\label{odesNd}
		{\mathbf{L}^i_{s,k_{s,i}}}[y_{{s},i}(x_i)]\equiv &\dfrac{d^{k_{s,i}}y_{{s},i}(x_i)}{d x_i^{k_{s,i}}}+ {\varphi^i_{s,k_{s,i}-1}}(x_i) \dfrac{d^{k_{s,i}-1}y_{{s},i}(x_i)}{dx_i^{k_{s,i}-1}} +\dots
		\\&+{\varphi^i_{s,0}}(x_i)y_{{s},i}(x_i)=0,
	\end{aligned}\end{eqnarray}
	where  $ {k_{s,i}}\in \mathbb{N},$  and  $ {\varphi^i_{s,{j}}(x_i) ,j= 0,1,\dots,{k_{s,i}}-1}$ are continuous functions of $ x_i .$
	Here, we consider  the basis of the solution space (or the fundamental set of solutions) for given linear homogeneous ODEs \eqref{odesNd} as
	\begin{eqnarray}\label{fsNd}
		{H_{s,i}}=\left\{{\xi^i_{s,1}(x_i)},{\xi^i_{s,2}(x_i)},\ldots,{\xi^i_{s,{{k_{s,i}}}}}(x_i)\right\}, \ i=1,2,\dots,N,s=1,2,\dots,m.
	\end{eqnarray}
{	Next, we provide how to construct the  linear component spaces ${{W}_{s,k_s}^p}$ for $p=1,2,\dots,N, $ using the basis of the solution space ${H_{s,i}}$ given in  \eqref{fsNd}.
	Additionally, we explain the details of determining $ N $-different types of invariant  product linear spaces $ {\mathbf{\hat{W}}_k^p},p=1,2,\dots,N $ and the generalized separable
	{analytical} solutions given in \eqref{Ndsolnt1}-\eqref{Ndsolntp} for the given coupled NTFPDEs \eqref{Nd1.1} using  $ {\mathbf{\hat{W}}_k^p} .$ }
\subsubsection{Type-1 invariant  product linear space ${\mathbf{\hat{W}}_k^1}$ for the given vector differential operator  $ \mathbf{\hat{A}}\big(\mathbf{U}\big)$: }
	Let us \textbf{consider} the  linear component spaces  $ {{W}^1_{s,k_s}} $  of the type-1  product linear space  ${\mathbf{\hat{W}}^1_k}$
	as the linear span of the set
	\begin{equation}\label{basisNd}
		\mathit{B}^1_{s}=\Big\{\prod_{i=1}^N{{\xi^i_{s,r_i}}(x_i)}|r_i=1,\dots,{k_{s,i}},i=1,2,\dots,N \Big\},s=1,2,\dots,m,
	\end{equation}  where the functions ${{\xi^i_{s,r_i}}(x_i)}\in 	{H_{s,i}}$ is given in \eqref{fsNd}.
	Using the above set $\mathit{B}^1_{s}$, the type-1  linear component space  $ {{W}^1_{s,k_s}} $  can be defined as
	\begin{eqnarray}
		\begin{aligned}	\label{Ndt1c}
			{{W}^1_{s,k_s}}
		=&\text{Span}\left(\mathit{B}^1_{s}\right)
		\\	 =&
			\Big\{{w_s(x_1,\dots,x_N)}= \sum\limits_{r_1,\dots,r_N}{\mu_{s,r_1,\dots,r_N}}{\prod\limits_{i=1}^{N}}{{\xi^i_{s,r_i}}(x_i)}
			\Big{|} \mu_{s,r_1,\dots,r_N}\in  \mathbb{R},
		\\&	r_i =1,\dots,{k_{s,i}},i=1,\dots,N\Big\}.
		\end{aligned}
	\end{eqnarray}
	Thus,  the type-1
	 product linear space ${\mathbf{\hat{W}}^1_k }$ with linear component spaces $ {{W}^1_{s,k_s}},s=1,2,\dots,m, $ given in \eqref{Ndt1c} {can be written} as follows:
	\begin{eqnarray}
		\begin{aligned}\label{type-1}
			{\mathbf{\hat{W}}^1_k }=& {{W}^1_{1,k_1}}\times {{W}^1_{2,k_2}}\times {{W}^1_{3,k_3}}\times \cdots \times {{W}^1_{m,k_m}}
			\\
			=&\Big\{{\mathbf{w}(x_1,\dots,x_N)}=({w_1}(x_1,\dots,x_N),\dots,{w_m}(x_1,\dots,x_N))\,
			\big{|}\,
		\\&	 {w_s(x_1,\dots,x_N)}= \sum\limits_{r_1,\dots,r_N}{\mu_{s,r_1,\dots,r_N}}{\prod_{i=1}^{N}}{{\xi^i_{s,r_i}}(x_i)},
		 {{\xi^i_{s,r_i}}(x_i)}\in H_{s,i},
		 	\\&
		 	\mu_{s,r_1,\dots,r_N}\in  \mathbb{R},
			r_i =1,\dots,{k_{s,i}},i=1,\dots,N,  s=1,2,\dots,m
			\Big\}.
		\end{aligned}
	\end{eqnarray}
	Next, we show that the set $\mathit{B}^1_{s}$ is a basis for the  linear component space $ 	{{W}^1_{s,k_s}} ,$
	for all $ s=1,2,\dots,m. $
	\begin{proposition}
		Consider the set $\mathit{B}^1_{s}$  given  in \eqref{basisNd} with
		$ {{\xi^i_{s,r_i}}(x_i)}\in {H_{s,i}}$  for  $i=1, \dots,N,\text{ and }s=1,2,\dots,m, $ where
		${H_{s,i}}$  given in \eqref{fsNd}  is the basis of the solution spaces for the given  linear homogeneous ODEs \eqref{odesNd}.
		Then, the set $\mathit{B}^1_{s}$  forms  a  basis for the linear component space
		${{W}^1_{s,k_s}} =\text{Span}(\mathit{B}^1_{s})$
		with the dimension $k_s=\prod\limits_{i=1}^Nk_{s,i},$ for all $ s=1,2,\dots,m.$
		Additionally, the $k_s$-dimensional linear spaces $ {{W}^1_{s,k_s}} $  form the solution spaces for system of linear homogeneous ODEs \eqref{odesNd} for $s=1,2,\dots,m.$\label{proposition-basis}
	\end{proposition}
	\begin{proof}
		Assume  that ${{W}^1_{s,k_s}} =\text{Span}(\mathit{B}^1_{s}),$ where $\mathit{B}^1_{s}=\Big\{\prod\limits_{i=1}^N{{\xi^i_{s,r_i}}(x_i)}|r_i=1,\dots,{k_{s,i}},i=1, \dots,N \Big\},$
		and $ {{\xi^i_{s,r_i}}(x_i)}\in {H_{s,i}},$ $i=1, \dots,N,$ $s=1,2,\dots,m.$
		Now, we  show that ${{W}^1_{s,k_s}} =\text{Span}(\mathit{B}^1_{s})$ is a $k_s$-dimensional vector space with basis $\mathit{B}^1_{s}$
		under the vector addition ($+$)
		$$(w_{1} \ {+} \ w_2)(x_1,\dots,x_N)=
		w_{1} (x_1,\dots,x_N)+w_{2} (x_1,\dots,x_N),$$
		and the  scalar multiplication ($ {\cdot} $) as
		$$({\mu_s}\cdot{w})(x_1,\dots,x_N)={\mu_s}\cdot{w}(x_1,\dots,x_N),
	$$
	for every $ \mu{_s}\in \mathbb{R},{w}(x_1,\dots,x_N), w_{i} (x_1,\dots,x_N) \in{{W}^1_{s,k_s}},i=1,2 ,s=1,2,\dots,m.$
		It is sufficient to prove that $ {\mathit{B}^1_{s}}$
		is a linearly independent set. Let us consider the following linear combination of the elements in the set $\mathit{B}^1_{s}$ as
		\small	{$$
	\sum\limits_{r_1,\dots,r_N}{\mu_{s,r_1,\dots,r_N}} \prod\limits_{j=1}^{N}{\xi^j_{s,r_j}(x_j)}=	\sum\limits_{r_N=1}^{k_{s,N}}\Big( \sum\limits_{r_1,\dots,r_{N-1}}{\mu_{s,r_1,\dots,r_N}} \prod_{j=1}^{N-1}{\xi^j_{s,r_j}(x_j)}\Big){\xi^N_{s,r_N}(x_N)}=0,
		$$}
		where ${\xi^i_{s,r_i}(x_i)} \in {H_{s,i}}$ for $ x_i, {\mu_{s,r_1,\dots,r_N}}\in  \mathbb{R},$ $ r_i =1,2,\dots,{k_{s,i}},$ $k_{s,i}\in\mathbb{N},$   $i=1,2, \dots,N,$ and $ s=1,2,\dots,m. $
		Since the functions ${\xi^N_{s,r_N}(x_N)},r_N=1,2,\dots,k_{s,N},$ are linearly independent, we obtain
		$$\sum\limits_{r_1,\dots,r_{N-1}}
		{\mu_{s,r_1,\dots,r_N}} \prod\limits_{j=1}^{N-1}{\xi^j_{s,r_j}(x_j)}=0,$$ $r_i=1,2,\dots,k_{s,i},i=1,2,\dots,N-1.$
		Thus, {proceeding in a similar manner for $N-1$ successive steps, we get }
		$${\mu_{s,r_1,\dots,r_N}}=0,  \forall\,  r_i =1,2,\dots,{k_{s,i}},i=1,2,\dots,N,
		$$
		which gives that  ${\mathit{B}^1_{s}}$ is a linearly independent set.
		Hence, we conclude that $  {{W}^1_{s,k_s}} $ is a $ k_s$-dimensional linear space  over  $ \mathbb{R} $ with  basis ${\mathit{B}^1_{s}},$ where $k_s=\prod\limits_{i=1}^{N}k_{s,i},$ for each $ s=1,2,\dots,m.$
		
		Moreover, we observe that $ \prod\limits_{i=1}^{N}{\xi^i_{s,r_i}(x_i)}$  is a solution  {for the given} system of linear homogeneous ODEs   \eqref{odesNd}  for  $  r_i \in\{1,2,\dots,k_{s,i}\},i=1,2,\dots,N.$
		Since the basis set ${H_{s,i}}$ spans the solution space of the given linear homogeneous ODEs   \eqref{odesNd},
		we can easily observe that
		any solution of the system of linear homogeneous ODEs \eqref{odesNd} has the form
		$ {w_s(x_1,x_2,\dots,x_N)}=\sum\limits_{r_1,\dots,r_{N}}
		{\mu_{s,r_1,\dots,r_N}} \prod\limits_{j=1}^{N}{\xi^j_{s,r_j}(x_j)},$
		for any real arbitrary constants $ \mu_{s,r_1,\dots,r_N}, r_i =1,2,\dots,{k_{s,i}},i=1,2,\dots,N, $ and $s=1,2,\dots,m. $
		Thus, the  linear component spaces $ {{W}^1_{s,k_s}} $ with basis ${\mathit{B}^1_{s}} $ constitute the solution spaces of the given system of linear homogeneous ODEs \eqref{odesNd} for $s=1,2,\dots,m,$ which completes the proof.
	\end{proof}
	From Proposition \ref{proposition-basis} and Theorem \ref{axlerprop}, we note that the dimension of  type-1  product linear space ${\mathbf{\hat{W}}^1_k }$ given in \eqref{type-1} is
	$k=\sum\limits_{s=1}^m\prod\limits_{i=1}^{N}k_{s,i}.$
	Now, we say that the   type-1  product linear space ${\mathbf{\hat{W}}^1_k }$  given  in \eqref{type-1} is an invariant  product linear space for the given vector differential operator $ \mathbf{\hat{A}}\big(\mathbf{U}\big)$  if $\mathbf{U}\in {\mathbf{\hat{W}}^1_k }$ implies that $\mathbf{\hat{A}}\big(\mathbf{U}\big)\in{\mathbf{\hat{W}}^1_k } ,$ or equivalently ${{A}}_s(\mathbf{U})\in {{W}^1_{s,k_s}},s=1,2,\dots,m.$
	Thus, we can describe the invariant criteria for the  product linear space ${\mathbf{\hat{W}}^1_k }$ under the  vector differential operator $ \mathbf{\hat{A}}\big(\mathbf{U}\big)$  through the given system of linear homogeneous ODEs \eqref{odesNd} as follows:

	Additionally, from the above-invariant conditions,  the type-1 invariant  product linear space ${\mathbf{\hat{W}}^1_k }$  can be written as follows:
	$$
	\begin{aligned}
		{\mathbf{\hat{W}}^1_k }
		=&\Big\{({w_1},\dots,{w_m})\
		\big{|}\, 	{\mathbf{L}^i_{s,k_{s,i}}}[{w_s}]=0,
		w_s={w_s(x_1,\dots,x_N)}\in{{W}^1_{s,k_s}}, \\& s=1,2,\dots,m, i=1,2,\dots,N
		\Big\}
		.
	\end{aligned}
	$$
	By  Theorem \ref{thm3}, we observe that whenever  ${\mathbf{\hat{W}}^1_k }$   given  in \eqref{type-1} is a type-1  invariant  product linear space for the given vector differential operator $ \mathbf{\hat{A}}\big(\mathbf{U}\big),$ then given coupled NTFPDEs \eqref{Nd1.1} admits an
	{analytical} solution of the form \eqref{Ndsolnt1}.
	
	Next, we present the other types of the  invariant  product linear space ${\mathbf{\hat{W}}^{p+1}_k },p=1,2,\dots,N-1.$ Also, we discuss how to obtain the generalized separable  {analytical} solution given in \eqref{Ndsolntp} for the given coupled NTFPDEs \eqref{Nd1.1} when $p=1,2,\dots,N-1$ using the other types of invariant  product linear spaces.
	
	\subsubsection{Type-2 invariant  product linear space $ 	{\mathbf{\hat{W}}^2_k } $ of the given vector differential operator  $ \mathbf{\hat{A}}\big(\mathbf{U}\big)$:}
	Now, we consider the  linear component spaces $ {{W}^2_{s,k_s}} $   of the type-2 invariant  product linear space $ 	{\mathbf{\hat{W}}^2_k } $ as
	${{W}^2_{s,k_s}} =\text{Span}(\mathit{B}^2_{s}) ,$  $s=1,2,\dots,m,$ { where }
	\begin{eqnarray}
		\begin{aligned}
			\label{typ2-basis}
			\mathit{B}^2_{s}=&\Big\{1,{{\xi^i_{s,r_i}}(x_i)},
			\dots, {{\xi^{i_1}_{s,r_{i_1}}}(x_{i_1})}
			{{\xi^{i_2}_{s,r_{i_2}}}(x_{i_2})}
			\cdots
			{{\xi^{i_{N-1}}_{s,r_{i_{N-1}}}}(x_{i_{N-1}})}  \ |\\& r_{i_j},r_i=1,2,\dots,{k_{s,i}}, i_j,i=1,2,\dots,N,  j=1,2,\dots,N-1\Big\},
		\end{aligned}
	\end{eqnarray}
	$ {{\xi^i_{s,r_i}}(x_i)}\in {H_{s,i}}$  given in \eqref{fsNd} for  $r_i=1,2,\dots,k_{s,i},$ and  ${{\xi^i_{s,k_{s,i}}}(x_i)}=1,$ for all $i=1,2,\dots,N,s=1,2,\dots,m. $
	\newpage
	Thus, we can define  the type-2  product linear space $ 	{\mathbf{\hat{W}}^2_k } $ as follows:
	\begin{eqnarray}
	&	\begin{aligned}\label{type-2}
			{\mathbf{\hat{W}}^2_k }=& {{W}^2_{1,k_1}}\times {{W}^2_{2,k_2}}\times {{W}^2_{3,k_3}}\times \cdots \times {{W}^2_{m,k_m}}
			\\
			=&\Big\{{\mathbf{w}(x_1,\dots,x_N)}=({w_1},w_2,\dots,{w_m})
			\big{|}
			w_s=\mu_{s,0}+
			\sum\limits_{i=1}^N \sum\limits_{r_i=1}^{k_{s,i}-1}{\mu^i_{s,r_i}}{{\xi^i_{s,r_i}}(x_i)}+
	\\		&	\dots
	+
			\sum\limits_{{i_1},\dots,{i_{N-1}}
			} \left( 	\sum\limits_{r_{i_1},\dots,r_{i_{N-1}}}  {\mu^{{i_1},\dots,{i_{N-1}}}_{s,r_{i_1},\dots,r_{i_{N-1}}}}
	 {\prod_{j=1}^{N-1}}{{\xi^{i_j}_{s,r_{i_j}}}(x_{i_j})}\right) ,\mu_{s,0},\mu^i_{s,r_i},
			\\&
			\dots,{\mu^{{i_1},\dots,{i_{N-1}}}_{s,r_{i_1},\dots,r_{i_{N-1}}}}\in  \mathbb{R},	w_s={w_s(x_1,\dots,x_N)}\in{{W}^2_{s,k_s}}, {{\xi^n_{s,r_n}}(x_n)}\in H_{s,n},
			\\&
			{{\xi^n_{s,k_{s,n}}}(x_n)}=1,
					r_{n}=1,2,\dots,{k_{s,n}}-1,
			n=i,i_j=1,\dots,N,
		j=1,2,\dots,N-1,
		\\&i_1<i_2<\dots<i_{N-1},
			s=1,2,\dots,m
			\Big\}
			.
		\end{aligned}
	\end{eqnarray}
		Now, we prove that the set $\mathit{B}^2_{s}$ is a basis for the linear component space $ 	{{W}^2_{s,k_s}} ,$
	for all $ s=1,2,\dots,m, $ in the following proposition.
	\begin{proposition}
		Suppose that the   set $\mathit{B}^2_{s}$   given  in \eqref{typ2-basis} spans the linear  space
		${{W}^2_{s,k_s}},$ for all $s=1,2,\dots,m,$ when
		$ {{\xi^i_{s,r_i}}(x_i)}\in {H_{s,i}},r_i=1,\dots,k_{s,i},$
		and ${\xi^i_{s,k_{s,i}}(x_i)}=1,$ $i=1,2,\dots,N,$ in the set	${H_{s,i}}$  given  in  \eqref{fsNd}.
		Then, the set $\mathit{B}^2_{s}$  forms  a  basis for the linear space
		${{W}^2_{s,k_s}} $
		with the dimension $k_s=1+\sum\limits_{i=1}^N(k_{s,i}-1)+\dots+	\sum\limits_{i_1,\dots,i_{N-1}}\prod\limits_{j=1}^{N-1}(k_{s,i_j}-1),$ for all $ i_j=1,2,\dots,N,j=1,2,\dots,N-1,s=1,2,\dots,m.$
		Additionally, the $k_s$-dimensional linear spaces $ {{W}^2_{s,k_s}} $  form the solution spaces for system of linear homogeneous ODEs \eqref{odesNd} with ${\varphi^i_{s,0}}(x_i)=0,i=1,2,\dots,N,s=1,2,\dots,m.$\label{proposition-basis-type2}
	\end{proposition}
	\begin{proof}
		Assume that ${{W}^2_{s,k_s}}=\text{Span}(\mathit{B}^2_{s}),s=1,2,\dots,m.$
		We know that the set $\mathit{B}^2_{s}$ is a subset of the linearly independent set $\mathit{B}^1_{s},$ when ${\xi^i_{s,k_{s,i}}(x_i)}=1,$ $i=1,2,\dots,N,$ in the  set	${H_{s,i}}$    given  in \eqref{fsNd}.
		Since subset of any linearly independent set is a linearly independent set, we can say that $\mathit{B}^2_{s}$ is a basis  of the linear space ${{W}^2_{s,k_s}},s=1,2,\dots,m.$
		
		Additionally, following the above similar arguments as given in the proof of Proposition \ref{proposition-basis}, we conclude that  the  linear spaces $ {{W}^2_{s,k_s}} $  form the solution spaces for the system of linear homogeneous ODEs \eqref{odesNd} with ${\varphi^i_{s,0}}(x_i)=0,i=1,2,\dots,N,$ and the dimension of $ {{W}^2_{s,k_s}} $ is $k_s=1+\sum\limits_{i=1}^N(k_{s,i}-1)+\dots+	\sum\limits_{i_1,\dots,i_{N-1}}\prod\limits_{j=1}^{N-1}(k_{s,i_j}-1),$ $i_j=1,2,\dots,N,j=1,2,\dots,N-1,$ for all $ s=1,2,\dots,m.$
	\end{proof}
	Based on  Propositions \ref{axlerprop} and  \ref{proposition-basis-type2}, we find that the dimension of the type-2  product linear space $	{\mathbf{\hat{W}}^2_k }$ is
	$$k=\sum\limits_{s=1}^m\Big[ 1+\sum\limits_{i=1}^N(k_{s,i}-1)+\dots+	\sum\limits_{i_1,\dots,i_{N-1}}\prod\limits_{j=1}^{N-1}(k_{s,i_j}-1)\Big],$$
where $i_j=1,2,\dots,N,j=1,2,\dots,N-1 .$	Also, we known  that $\mathit{B}^2_{s}\subset \mathit{B}^1_{s}.$ Thus, we conclude  that the  type-2   product linear space $	{\mathbf{\hat{W}}^2_k }$ is a subspace of the    type-1   product linear space $	{\mathbf{\hat{W}}^1_k }$  whenever $1\in {{W}^{i}_{s,k_s}},i=1,2,s=1,2,\dots,m.$
	
	Now, we say that  the  product linear space $	{\mathbf{\hat{W}}^2_k }$ is a  type-2 invariant  product linear space of the given nonlinear vector differential operator
	$ \mathbf{\hat{A}}\big(\mathbf{U}\big) $ if $ \mathbf{\hat{A}}\left(  \mathbf{U}\right) \in{\mathbf{\hat{W}}^2_k } $ whenever 	$\mathbf{U} \in {\mathbf{\hat{W}}^2_k },$ 	or
	${{A}}_s\left[  {\mathbf{\hat{W}}^2_k }\right] \subseteq {{W}^{2}_{s,k_s}},s=1,2,\dots,m.$  That is,
	for each $ u_s\in {{W}^{2}_{s,k_s}}, s=1,2,\dots,m,$ we have ${{A}_s}( \mathbf{U})\in  {{W}^{2}_{s,k_s}},s=1,2,\dots,m .$
	For this case, the invariant criteria for the type-2 invariant  product linear space  $ 	{\mathbf{\hat{W}}^2_k } $ of the given  vector differential operator
	$ \mathbf{\hat{A}}\big(\mathbf{U}\big) $  read as, for  $ s=1,2,\dots,m,$
	\begin{eqnarray*}
	\begin{aligned}	&\dfrac{\partial^N {{{A}}_s}\big(\mathbf{U}\big)}{\partial x_1\partial x_2 \cdots\partial x_N}=0, \, {\mathbf{L}^i_{s,k_{s,i}}}\Big[ {{{A}}_s}\big(\mathbf{U}\big) \Big]=0,
		\text{ and }
		{\varphi^i_{s,0}}(x_i)=0, i=1,2,\dots,N
		,
		\\
		\text{ if }
		&\dfrac{\partial^N u_s}{\partial x_1\partial x_2 \cdots\partial x_N}=0, \, {\mathbf{L}^i_{s,k_{s,i}}}[u_s]=0, \text{ and } {\varphi^i_{s,0}}(x_i)=0, i=1,2,\dots,N.
\end{aligned}	\end{eqnarray*}
	It is interesting to note that based on the above-invariant criteria,   the type-2 invariant  product linear space  $ 	{\mathbf{\hat{W}}^2_k } $ can be read as
	$$ \begin{aligned}	{\mathbf{\hat{W}}^2_k }=&\Big\{{\mathbf{w}(x_1,\dots,x_N)}=({w_1}, \dots,{w_m})
		\big{|}	{\mathbf{L}^i_{s,k_{s,i}}}[{w_s}]=0,{\varphi^i_{s,0}}(x_i)=0,
		\dfrac{\partial^N w_s}{\partial x_1\cdots\partial{x_N}}=0,
		\\&w_s=	{w_s(x_1,x_2,\dots,x_N)}\in{{W}^2_{s,k_s}},
		i=1,2,\dots,N,  s=1,2,\dots,m
		\Big\} .\end{aligned}$$
	Additionally, we observe that  the given coupled NTFPDEs \eqref{Nd1.1} admit the  following form of the {analytical} solution
	\begin{eqnarray}\label{p=1}
		\begin{aligned}
				u_s	=&
			{C_{s,0}}(t)+\sum\limits_{i=1}^N \sum\limits_{r_i=1}^{k_{s,i}-1}
			C^i_{s,r_i}(t)\
			{\xi^{i}_{s,r_i}(x_i)}+\ldots
		\\&	 +
		\sum\limits_{{{i_1},\dots,{i_{N-1}}}}	\left( \sum\limits_{{r_{i_1},\dots,r_{i_{N-1}}}}{C^{i_1\dots,{i_{N-1}}}_{s,{r_{i_1},\dots, r_{i_{N-1}}}}}(t)\prod\limits_{j=1}^{N-1}{\xi^{i_j}_{s,{r_{i_j}}}(x_{i_j})}\right)  ,
				\end{aligned}
	\end{eqnarray}
where $	r_{i_j}=1,2,\dots,{k_{s,i_j}-1}, $ $i_j=1,2,\dots,N,$ $j=1,2,\dots,N-1,$ $ i_1<\dots<i_{N-1},$ $s=1,2,\dots,m,$	if the type-2  product linear space $	{\mathbf{\hat{W}}^2_k } $ is invariant under 	$ \mathbf{\hat{A}}\big(\mathbf{U}\big) .$
	The above solution \eqref{p=1} is exactly the same as the {analytical} solution \eqref{Ndsolntp}  when $ p=1 .$
	\begin{note}
		We conclude  that the  type-2  product linear space $	{\mathbf{\hat{W}}^2_k }$ is a subspace of the  type-1  product linear space $	{\mathbf{\hat{W}}^1_k },$  whenever $1\in {{W}^{1}_{s,k_s}},i=1,2,s=1,2,\dots,m,$
		which does not imply that  the nonlinear vector differential operator	$ \mathbf{\hat{A}}\big(\mathbf{U}\big) $  admits the type-2 invariant  product linear space $	{\mathbf{\hat{W}}^2_k },$ and the type-1 invariant  product linear space $	{\mathbf{\hat{W}}^2_k },$ simultaneously.
			In other words, the nonlinear vector differential operator	$ \mathbf{\hat{A}}\big(\mathbf{U}\big) $ may admit the type-2 invariant  product linear space $ 		{\mathbf{\hat{W}}^2_k } $, but the operator $ \mathbf{\hat{A}}\big(\mathbf{U}\big) $ need not  admit the type-1 invariant  product linear space $ 		{\mathbf{\hat{W}}^1_k }, $ and vice versa.
	\end{note}

	Next, we present the  type-$(p+1)  $ invariant  product linear spaces  $ 	{\mathbf{\hat{W}}^{p+1}_k }  $ for $ \mathbf{\hat{A}}\big(\mathbf{U}\big) $ along with their invariant criteria when $p=1,2,\dots,N-1$. Also, we discuss how to derive the corresponding generalized separable {analytical} solution for the given  {multi-component}   coupled NTFPDEs \eqref{Nd1.1} using $ 	{\mathbf{\hat{W}}^{p+1}_k } ,p=1,2,\dots,N-1. $
\subsubsection{Type-$ (p+1) $ invariant  product linear space $ 	{\mathbf{\hat{W}}^{p+1}_k } ,p=1,2,\dots,N-1$ for the given vector differential operator $ \mathbf{\hat{A}}\big(\mathbf{U}\big) $:}
	Let us consider the set $\mathit{B}^{p+1}_{s},$ for each   $s=1,2,\dots,m,$  as follows:
	\begin{eqnarray}
		\begin{aligned}	\label{typp-basis}
			\mathit{B}^{p+1}_{s}=&\Big\{1,{{\xi^i_{s,r_i}}(x_i)},
			\dots, {{\xi^{i_1}_{s,r_{i_1}}}(x_{i_1})}
			{{\xi^{i_2}_{s,r_{i_2}}}(x_{i_2})}
			\cdots
			{{\xi^{i_{N-p}}_{s,r_{i_{N-p}}}}(x_{i_{N-p}})}  \ |
			\\& r_{i_j},r_i=1,2,\dots,{k_{s,i}}, i_j,i=1,2,\dots,N,  j=1,2,\dots,N-p\Big\},
		\end{aligned}
	\end{eqnarray}
	where $ {{\xi^i_{s,r_i}}(x_i)}\in {H_{s,i}}$  given in \eqref{fsNd} for  $r_i=1,2,\dots,k_{s,i},$ and  ${{\xi^i_{s,k_{s,i}}}(x_i)}=1,$ for all $i=1,2,\dots,N,s=1,2,\dots,m. $
	Now, we define the  linear component spaces $ {{W}^{p+1}_{s,k_s}} $   of the type-$(p+1)$ invariant  product linear space $ 	{\mathbf{\hat{W}}^{p+1}_k } $ as
	${{W}^{p+1}_{s,k_s}} =\text{Span}(\mathit{B}^{p+1}_{s}),$  $s=1,2,\dots,m.$
	Thus, we construct the type-$(p+1)$  product linear space $ 	{\mathbf{\hat{W}}^{p+1}_k } $ as follows:
	\begin{eqnarray}
		\begin{aligned}\label{type-p}
			{\mathbf{\hat{W}}^{p+1}_k }=& {{W}^{p+1}_{1,k_1}}\times {{W}^{p+1}_{2,k_2}}\times {{W}^{p+1}_{3,k_3}}\times \cdots \times {{W}^{p+1}_{m,k_m}}
			\\
			=&\Big\{{\mathbf{w}(x_1,\dots,x_N)}=({w_1},\dots,{w_m})
			\big{|}
			w_s=\mu_{s,0}+
			\sum\limits_{i=1}^N \sum\limits_{r_i=1}^{k_{s,i}-1}{\mu^i_{s,r_i}}{{\xi^i_{s,r_i}}(x_i)}+
		\\&	\dots
			+
			\sum\limits_{{i_1},\dots,{i_{N-p}}
		} \left( 	\sum\limits_{r_{i_1},\dots,r_{i_{N-p}}}  {\mu^{{i_1},\dots,{i_{N-p}}}_{s,r_{i_1},\dots,r_{i_{N-p}}}}
	 {\prod_{j=1}^{N-p}}{{\xi^{i_j}_{s,r_{i_j}}}(x_{i_j})}\right) ,\,	\mu_{s,0},\mu^i_{s,r_i},\\&
	 \dots,
	 {\mu^{{i_1},\dots,{i_{N-p}}}_{s,r_{i_1},\dots,r_{i_{N-p}}}}\in  \mathbb{R},		
		w_s={w_s(x_1,\dots,x_N)}\in{{W}^{p+1}_{s,k_s}}, 
		\\&{{\xi^n_{s,r_n}}(x_n)}\in H_{s,n},	{{\xi^i_{s,k_{s,i}}}(x_i)}=1,		 r_{n}=1,2,\dots,{k_{s,n}}-1,\\&	n=	i,i_j=1,2,\dots,N,
			j=1,2,\dots,N-p,i_1<\dots <i_{N-p}	,
			s=1,2,\dots,m	
			\Big\},
				\end{aligned}
	\end{eqnarray}
for all $	p=1,2,\dots,N-1
.$	Next, we  show that the set $\mathit{B}^{p+1}_{s}$ is a basis for the linear component space $ 	{{W}^{p+1}_{s,k_s}} $
	for all $ s=1,2,\dots,m. $
	\begin{proposition}
		Suppose that the   linear space
		${{W}^{p+1}_{s,k_s}}$  is spanned by set $\mathit{B}^{p+1}_{s}$   given  in \eqref{typp-basis} for all $s=1,2,\dots,m,$ when
		$ {{\xi^i_{s,r_i}}(x_i)}\in {H_{s,i}},r_i=1,\dots,k_{s,i},$
		and ${\xi^i_{s,k_{s,i}}(x_i)}=1,$ $i=1,2,\dots,N,$ in the set	${H_{s,i}}$   given  in \eqref{fsNd}.
		Then, the set $\mathit{B}^{p+1}_{s}$  forms  a  basis for the linear space
		${{W}^{p+1}_{s,k_s}} $
		with the dimension $k_s=1+\sum\limits_{i=1}^N(k_{s,i}-1)+\dots+	\sum\limits_{i_1,\dots,i_{N-p}}\prod\limits_{j=1}^{N-p}(k_{s,i_j}-1),$ $i_j=1,2,\dots,N,j=1,2,\dots,N-p,$ for all $ s=1,2,\dots,m.$
		Additionally,   the $k_s$-dimensional linear spaces $ {{W}^{p+1}_{s,k_s}} $  form  solution spaces for  system of linear homogeneous ODEs \eqref{odesNd} with ${\varphi^i_{s,0}}(x_i)=0,i=1,2,\dots,N,s=1,2,\dots,m.$ \label{proposition-basis-typep}
	\end{proposition}
	\begin{proof}
		Similar  to the proof of Proposition \ref{proposition-basis-type2}.
	\end{proof}
	Based on  Propositions \ref{axlerprop} and  \ref{proposition-basis-typep}, we find that the dimension of the type-$({p+1})$  product linear space $	{\mathbf{\hat{W}}^{p+1}_k }$ is
	$$k=\sum\limits_{s=1}^m\Big[ 1+\sum\limits_{i=1}^N(k_{s,i}-1)+\dots+	\sum\limits_{i_1,\dots,i_{N-p}}\prod\limits_{j=1}^{N-p}(k_{s,i_j}-1) \Big],$$
where $i_j=1,2,\dots,N,j=1,2,\dots,N-p,p=1,2,\dots,N-1 .$	Additionally, we note that  $\mathit{B}^{p+1}_{s}\subset \mathit{B}^p_{s}$ whenever $1\in \mathit{B}^1_s,$ for $p=1,2,\dots,N-1,s=1,2,\dots,m.$ Thus, we say that the  type-$({p+1})$   product linear space $	{\mathbf{\hat{W}}^{p+1}_k }$ is a subspace of the    type-$p$  product linear space $	{\mathbf{\hat{W}}^p_k }$  whenever $1\in {{W}^{p}_{s,k_s}},p=1,2,\dots,N-1,s=1,2,\dots,m.$
	That is, when  $1\in {{W}^{i}_{s,k_s}},i=1,2,\dots,N,s=1,2,\dots,m,$
	we have
	$$	{\mathbf{\hat{W}}^N_k }\subset	{\mathbf{\hat{W}}^{N-1}_k }\subset\dots\subset	{\mathbf{\hat{W}}^2_k }	\subset{\mathbf{\hat{W}}^1_k }.
	$$
	A finite-dimensional   product linear space
	$	{\mathbf{\hat{W}}^{p+1}_k }$ is a  type-$({p+1})$ invariant  product linear space for the given vector differential operator
	$ \mathbf{\hat{A}}\big(\mathbf{U}\big) $ if $ \mathbf{\hat{A}}\left(  \mathbf{U}\right) \in{\mathbf{\hat{W}}^{p+1}_k } $	whenever 	$\mathbf{U} \in {\mathbf{\hat{W}}^{p+1}_k } $ or
	${{A}}_s\left[  {\mathbf{\hat{W}}^{p+1}_k }\right] \subseteq {{W}^{{p+1}}_{s,k_s}},s=1,2,\dots,m .$ That is,
	for each $ u_s\in {{W}^{{p+1}}_{s,k_s}} $ we have ${{A}_s}( \mathbf{U})\in  {{W}^{{p+1}}_{s,k_s}},p=1,2,\dots,N-1,s=1,2,\dots,m.$
	\\
	Thus, the invariant criteria for the type-$({p+1})$ invariant  product linear space  $ 	{\mathbf{\hat{W}}^{p+1}_k } $ of the given  vector differential operator
	$ \mathbf{\hat{A}}\big(\mathbf{U}\big) $  can be read as, for every $ s=1,2,\dots,m,$
	\begin{eqnarray*}
		\begin{aligned}
			&\dfrac{\partial^N {{{A}}_s}\big(\mathbf{U}\big)}{\partial x_{i_1}\partial x_{i_2}\cdots\partial x_{i_{N-p-1}}}=0, \, {\mathbf{L}^j_{s,k_{s,j}}}\Big[ {{{A}}_s}\big(\mathbf{U}\big) \Big]=0,
			\text{ and }
			{\varphi^j_{s,0}}(x_j)=0, 
		 \text{ whenever }
			\\&
			\dfrac{\partial^N {u_s}}{\partial x_{i_1}\partial x_{i_2}\cdots\partial x_{i_{N-p-1}}}=0, \, {\mathbf{L}^j_{s,k_{s,j}}}[ {u_s} ]=0,
			\text{ and }
			{\varphi^j_{s,0}}(x_j)=0, 
		\end{aligned}
	\end{eqnarray*}
where $i_r,j=1,2,\dots,N,r=1,2,\dots,N-(p+1).$	Also,  from the above-invariant criteria, the type-$({p+1})$ invariant  product linear space  $ 	{\mathbf{\hat{W}}^{p+1}_k } $ can be written as
	$$
	\begin{aligned}	{\mathbf{\hat{W}}^{p+1}_k }=&\Big\{{\mathbf{w}(x_1,\dots,x_N)}=({w_1},w_2,\dots,{w_m})
		\,	\big{|}	\,\dfrac{\partial^N {w_s}}{\partial x_{i_1}\partial x_{i_2}\cdots\partial x_{i_{N-p-1}}}=0, 
		\\&  {\mathbf{L}^j_{s,k_{s,j}}}[ {w_s} ]=0, {\varphi^j_{s,0}}(x_j)=0,
		 w_s=	{w_s(x_1,\dots,x_N)}\in{{W}^{p+1}_{s,k_s}},
	\\& 	i_r,j=1,2,\dots,N,  s=1,2,\dots,m,r=1,2,\dots,N-(p+1)
		\Big\} .\end{aligned}$$
	
	It is important to observe that the  given coupled NTFPDEs \eqref{Nd1.1} obtain the  {analytical} solution of the form given in  \eqref{Ndsolntp} whenever the given  vector differential operator
	$ \mathbf{\hat{A}}\big(\mathbf{U}\big) $ admits $ 	{\mathbf{\hat{W}}^{p+1}_k } ,p=1,2,\dots,N-1.$
	Additionally, {it is noteworthy to} mention that   the type-2  product linear space $	{\mathbf{\hat{W}}^2_k } $ can be obtained from the  type-$({p+1})$ invariant  product linear space  $ 	{\mathbf{\hat{W}}^{p+1}_k } $ when $ p=1 .$
	\begin{note}
		We observe that the relation is obtained as 	${\mathbf{\hat{W}}^N_k }\subset	{\mathbf{\hat{W}}^{N-1}_k } \subset \dots\subset	{\mathbf{\hat{W}}^2_k }	\subset{\mathbf{\hat{W}}^1_k }$  from all $N$-types of   product linear spaces ${\mathbf{\hat{W}}^p_k }$,  whenever $1\in {{W}^{p}_{s,k_s}},s=1,2,\dots,m,p=1,2,\dots,N.$ However, we can say that, in general, it need not be provided with
		the  product linear spaces   $	{\mathbf{\hat{W}}^1_k },$ $	{\mathbf{\hat{W}}^2_k },\dots,{\mathbf{\hat{W}}^N_k },$   are  invariant under the nonlinear vector differential operator	$ \mathbf{\hat{A}}\big(\mathbf{U}\big), $  simultaneously.
		Thus, it can be simply said that the nonlinear vector differential operator $ \mathbf{\hat{A}}\big(\mathbf{U}\big) $ may admit the type-${p_1}$ invariant  product linear space $ 		{\mathbf{\hat{W}}^{p_1}_k } $, but the operator $ \mathbf{\hat{A}}\big(\mathbf{U}\big) $ need not  admit the type-${p_2}$ invariant  product linear space $ 		{\mathbf{\hat{W}}^{p_2}_k } $ whenever ${p_1}\neq{p_2,} p_1,p_2=1,2,\dots,N.$
	\end{note}
	Next, we show the efficient application of the invariant subspace method by determining various dimensional type-1 and type-2 invariant  product linear spaces and their generalized separable  {analytical} solutions for the  {multi-component} $(2+1)$-dimensional coupled  nonlinear  {TFDCWEs} \eqref{cdsystemeqn} in the following section.
	\section{Invariant subspace method for the  {multi-component} $(2+1)$-dimensional coupled  nonlinear {TFDCWEs} $(7)$}\label{sec5}
	In this section, we explain how to apply the invariant subspace method to obtain the  {analytical} solutions of the  {multi-component} $(2+1)$-dimensional coupled  nonlinear {TFDCWEs} \eqref{cdsystemeqn}.
	So, let us first present how to find different type-1 and type-2 invariant  product linear spaces  for the following   vector differential operators of the given coupled equations \eqref{cdsystemeqn}:
	\begin{eqnarray}\label{cdsystem}
		\mathbf{F}\left(\mathbf{U}\right)=(	F_1(u_1,u_2),	F_2(u_1,u_2)),
	\end{eqnarray}
where  	\begin{eqnarray*}
		\begin{aligned}
	F_1(u_1,u_2)=&
			\sum\limits_{i=1}^2	
			{\frac {\partial }{\partial x_i}}
			\left(  \beta_{1i}(u_1,u_2) {\frac {\partial u_{{1}} }{	\partial x_i}}
			+ \rho_{1i} (u_1,u_2){\frac {\partial u_{{2}} }{
					\partial x_i}}  \right)
		\\& 	+\sum\limits_{i=1}^2 \left( \zeta_{1i}(u_1,u_2) {
				\frac {\partial u_{{1}} }{\partial x_i}}
			+ \upsilon_{1i} (u_1,u_2){\frac {\partial u_{{2}} }{\partial x_i}}\right),
			\\
	\text{and }		F_2(u_1,u_2)=&
			\sum\limits_{i=1}^2	
			{\frac {\partial }{\partial x_i}}
			\left(  \beta_{2i}(u_1,u_2) {\frac {\partial u_{{1}} }{	\partial x_i}}
			+ \rho_{2i} (u_1,u_2){\frac {\partial u_{{2}} }{
					\partial x_i}}  \right)
		\\&	+\sum\limits_{i=1}^2 \left( \zeta_{2i}(u_1,u_2) {
				\frac {\partial u_{{1}} }{\partial x_i}}
			+ \upsilon_{2i} (u_1,u_2){\frac {\partial u_{{2}} }{\partial x_i}}\right)
			.
		\end{aligned}
	\end{eqnarray*}
	In general, we cannot find the invariant   product linear spaces for  arbitrary  differential operators $F_s(u_1,u_2),s=1,2.$ Thus, we consider the   {multi-component} $(2+1)$-dimensional coupled  {TFDCWEs} \eqref{cdsystemeqn}  with cubic nonlinearity.
	So, let us consider the following forms of  diffusion coefficients and convection terms as
	$$
	\begin{aligned}
	&	\beta_{s1}(u_1,u_2)=k_{{s5}}  u_1  ^{2}+k_{{s4}}  u_{{2}}^{2}+k_{{s3}} u_{{1}} u_2		+k_{{s2}}u_{{2}} +k_{{s1}}u_{{1}} 		+k_{{s0}} ,
\\	&
		\rho_{s1}(u_1,u_2)= q_{{s5}}  u_{{1}}^{2}		+q_{{s4}}  u_{{2}}  ^{2}		+q_{{s3}}u_{{1}}u_{{2}} +q_{{s2}}u_{{2}}+q_{{s1}}u_{{1}}	+q_{{s0}} ,
		\\
	&	\beta_{s2}(u_1,u_2)= p_{{s5}}  u_{{1}}  ^{2}	+p_{{s4}}  u_{{2}}  ^{2}		+p_{{s3}}u_{{1}} u_{{2}}		+p_{{s2}}u_{{2}}+p_{{s1}}u_{{1}}		+p_{{s0}} ,
\\
	&
		\rho_{s2}(u_1,u_2)= c_{{s5}}  u_{{1}}^{2}		+c_{{s4}}  u_{{2}}^{2}		+c_{{s3}}u_{{1}} u_{{2}}		+c_{{s2}}u_{{2}} +c_{{s1}}u_{{1}} 		+c_{{s0}} ,
\\&	\zeta_{s1}(u_1,u_2)=\eta_{{s5}}u_{{1}} ^{2}		+\eta_{{s4}}  u_{{2}}^{2}		+\eta_{{s3}}u_{{1}}  u_{{2}} 		+\eta_{{s2}}u_{{2}} 		+\eta_{{s1}}u_{{1}} 		+\eta_{{s0}} ,
\end{aligned}
$$ 	$$
\begin{aligned}	&	\upsilon_{s1}(u_1,u_2)= g_{{s5}}  u_{{1}}^{2}
		+g_{{s4}}	 u_{{2}}^{2}
		+g_{{s3}}u_{{1}}u_{{2}}
		+g_{{s2}}u_{{2}}
		+g_{{s1}}u_{{1}}
		+g_{{s0}}
		,
		\\&	\zeta_{s2}(u_1,u_2)= d_{{s5}}  u_{{1}}^{2}
		+d_{{s4}}  u_{{2}}^{2}
		+d_{{s3}}u_{{1}} u_{{2}}
		+d_{{s2}}u_{{2}}
		+d_{{s1}}u_{{1}}
		+d_{{s0}},
\\	&
		\upsilon_{s2}(u_1,u_2)= h_{{s5}}  u_{{1}}^{2}
		+h_{{s4}}  u_{{2}}^{2}
		+h_{{s3}}u_{{1}} u_{{2}}
		+h_{{s2}}u_{{2}} +h_{{s1}}u_{{1}} +h_{{s0}}
		,
	\end{aligned}
	$$ where $k_{si},q_{si},p_{si},c_{si},\eta_{si},g_{si},d_{si},h_{si}\in\mathbb{R},$ for  $s=1,2,$ and $i=0,1,2,\dots,5.$
	{In the following subsection, we  present a detailed description of  finding two different types of  invariant   product linear spaces $ {{	\mathbf{\hat{W}}^{1}_k}} $ and $ {{	\mathbf{\hat{W}}^{2}_k}}$ for the  cubic nonlinear vector differential operator $  \mathbf{F}\left(\mathbf{U}\right). $  }
	\subsection{Finding the invariant  product linear  space ${{	\mathbf{\hat{W}}_k}}$ for the considered operator  $  \mathbf{F}\left(\mathbf{U}\right)  $ }
	Let us describe how to construct  two types of various finite-dimensional invariant  product linear spaces of the cubic nonlinear vector differential operator $ 	\mathbf{F}\left(\mathbf{U}\right) $ given in  \eqref{cdsystem}.
	Thus, we consider the  system of $k_{s,i}$-order linear homogeneous ODEs of the form
	\begin{eqnarray}
		\begin{aligned}
			\label{odesrdw}
				{\mathbf{L}^i_{s,k_{s,i}}}[y_{{s},i}(x_i)]\equiv& \dfrac{d^{k_{s,i}}y_{{s},i}(x_i)}{d x_i^{k_{s,i}}}+ {\varphi^i_{s,k_{s,i}-1}}(x_i) \dfrac{d^{k_{s,i}-1}y_{{s},i}(x_i)}{dx_i^{k_{s,i}-1}} +\ldots
				\\&+{\varphi^i_{s,0}}(x_i)y_{{s},i}(x_i)=0,k_{s,i}\in \mathbb{N},s,i=1,2
		\end{aligned}
	\end{eqnarray}
	with the basis of the solution spaces for  ${\mathbf{L}^i_{s,k_{s,i}}}[y_{{s},i}(x_i)]=0$ as follows:
	\begin{eqnarray}\label{fsrdw}
		{H_{s,i}}=\left\{{\xi^i_{s,1}(x_i)},{\xi^i_{s,2}(x_i)},\ldots,{\xi^i_{s,{{k_{s,i}}}}}(x_i)\right\}, \ s,i=1,2.
	\end{eqnarray}
	Note that the given cubic nonlinear vector differential operator $ 	\mathbf{F}\left(\mathbf{U}\right) $ involves two space variables $x_1$ and $x_2.$ Thus, we consider $N=2.$
	Hence, we can define two types of invariant  product linear spaces, which are named  type-1  product linear spaces $ {{	\mathbf{\hat{W}}^{1}_k}} $  and  type-2  product linear spaces $ {{	\mathbf{\hat{W}}^{2}_k}} .$
	Now, let us first explain how to find the type-1  product linear spaces $ {	{	\mathbf{\hat{W}}^{1}_k}} $  for  $ 	\mathbf{F}\left(\mathbf{U}\right) $ in the following.
	\\
	\\
	\textbf{Type-1 invariant  product linear spaces  $ {	\mathbf{\hat{W}}^{1}_k} $ preserved by the given  coupled equations \eqref{cdsystem}:}
	From the definition of type-1 invariant  product linear  spaces given in \eqref{type-1}, we can write the type-1 invariant  product linear spaces $ {	{	\mathbf{\hat{W}}^{1}_k}} $ for the given coupled equations  \eqref{cdsystem}  as follows:
	\begin{eqnarray}
		\begin{aligned}
			{	\mathbf{\hat{W}}^{1}_k} =&{{W^1_{1,k_1}}}\times {{W^1_{2,k_2}}}
			\\
			=&\Big\{w(x_1,x_2)=(w_1(x_1,x_2),w_2(x_1,x_2))|\ {\mathbf{L}^i_{s,k_{s,i}}}[w_s(x_1,x_2)]=0,s,i=1,2 \Big\}	
			\\
			=&	\Big\{(w_1,w_2)|\
			w_s:={w_s(x_1,x_2)}=\sum\limits_{{r_1,r_2}}{\mu^{s}_{r_1,r_2}}{\xi^{1}_{s,r_1}(x_1)}{\xi^2_{s,r_2}(x_2)}\in {{W^1_{s,k_s}}} ,
		\\&	{\mu^{s}_{i_1,i_2}}\in\mathbb{R},  r_i=1,\dots,k_{s,i} s,i=1,2 \Big\},
			\label{vm2d}
		\end{aligned}
	\end{eqnarray}
	where  ${\mathbf{L}^i_{s,k_{s,i}}}(\cdot)$ is
	given in \eqref{odesrdw} and ${\xi^i_{s,1}(x_i)}\in H_{s,i}$ is given in \eqref{fsrdw} for $s,i=1,2.$
	Note that the dimension of
	$ {	\mathbf{\hat{W}}^{1}_k} $ is
	$ {k}=$dim$({{	\mathbf{\hat{W}}^{1}_k}})=k_{1,1}k_{1,2}+k_{2,1}k_{2,2}.$
	
	The above  product linear space \eqref{vm2d} is invariant under $ \mathbf{F}\left(\mathbf{U}\right) $ if
	$ \mathbf{F} \left[ {{	\mathbf{\hat{W}}^{1}_k}} \right]\subseteq {{	\mathbf{\hat{W}}^{1}_k}} , $
	which means that,  for every
	$ \mathbf{U} \in {{	\mathbf{\hat{W}}^{1}_k}} ,$
	${F_s}({{	\mathbf{\hat{W}}^{1}_k}})= {F_s}({{W^1_{1,k_1}}}\times {{W^1_{2,k_2}}} ) \subseteq{{W^1_{s,k_s}}},\forall s=1,2. $
	These invariant conditions can be rewritten in the following equivalent form
	\begin{eqnarray}\label{t1invc}
		{\mathbf{L}^i_{s,k_{s,i}}}\Big[ {F_s}\big(u_1,u_2\big)\Big]=0 \text{ whenever } 	{\mathbf{L}^i_{s,k_{s,i}}}[u_s]=0 ,\, s,i=1,2.
	\end{eqnarray}
	Observe that if the type-1  product linear space  $ {	{	\mathbf{\hat{W}}^{1}_k}} $ is invariant under $ \mathbf{F}\left(\mathbf{U}\right)$, then  we expect the {analytical} solution for the given coupled equations \eqref{cdsystemeqn}  as
	$ 	\mathbf{U}=(u_1,u_2) $ with component functions
	\begin{eqnarray*}
		\begin{aligned}
			u_s	= \sum\limits_{r_1=1}^{k_{s,1}}\sum\limits_{r_2=1}^{k_{s,2}}	{C_{s,r_1,r_2}}(t) {\xi^{1}_{s,r_1}(x_1)}{\xi^2_{s,r_2}(x_2)},{k_{s,i}}\in\mathbb{N},s,i=1,2.
		\end{aligned}
	\end{eqnarray*}
		Now, we present a detailed illustration for finding the type-1  product linear spaces  $ {	{	\mathbf{\hat{W}}^{1}_k}} $ for the given coupled equations \eqref{cdsystem} using the given system of  ODEs \eqref{odesrdw} with ${k_{s,i}}=2.$ Thus, for ${k_{s,i}}=2,$ the  given system of  ODEs  \eqref{odesrdw}  is obtained  as follows:
	\begin{eqnarray}
		\begin{aligned}
			\label{cdsysodes}
			&{\mathbf{L}^i_{1,k_{1,i}}}[y_{1,i}(x_i)]\equiv \Big[\dfrac{d^{{2}}}{d x_i^{{2}}}+ a_{i1}\dfrac{d}{dx_i}+{a_{i0}}\Big]y_{1,i}(x_i)=0,
		\\	\text{ and } &	{\mathbf{L}^i_{2,k_{2,i}}}[y_{2,i}(x_i)]\equiv \Big[\dfrac{d^{{2}}}{d x_i^{{2}}}+ b_{i1}\dfrac{d}{dx_i}+{b_{i0}}\Big]y_{2,i}(x_i)=0,
		\end{aligned}
	\end{eqnarray}
	where $\varphi^i_{1,{j}}(x_i) =a_{ij},\varphi^i_{2,{j}}(x_i) =b_{ij} \in\mathbb{R},$ and the fundamental set of solutions of ${\mathbf{L}^i_{1,k_{1,i}}}[y_{1,i}(x_i)]=0$ is $ 	H_{s,i}=\left\{{\xi^i_{s,1}(x_i)},{\xi^i_{s,2}(x_i)}\right\}, $ $j=0,1, s,i=1,2.$
	Thus, the given invariant conditions \eqref{t1invc} for the type-1 invariant  product linear space  $ {{	\mathbf{\hat{W}}^{1}_k}} $ of $ \mathbf{F}\left(\mathbf{U}\right) $ are reduced as
\small{	\begin{eqnarray}
		&\Big( \dfrac{d^{{2}}}{d x_1^{{2}}}+ a_{11}\dfrac{d}{dx_1}+{a_{10}}\Big)  F_1\left(u_1,u_2\right) =0,\label{eq1}
		&\Big( \dfrac{d^{{2}}}{d x_2^{{2}}}+ a_{21}\dfrac{d}{dx_2}+{a_{20}}\Big)  F_1\left(u_1,u_2\right) =0,
		\\&
		\Big(\dfrac{d^{{2}}}{d x_1^{{2}}}+ b_{11}\dfrac{d}{dx_1}+{b_{10}}\Big) F_2\left(u_1,u_2\right) =0,\label{eq3}
		&\&\,
		\Big(\dfrac{d^{{2}}}{d x_2^{{2}}}+ b_{21}\dfrac{d}{dx_2}+{b_{20}}\Big) F_2\left(u_1,u_2\right) =0,
\\&
		\nonumber\text{whenever }\qquad\qquad\qquad\qquad\qquad&
		\\ \label{2,2t1u}&	\left( \dfrac{d^{{2}}}{d x_i^{{2}}}+ a_{i1}\dfrac{d}{dx_i}+{a_{i0}}\right)  u_1 =0,
		\text{ and }&	\left( \dfrac{d^{{2}}}{d x_i^{{2}}}+ b_{i1}\dfrac{d}{dx_i}+{b_{i0}}\right)u_2=0,i=1,2.
	\end{eqnarray}}
	Here, we consider  the cubic nonlinear vector differential operators $F_s(u_1,u_2),s=1,2,$ as follows:
	\begin{eqnarray}\label{cubicnonlinearcomp}
		\begin{aligned}
			F_s(u_1,u_2)=&
			{\frac {\partial }{\partial x_1}}
			\left[  \left( k_{{s5}}  u_1^{2}
			+k_{{s4}} u_{{2}}^{2}
			+k_{{s3}} u_{{1}} u_2
			+k_{{s2}}u_{{2}} +k_{{s1}}u_{{1}}
			+k_{{s0}} \right) {\frac {\partial u_{{1}} }{	\partial x_1}}
\right.\\ & \left.	+ \left( q_{{s5}}  u_{{1}}^{2}
			+q_{{s4}}  u_{{2}}^{2}
			+q_{{s3}}u_{{1}} u_{{2}}
			+q_{{s2}}u_{{2}}		 +q_{{s1}}u_{{1}}
			+q_{{s0}} \right) {\frac {\partial u_{{2}} }{
					\partial x_1}}  \right]
		\\&	+{\frac {\partial }{	\partial x_2}}
			\left[  \left( p_{{s5}}  u_{{1}}^{2}
			+p_{{s4}} u_{{2}}^{2}
			+p_{{s3}}u_{{1}} u_{{2}}
			+p_{{s2}}u_{{2}}+p_{{s1}}u_{{1}}
				+p_{{s0}} \right) {\frac {\partial u_{{1}}}{\partial x_2}}
		\right.\\& \left.	+ \left( c_{{s5}} u_{{1}}^{2}
			+c_{{s4}} u_{{2}}^{2}
			+c_{{s3}}u_{{1}} u_{{2}}
			+c_{{s2}}u_{{2}} 	+c_{{s1}}u_{{1}}
			+c_{{s0}} \right) {\frac {\partial u_{{2}}}{
					\partial x_2}}  \right]
	\end{aligned}
\end{eqnarray}
\begin{eqnarray}
\begin{aligned}\nonumber
&	+ \left( \eta_{{s5}} u_{{1}}^{2}
					+\eta_{{s4}} u_{{2}}^{2}
			+\eta_{{s3}}u_{{1}}  u_{{2}}
			+\eta_{{s2}}u_{{2}}
		+\eta_{{s1}}u_{{1}}
					+\eta_{{s0}} \right) {
				\frac {\partial u_{{1}} }{\partial x_1}}
		\\&	+ \left( g_{{s5}} u_{{1}}^{2}
			+g_{{s4}}	 u_{{2}}^{2}
			+g_{{s3}}u_{{1}}u_{{2}}
			+g_{{s2}}u_{{2}}
				+g_{{s1}}u_{{1}}
			+g_{{s0}}
			\right) {\frac {\partial u_{{2}}}{\partial x_1}}
		\\&	+
			\left( d_{{s5}}  u_{{1}}^{2}
			+d_{{s4}} u_{{2}}^{2}
		+d_{{s3}}u_{{1}} u_{{2}}
		 	+d_{{s2}}u_{{2}}
			+d_{{s1}}u_{{1}}
			+d_{{s0}}
			\right) {\frac {\partial u_{{1}}  }{\partial x_2}}
		\\&	+
			\left( h_{{s5}}u_{{1}}^{2}
				+h_{{s4}} u_{{2}}^{2}
			+h_{{s3}}u_{{1}} u_{{2}}
			+h_{{s2}}u_{{2}} +h_{{s1}}u_{{1}} +h_{{s0}}
			\right) {\frac {\partial u_{{2}} }{\partial x_2}}.
		\end{aligned}
	\end{eqnarray}
	Substituting the above operators  $F_s(u_1,u_2),s=1,2,$ 
	into the
	above invariant conditions  \eqref{eq1}-\eqref{eq3}, and equating each the coefficients of $u_s,$ 
	$\dfrac{\partial u_s}{\partial x_i},$
	$u_1\left( \dfrac{\partial u_s}{\partial x_i}\right), $
	$u_2 \left( \dfrac{\partial u_s}{\partial x_i}\right), $
	$\left( \dfrac{\partial u_s}{\partial x_i}\right)^2 ,s,i=1,2,$ etc.  to zero, we get the {subsequent} set of equations:
	\\
	(i)\quad From the first  equation
	of	\eqref{eq1}, we get
	$$
	\begin{aligned}
		&	c_{{11}}+p_{{12}}=0,	
		a_{{11}}b_{{10}}	(2k_{{12}}+q_{{11}})
		+b_{{10}}b_{{11}}(k_{{12}}
		+3q_{{11}})-b_{{10}}(\eta_{{12}}-2g_{{11}})=0,
		\\&	b_{{10}}(b_{{21}}c_{{11}}-h_{{11}})=0,
		b_{{10}}(a_{{21}}p_{{12}}-d_{{12}})=0,	
		d_{{12}}-a_{{21}}p_{{12}}=0,
		2h_{{11}}-b_{{21}}c_{{11}}=0,			
			\\&
	p_{1i}=p_{11}=
		d_{1i}=d_{11}=
		c_{1i}=c_{12}=
		h_{1i}=h_{12}=0
		,i=3,4,5,	(b_{{11}}-a_{{11}})(b_{{21}}c_{{11}}
		\\&-h_{{11}})=0,
						(b_{{10}}-a_{{10}})(b_{{21}}c_{{10}}
		-h_{{10}})=0,
		b_{{10}}(	(a_{{11}}-{b_{{11}}})b_{{11}}q_{{11}}
		+(3k_{{12}}+2q_{{11}})a_{{10}}
		\\&		-a_{{11}}g_{{11}}	+b_{{11}}g_{{11}}
		+a_{{20}}p_{{12}}
		{b_{{10}}}q_{{11}}
		+b_{{20}}c_{{11}})=0,
		(2b_{{11}}-a_{{11}})b_{{10}}q_{{10}}+(a_{{10}}
		-b_{{10}})g_{{10}}	
		\\&-a_{{10}}b_{{11}}q_{{10}}
+({b_{{11}}}g_{{10}}
			+b_{{20}}c_{{10}}
		-b_{{11}}^2q_{{10}})(b_{{11}}-a_{{11}})=0,
		{a_{{11}}}( k_{{12}}({a_{{11}}}+2b_{{11}})
		+b_{{11}}q_{{11}}
	\\&		-\eta_{{12}})
		+{b^2_{{11}}}(k_{{12}}+3q_{{11}})
	-(\eta_{{12}}+2g_{{11}})b_{{11}}
		-b_{{10}}(3q_{{11}}+k_{{12}})
		-2(	a_{{10}}k_{{12}}	+a_{{20}}p_{{12}}
		\\&+b_{{20}}c_{{11}})=0,
		(a_{{11}}-{b_{{11}}}){b_{{11}}}^{2}q_{{11}}
		+a_{{10}}k_{{12}}(a_{{11}}
		+3b_{{11}})
			+(2a_{{10}}q_{{11}}
		+a_{{20}}p_{{12}}
		+2b_{{10}}q_{{11}}
		\\
		&+{b_{{11}}}g_{{11}} +b_{{20}}c_{{11}})b_{{11}}
		-
		a_{{11}}(a_{{20}}p_{{12}}
		+b_{{10}}q_{{11}}+b_{{11}}g_{{11}}
		+b_{{20}}c_{{11}})
		-2a_{{10}}\eta_{{12}}-b_{{10}}g_{{11}}=0,
			\\&
		b_{{10}}(	a_{{11}}b_{{11}}q_{{10}}-{b^2_{{11}}}q_{{10}}-a_{{10}}q_{{10}} -a_{{11}}g_{{10}}+{b_{{10}}}q_{{10}}
		+b_{{11}}g_{{10}}+b_{{20}}c_{{10}})
		-a_{{10}}b_{{20}}c_{{10}}=0,
		\\&	
		(b_{{11}}-a_{{11}})(a_{{21}}p_{{12}}-d_{{12}})=0,	(a_{{11}}-b_{{11}})(p_{{12}}+c_{{11}})=0,		
		(b_{{11}}-a_{{11}})(b_{{21}}c_{{10}}-h_{{10}})=0.
	\end{aligned}
	$$
	(ii) \quad From the second equation of
	\eqref{eq1}, we get
	$$ 	\begin{aligned}
	&	{b^{2}_{{21}}}	(a_{{21}}-{b_{{21}}})c_{{11}}
	+
	(b_{{21}}
	-a_{{21}})a_{{10}}k_{{12}}
	+a_{{20}}p_{{12}}(a_{{21}}+3b_{{21}})
	+(2a_{{20}}c_{{11}}-a_{{21}}h_{{11}} \\&	+b_{{10}}q_{{11}}+2b_{{20}}c_{{11}} +{b_{{21}}}h_{{11}})
	-a_{{21}}(b_{{10}}q_{{11}}
	+b_{{20}}c_{{11}})
	-b_{{20}}h_{{11}}-2a_{{20}}d_{{12}}=0,\\ &
		q_{1i}=q_{12}=
		k_{1i}=k_{11}=
		\eta_{1i}=\eta_{11}=
		g_{1i}=g_{12}=0
		,i=3,4,5,(b_{{21}}-a_{{21}})(b_{{11}}q_{{10}}
		\\&
		-g_{{10}})=0,
			b_{{20}}b_{{21}}c_{{10}}(a_{{21}}-{b_{{21}}})
	+a_{{20}}(-b_{{10}}q_{{10}}-b_{{20}}c_{{10}})+b_{{20}}(
	-a_{{21}}h_{{10}}+b_{{10}}q_{{10}}
	\\&+{b_{{20}}}c_{{10}}
	+b_{{21}}h_{{10}})=0,
		b_{{20}}b_{{21}}c_{{11}}(a_{{21}}-{b_{{21}}})
	+b_{{20}}(a_{{10}}k_{{12}}+3a_{{20}}p_{{12}}+2a_{{20}}c_{{11}}-a_{{21}}h_{{11}}
	\\&+b_{{10}}q_{{11}}
	+{b_{{20}}}c_{{11}}+b_{{21}}h_{{11}})=0,
		b_{{20}}((2p_{{12}}+c_{{11}})a_{{21}}
		+b_{{21}}(p_{{12}}+3bc_{{11}})
		-d_{{12}}+2h_{{11}})=0,
	\\
&	(b_{{11}}q_{{11}}
	-g_{{11}})	b_{{20}}=0,(b_{{20}}c_{{10}}
	-b_{{21}}^2c_{{10}}+b_{{21}}h_{{10}})(b_{{21}}-a_{{21}})+b_{{10}}q_{{10}}(2b_{{21}}-a_{{21}})
		 +(a_{{20}}	\\& -b_{{20}})h_{{10}}
		 -a_{{20}}b_{{21}}c_{{10}}=0,
		b_{{20}}(k_{{12}}+q_{{11}})=0,		{a_{{21}}}p_{{12}}({a_{{21}}}+2b_{{21}})
	 	+({b_{{21}}}-b_{{20}}-2a_{{20}})p_{{12}}
	\\&	+b_{{21}}(a_{{21}}
		+3{b_{{21}}})c_{{11}}
	-2a_{{10}}k_{{12}}
				-d_{{12}}(b_{{21}}
	+a_{{21}})
			-2b_{{10}}q_{{11}}
		-3b_{{20}}c_{{11}}-2b_{{21}}h_{{11}}=0, 
		\\	&
		g_{{11}}-b_{{11}}q_{{11}}=0,(b_{{20}}-a_{{20}})(b_{{11}}q_{{10}}+g_{{10}})=0,
		b_{{20}}(a_{{11}}k_{{12}}-\eta_{{12}})=0,q_{{11}}+k_{{12}}=0,
		\\&
		(a_{{21}}-b_{{21}})(k_{{12}}+q_{{11}})=0,
		(b_{{21}}-a_{{21}})(b_{{11}}q_{{11}}-g_{{11}})=0,
		(\eta_{{12}}-a_{{11}}k_{{12}})(a_{{21}}-b_{{21}})=0.
	\end{aligned}$$
	(iii)\quad  From the first equation of
	\eqref{eq3}, we get
	$$
	\begin{aligned}
		&
		p_{2i}=p_{21}=
		d_{2i}=d_{21}=
		c_{2i}=c_{22}=
		h_{2i}=h_{22}=0
		,i=3,4,5,
	d_{{22}}-a_{{21}}p_{{22}}=0,	\\ &	(p_{{22}}+c_{{21}})(b_{{11}} 
		-a_{{11}})=0,	
	 (b_{{21}}c_{{21}}
			-h_{{21}})
		(a_{{11}}-b_{{11}})=0,
		a_{{10}}(b_{{21}}c_{{21}}-h_{{21}})=0,
	\end{aligned}$$
$$	\begin{aligned}& 	a_{{10}}(a_{{21}}p_{{22}}-d_{{22}})=0,		 (a_{{11}}-b
		_{{11}})(a_{{21}}p_{{20}}-d_{{20}})=0,
		(a_{{11}}-b_{{11}})(a_{{21}}p_
		{{22}}-d_{{22}})=0,
		\\
		&	(a_{{21}}p_{{20}}-d_{{20}})(a_{{10}}
			-b_{{10}})=0,
		a_{{10}}(a_{{11}}(3k_{{22}}+q_{{21}})+b_{{11}}(k_{{22}}+2q_{{21}})+g_{{21}}
	\\&	-2\eta_{{22}})=0,
	a_{{10}}(a_{{11}}b_{{11}}+{a_{{10}}}k_{{20}}
			+a_{{11}}\eta_{{20}}		-{a_{{11}}}^{2}k_{{10}} k_{{20}}
	 +a_{{20}}p_{{20}}-b_{{10}}k_{{20}}-b_{{11}}\eta_{{20}})\\&-a_{{20}}b_{{10}}p_{{20}}=0,
(a_{{10}}{a_{{11}}}k_{{22}}-\eta_{{22}}a_{{10}})(b_{{11}}-{a_{{11}}}) +({a_{{10}}}k_{{22}} +a_{{20}}p_{{22}}
				+b_{{10}}
		+b_{{20}}c_{{21}})a_{{10}}
		\\&	+a_{{20}}(2k_{{22}}
	+3q_{{21}})=0,
			(b_{{11}}-{a_{{11}}})({a_{{11}}}^{2}k_{{20}}-a_{{11}}\eta_{{20}}
				-a_{{20}}p_{{20}})
		+(2a_{{11}}-b_{{11}})a_{{10}}k_{{20}}
			\\&
			+b_{{10}}(\eta_{{20}}
		-a_{{11}}k_{{20}})
		-a_{{10}}\eta_{{20}}=0,
		{a_{{11}}}^{2}(3k_{{22}}+q_{{21}})
		+
		a_{{11}}b_{{11}}(k_{{22}}
				+2q_{{21}})
		+({b_{{11}}}^{2}-a_{{10}}
		\\&	-2b_{{10}})q_{{21}}
	-3a_{{10}}k_{{22}}
		-a_{{11}}(2\eta_{{22}}+g_{{21}})
		-2a_{{20}}p_{{22}} -b_{{11}}g_{{21}}
			-2b_{{20}}c_{{21}}=0, p_{{22}}+c_{{21}}=0.		 	\end{aligned}\label{o2e3cc}
	$$
	(iv) \quad From the second equation
	of	\eqref{eq3}, we get
	$$
	\begin{aligned}
&		 (b_{{21}}-{a_{{21}}}p_{{22}}){a_{{21}}}^{2}p_{{22}}
+(a_{{10}}k_{{22}}			+b_{{10}}q_{{21}})(a_{{21}}-b_{{21}})
-b_{{20}}(2h_{{21}}+c_{{21}})
-b_{{21}}\\&		+p_{{22}}(2a_{{21}}
-b_{{21}})a_{{20}} +a_{{21}}({a_{{21}}}d_{{22}}
+2b_{{20}}p_{{22}}+3b_{{20}}c_{{21}}
-b_{{21}}d_{{22}}) 	-a_{{20}}d_{{22}} =0
,	\\	&
		q_{2i}=q_{22}=
		k_{2i}=k_{21}=
		\eta_{2i}=\eta_{21}=
		g_{2i}=g_{22}=0
		,i=3,4,5,	\eta_{{22}}-a_{{11}}k_{{22}}=0,	\\
		&(a_{{11}}k_{{20}}-\eta_{{20}})(a_{{2i}}-b_{{2i}})=0,
		i=0,1,
		2g_{{21}}	-2b_{{11}}q_{{21}}=0,
		a_{{20}}b_{{11}}q_{{21}}-a_{{20}}g_{{21}}=0,
	\\&	(a_{{11}}k_{{22}}-\eta_{{22}})(a_{{21}}-b_{{21}})=0,
		(a_{{21}}-b_{{21}})(b_{{11}}q_{{21}}-g_{{21}})=0,
				a_{{20}}a_{{21}} 	(3p_{{22}}
		+c_{{21}})
		\\&+
		a_{{20}}b_{{21}}(p_{{22}}+2c_{{21}}) -a_{{20}}(2d_{{22}}+h_{{21}})=0,	
		(b_{{21}}-a_{{21}})(k_{{22}}+q_{{21}})=0,
		a_{{20}}(k_{{22}}
		\\&
		+q_{{21}})=0,
		a_{{20}}(a_{{21}}b_{{21}}p_{{20}}
		-{a_{{21}}}^{2}p_{{20}}  +a_{{10}}k_{{20}}+{a_{{20}}}p_{{20}}+a_{{21}}d_{{20}} -b_{{20}}p_{{20}}
		-b_{{21}}d_{{20}}) \\	&
		-a_{{10}}b_{{20}}k_{{20}}
		=0,
		a_{{20}}[	(b_{{21}}-{a_{{21}}})a_{{21}}p_{{22}} +a_{{10}}k_{{22}}+{a_{{20}}}p_{{22}}
		 +a_{{21}}d_{{22}}
		+b_{{10}}q_{{21}}  -b_{{21}}d_{{22}}
		\\&+2b_{{20}}p_{{22}}
		+3b_{{20}}c_{{21}}]=0,
		({a_{{21}}}^{2}p_{{20}}-a_{{10}}k_{{20}}-a_{{21}}d_{{20}})(b_{{21}}-{a_{{21}}})  +a_{{20}}p_{{20}}(2a_{{21}}
		 -b_{{21}})\\&
			+b_{{20}}(d_{{20}}-a_{{21}}p_{{20}}) 	-a_{{20}}d_{{20}} =0,
		{a_{{21}}}^{2}(3p_{{22}}+c_{{21}}) +a_{{21}}b_{{21}}(p_{{22}}+2c_{{21}}) +({b_{{21}}}^{2}-a_{{20}} \\&
		-2b_{{20}})c_{{21}}
		-2a_{{10}}k_{{22}}
		-3a_{{20}}p_{{22}}
				-b_{{21}}h_{{21}}-2b_{{10}}
		q_{{21}}
		-a_{{21}}(2d_{{22}} +h_{{21}})
		=0,q_{{21}}+k_{{22}}=0.
		\end{aligned}\label{o2e4cc}
	$$
	Solving the above sets of equations  (i)-(iv) together, we obtain the eight-dimensional type-1 invariant  product linear spaces $ {\mathbf{\hat{W}}^{1}_8} $ with their corresponding {multi-component} cubic nonlinear vector differential operators $\mathbf{F}\left( \mathbf{U}\right) .$
	However, we cannot generally solve them for the above set of equations (i)-(iv).
	Thus, let us first consider the  case $ a_{11}=b_{21}=0 .$
	{For this case, the obtained eight-dimensional type-1 invariant  product linear spaces $ {\mathbf{\hat{W}}^{1}_8} $ and their  corresponding non-zero parameter restrictions in the differential operators $  F_s\left(u_1,u_2\right) ,s=1,2,$ are listed from
		Table 1 to Table 5.
		Additionally, we wish to point out that  some parameters of the  given component differential operators $  F_s\left(u_1,u_2\right),s=1,2, $ \eqref{cubicnonlinearcomp} are zero that are not listed in Table 1-Table 5 due to the  large number of parameters and cases.
		Hence, only the  non-zero parameters of   $  F_s\left(u_1,u_2\right),s=1,2, $  given in \eqref{cubicnonlinearcomp}
		are  listed  in Table 1-Table 5, and  the parameters, which are zero, are not shown explicitly.}

	Similarly,  we can compute other type-1 invariant  product linear spaces $ {{	\mathbf{\hat{W}}^{1}_8}} $  with their vector differential operators $\mathbf{F}\left( \mathbf{U}\right) $ for other cases such as $ a_{10}=b_{20}=0,$ $ a_{20}=b_{21}=0,$ etc. Moreover, in the same way, we can determine the other $k$-dimensional type-1 invariant  product linear spaces $ {{	\mathbf{\hat{W}}^{1}_k}} $ for the vector differential operators $\mathbf{F}\left( \mathbf{U}\right). $
	\begin{remark} It is important to mention that using the definition of type-2 invariant  product linear spaces  given in \eqref{type-2}, we can consider the type-2 invariant  product linear spaces $ {{	\mathbf{\hat{W}}^{2}_k}} $ for the given coupled equations \eqref{cdsystem}  as follows:
		\begin{eqnarray}
			\begin{aligned}
				&	{{	\mathbf{\hat{W}}^{2}_k}} = {{W}^{2}_{1,k_1}}\times {{W}^{2}_{2,k_2}},
				\label{t22d}
			\end{aligned}
		\end{eqnarray}
		where the component spaces $ {{W}^{2}_{s,k_s}},s=1,2, $ are given by
	\begin{eqnarray*}
		\small		\begin{aligned}
				{{W}^{2}_{s,k_s}}=&\Big\{{w_s(x_1,x_2)}|\ {\mathbf{L}^i_{s,k_{s,i}}}[{w_s(x_1,x_2)}]=0, {\varphi^i_{s,0}}(x_i)=0,
				\dfrac{\partial^2 {w_s(x_1,x_2)}}{\partial x_1\partial x_2 }=0,i=1,2 \Big\}	\\
				=&\Big\{{w_s(x_1,x_2)}= {\mu_{s,0}}+\sum\limits_{i,r_i}{\mu^{i}_{s,r_i}}{\xi^i_{s,{{r_{i}}}}}(x_i)	\Big{|} {\mu_{s,0}},{\mu^{i}_{s,r_i}}\in  \mathbb{R},
				r_i =1,\dots,k_{s,i}-1,i=1,2\Big\}.
			\end{aligned}
		\end{eqnarray*}
		Note that the dimension of $ 	{{	\mathbf{\hat{W}}^{2}_k}} $ is $ k= $dim$( {{	\mathbf{\hat{W}}^{2}_k}})=
		\sum\limits_{s=1}^{2}\left( {k_{s,1}}+{k_{s,2}}-1 \right).$
		The above  product linear space $ 	{{	\mathbf{\hat{W}}^{2}_k}} $ is invariant under $\mathbf{F}\left(\mathbf{U}\right)$ if
		${F_s}({{	\mathbf{\hat{W}}^{2}_k}})= {F_s}({{W}^{2}_{1,k_1}}\times {{W}^{2}_{2,k_2}} ) \subseteq{{W}^{2}_{s,k_s}},\forall s=1,2, $ whenever $ \mathbf{U}=(u_1,u_2) \in {{	\mathbf{\hat{W}}^{2}_k}} ,$
		which reduce to the following forms of invariant conditions:
		\begin{eqnarray}
			\begin{aligned}
				\label{t2invc}
				&{\mathbf{L}^i_{s,k_{s,i}}}[u_s]=0,{ \varphi^i_{s,0}}(x_i)=0,\dfrac{\partial^2 u_s}{\partial x_1\partial x_2}=0,
				\text{  implies that }
			\\&	{\mathbf{L}^i_{s,k_{s,i}}}\big[ F_s\left(u_1,u_2\right)\big]=0, { \varphi^i_{s,0}}(x_r)=0,\dfrac{\partial^2 F_s\left(u_1,u_2\right)}{\partial x_1\partial x_2 }=0, s,i=1,2.
			\end{aligned}
		\end{eqnarray}
		Observe that for $k_{s,i}=2,\varphi^i_{1,{j}}(x_i) =a_{ij}\in\mathbb{R},$  and  $\varphi^i_{2,{j}}(x_i) =b_{ij}\in\mathbb{R},$  $j=0,1,s,i=1,2,$ the above invariant conditions \eqref{t2invc} read as
		\begin{eqnarray}
			\begin{aligned}
				\label{2invct2}
				& \left( \dfrac{d^{{2}}}{d x_r^{{2}}}+ a_{r1}\dfrac{d}{dx_r}\right)  F_1\left(u_1,u_2\right) =0,
				\dfrac{\partial^2 F_1\left(u_1,u_2\right)}{\partial x_1\partial x_2}=0,
			\\&
				\left(\dfrac{d^{{2}}}{d x_r^{{2}}}+ b_{r1}\dfrac{d}{dx_r}\right)F_2\left(u_1,u_2\right)=0,
				\dfrac{\partial^2 F_2\left(u_1,u_2\right)}{\partial x_1\partial x_2}=0,
				\label{eq2t2}
				\\&
				\text{whenever }
				\left( \dfrac{d^{{2}}}{d x_r^{{2}}}+ a_{r1}\dfrac{d}{dx_r}\right)  u_1 =0,
				\dfrac{\partial^2 u_1}{\partial x_1\partial x_2}=0,
			\\&	\text{ and }	\left( \dfrac{d^{{2}}}{d x_r^{{2}}}+ b_{r1}\dfrac{d}{dx_r}\right)u_2=0,\dfrac{\partial^2 u_2}{\partial x_1\partial x_2}=0,r=1,2,
			\end{aligned}
		\end{eqnarray}
		where the differential operators $F_s\left(u_1,u_2\right),s=1,2, $ are considered as given in \eqref{cubicnonlinearcomp}.
		Proceeding in a way similar to  finding type-1 invariant  product linear spaces earlier, we can obtain the various six-dimensional type-2 invariant  product linear spaces $ {\mathbf{\hat{W}}^{2}_6} $ and their corresponding cubic nonlinear vector differential operator $  \mathbf{F}\left(\mathbf{U}\right) $  using the given invariant conditions \eqref{2invct2}.
		Also, we note that the type-2 invariant  product linear space $ {\mathbf{\hat{W}}^{2}_6} $ is a subspace of the type-1 invariant  product linear space  $ {\mathbf{\hat{W}}^{1}_8} $ when $a_{r0}=b_{r0}=0,r=1,2.$ That means  $ {{\mathbf{\hat{W}}^{2}_6}}\subset {{\mathbf{\hat{W}}^{1}_8}}, $ if $a_{r0}=b_{r0}=0,r=1,2.$   However, the type-2   product linear  space $ {\mathbf{\hat{W}}^{2}_6}$ is {preserved} under the nonlinear vector differential operator $  \mathbf{F}\left(\mathbf{U}\right) $ need not imply the same operator $  \mathbf{F}\left(\mathbf{U}\right) $ admits the corresponding type-1 invariant  product linear space  $ {{\mathbf{\hat{W}}^{1}_8}} ,$ in general.
	$$
 {\mathbf{L}^{i}_{s,k_{s,i}}}[u_s]=0  , s=1,2,\dots,m, i=1,2,\dots,N,
	$$ $$	\text{ implies that }  {\mathbf{L}^i_{s,k_{s,i}}} [{{A}}_s(\mathbf{U})  ] =0,s=1,2,\dots,m, i=1,2,\dots,N.
 	$$		
		
		Additionally,  we can find the other possible type-2 invariant  product linear spaces $ {{	\mathbf{\hat{W}}^{2}_k}} $ and their corresponding cubic nonlinear vector differential operator $  \mathbf{F}\left(\mathbf{U}\right) $ similarly  using suitable values of $n_r^s,r,s=1,2$ in the given invariant conditions \eqref{t2invc}.
		Moreover, we observe that  the type-2 invariant  product linear spaces $ {{	\mathbf{\hat{W}}^{2}_k}} $ are more helpful in finding  generalized separable {analytical} solutions of the given coupled equations \eqref{cdsystemeqn}, especially when it is not possible to construct a type-1 invariant  product linear space $ {{	\mathbf{\hat{W}}^{1}_k}} $ for the given coupled equations  \eqref{cdsystemeqn}.
	\end{remark}
		\begin{table}[pt]
\small	\caption{Type-1  product linear spaces $ {{	\mathbf{\hat{W}}^{1}_k}} $ and  their corresponding parameter restrictions of  $  \mathbf{F}\left(\mathbf{U}\right)  $ with $ a_{11}=b_{21}=0$ and $d_{10},\eta_{10},k_{10},p_{10},g_{20},h_{20},q_{20},c_{20}\in\mathbb{R}.$}
	{\begin{tabular}{@{}llllc@{}} \toprule
			Cases &	Some non-zero parameter 
			  &  Type-1 invariant product linear space 
		\\	&	
		restrictions of
		 $  \mathbf{F}\left(\mathbf{U}\right)  $ &  
		 $ {{	\mathbf{\hat{W}}^{1}_8}}={{W^1_{1,4}}}\times {{W^1_{2,4}}} $ \\
		 \midrule
			1. &
			$\left\{\begin{array}{ll}
				 p_{12} = -c_{11},
				p_{22} = -c_{21} ,
			\\	c_{21},
				c_{11},
				c_{10}\in\mathbb{R}
			\end{array}\right.$
			&
			$\left\{\begin{array}{ll}
				{W}^1_{1,4}= \text{Span}\{
				\sin(\sqrt{\theta}x_1),
				\cos(\sqrt{\theta}x_1),
				x_2\sin(\sqrt{\theta}x_1),	
			\\	x_2\cos(\sqrt{\theta}x_1)\} ,\theta_1={b_{10}-b_{11}^2}  \\
				{W}^1_{2,4}=
				\text{Span}\{
				e^{-\theta_2x_1},
				e^{-\theta_3x_1},
		 	x_2e^{-\theta_2x_1},
				x_2e^{-\theta_3x_1}\},
				\\
				\theta_2=\frac{b_{11}+\sqrt{b_{11}^2-4b_{10}}}{2},\theta_3=\frac{b_{11}-\sqrt{b_{11}^2-4b_{10}}}{2}
			\end{array}\right.
			$
			\\
						2. &
			$\left\{ \begin{array}{ll}
			 p_{12} = -c_{11},
				p_{22} = -c_{21} ,
			\\	g_{10},
				q_{10},
				c_{21},
				c_{11}\in\mathbb{R}
			\end{array}\right.$
			&
			$\left\{\begin{array}{ll}
				{W}^1_{1,4}=  \text{Span}\{
				e^{b_{11}x_1},	
				e^{-b_{11}x_1},
				x_2e^{b_{11}x_1},	
				x_2e^{-b_{11}x_1}\}
				\\	{W}^1_{2,4}=  \text{Span}\{
				\sin(\sqrt{b_{20}}x_2),
							e^{-b_{11}x_1}\sin(\sqrt{b_{20}}x_2),
			\\\cos(\sqrt{b_{20}}x_2),	e^{-b_{11}x_1}\cos(\sqrt{b_{20}}x_2)\}
			\end{array}\right.
			$
			\\
						3.	&
			$\left\{\begin{array}{ll}
				 p_{12} = -c_{11},
				p_{22} = -c_{21} ,
			\\	q_{10},
				c_{21},
				c_{11},
				c_{10}\in\mathbb{R}
			\end{array}\right.$
			&
			$\left\{\begin{array}{ll}
				{W}^1_{1,4}=  \text{Span}\{
				\sin(\sqrt{a_{10}}x_1),
				\cos(\sqrt{a_{10}}x_1),
			\\	x_2\sin(\sqrt{a_{10}}x_1),
				x_2\cos(\sqrt{a_{10}}x_1)\}        		\\
				{W}^1_{2,4}=  \text{Span}\{1,x_1,x_2,x_1x_2\}
			\end{array}\right.
			$
			\\
						4.&
			$\left\{\begin{array}{ll} p_{12} = -c_{11},
			p_{22} = -c_{21},
		\\	c_{11},c_{21}\in\mathbb{R}	\end{array}\right.$
			&
			$\left\{\begin{array}{ll}
				{W}^1_{1,4}= \text{Span}\{
				\sin(\sqrt{a_{10}}x_1)\sin(\sqrt{b_{20}}x_2),
		\\		\cos(\sqrt{a_{10}}x_1)\sin(\sqrt{b_{20}}x_2),
					\cos(\sqrt{b_{20}}x_2)\sin(\sqrt{a_{10}}x_1),
			\\	\cos(\sqrt{b_{20}}x_2)\cos(\sqrt{a_{10}}x_1)\}
				\\
				{W}^1_{2,4}=
				\text{Span}\{
				\sin(\sqrt{b_{20}}x_2),
				\cos(\sqrt{b_{20}}x_2),
				x_1\sin(\sqrt{b_{20}}x_2),
			\\
				x_1\cos(\sqrt{b_{20}}x_2)\} \end{array}\right.$
						\\
								5.&
		$\left\{\begin{array}{ll}p_{12} = -c_{11},  p_{22} = -c_{21},
			\\	p_{20},c_{21},c_{11},c_{10}\in\mathbb{R}
		\end{array}\right.$
			&
			$\left\{\begin{array}{ll}
				{W}^1_{1,4}= \text{Span}\{\sin(\sqrt{a_{10}}x_1),\cos(\sqrt{a_{10}}x_1),x_2\sin(\sqrt{a_{10}}x_1),
				\\
				x_2\cos(\sqrt{a_{10}}x_1)\}
				\\{W}^1_{2,4}= \text{ Span}\{\sin(\sqrt{b_{10}}x_1),\cos(\sqrt{b_{10}}x_1),x_2\sin(\sqrt{b_{10}}x_1),
				\\
				x_2\cos(\sqrt{b_{10}}x_1)\}
			\end{array}\right.$
			\\
						6.&
			$\left\{\begin{array}{ll}
				k_{12} = -q_{11}, k_{22} = -q_{21},
			\\	q_{11},q_{21},q_{10},
				 c_{10},k_{20}\in\mathbb{R}
			\end{array}\right.$
			&
			$\left\{\begin{array}{ll}
				{W}^1_{1,4}= \text{Span}\{\sin(\sqrt{a_{20}}x_2),\cos(\sqrt{a_{20}}x_2),x_1\sin(\sqrt{a_{20}}x_2),
				\\
				x_1\cos(\sqrt{a_{20}}x_2)\}
				\\	{W}^1_{2,4}= \text{Span}\{1,x_1,
				x_2,x_1x_2\}
			\end{array}\right.$
			\\
					7.&
		$	\left\{\begin{array}{ll}
			  k_{12} = -q_{11}, k_{22} = -q_{21},
			\\	q_{11},q_{21},q_{10},k_{20}\in\mathbb{R}
				\end{array}\right.$
			&
			$\left\{\begin{array}{ll}
				{W}^1_{1,4}=  \text{Span}\{\sin(\sqrt{a_{20}}x_2),\cos(\sqrt{a_{20}}x_2),x_1\sin(\sqrt{a_{20}}x_2),
				\\
				x_1\cos(\sqrt{a_{20}}x_2)\}
				\\	{W}^1_{2,4}=  \text{Span}\{\sin(\sqrt{b_{20}}x_2),x_1\sin(\sqrt{b_{20}}x_2),\cos(\sqrt{b_{20}}x_2),
				\\ x_1\cos(\sqrt{b_{20}}x_2)\}
			\end{array}\right.$
			\\
			8.&
			$\left\{\begin{array}{ll}
			 k_{12} = -q_{11}, k_{22} = -q_{21},
			\\	q_{11},q_{21},q_{10},p_{20}\in\mathbb{R}
		\end{array}\right.$
			&
			$\left\{\begin{array}{ll}
				{W}^1_{1,4}=  
				\text{Span}\{1,x_1,x_2,x_1x_2\}
				\\
				{W}^1_{2,4}=   \text{Span}\{\sin(\sqrt{b_{20}}x_2),x_1\sin(\sqrt{b_{20}}x_2),
			\cos(\sqrt{b_{20}}x_2),
			\\
			x_1\cos(\sqrt{b_{20}}x_2)\}
			\end{array}\right.$
					\\
			9.
			&
			$\left\{\begin{array}{ll}
				p_{12} = -c_{11}, p_{22} = -c_{21},
			\\	q_{10} = -\frac{b_{20}c_{10}}{b_{10}},
				\\	c_{21},
				c_{11},c_{10}\in\mathbb{R}
			\end{array}\right.$
			&
			$\left\{\begin{array}{ll}
				{W}^1_{1,4}= 	\text{Span}\{\sin(\sqrt{a_{10}}x_1)\sin(\sqrt{b_{20}}x_2), \\
				\cos(\sqrt{a_{10}}x_1)\sin(\sqrt{b_{20}}x_2),
				\cos(\sqrt{a_{10}}x_1)	\cos(\sqrt{b_{20}}x_2),	\\
				\sin(\sqrt{a_{10}}x_1)	\cos(\sqrt{b_{20}}x_2)
		\}
				\\ {W}^1_{2,4}= \text{ Span}\{\sin(\sqrt{b_{10}}x_1)\sin(\sqrt{b_{20}}x_2),
			\\		
				\cos(\sqrt{b_{20}}x_2)\sin(\sqrt{b_{10}}x_1),
			 \cos(\sqrt{b_{10}}x_1)\sin(\sqrt{b_{20}}x_2),\\
			 \cos(\sqrt{b_{10}}x_1)\cos(\sqrt{b_{20}}x_2)\}
			\end{array}\right.$
			\\
			10. &
		$\left\{\begin{array}{ll}
			p_{12} = -c_{11},
			p_{22} = -c_{21} ,
			\\	c_{21},c_{11},
			c_{10},g_{10},q_{10}\in\mathbb{R}
		\end{array}\right.$
		&
		$\left\{\begin{array}{ll}
			{W}^1_{1,4}=    \text{Span}\{e^{b_{11}x_1},e^{-b_{11}x_1}, e^{b_{11}x_1}x_2,e^{-b_{11}x_1}x_2\}
			
			\\
			{W}^1_{2,4}= \text{Span}\{1,e^{-b_{11}x_1},x_2,x_2e^{-b_{11}x_1}\}
		\end{array}\right.
		$
			\\
		\bottomrule
	\end{tabular} }\label{table1}
	\end{table}
	%
	%
	%
	%
	%
	%
	%
	%
		\begin{table}[pt]
	\small	\caption{Type-1  product linear spaces $ {{	\mathbf{\hat{W}}^{1}_k}} $ and  their corresponding parameter restrictions of  $  \mathbf{F}\left(\mathbf{U}\right)  $ with $ a_{11}=b_{21}=0$ and $d_{10},\eta_{10},k_{10},p_{10},g_{20},h_{20},q_{20},c_{20}\in\mathbb{R}.$}
		{\begin{tabular}{@{}llllc@{}} \toprule
				Cases &	Some non-zero parameter 
				&  Type-1 invariant product linear space 
				\\	&	
				restrictions of
				$  \mathbf{F}\left(\mathbf{U}\right)  $ &  
				$ {{	\mathbf{\hat{W}}^{1}_8}}={{W^1_{1,4}}}\times {{W^1_{2,4}}} $ \\ \midrule
			11. &
			$\left\{\begin{array}{ll}
				p_{12} = -c_{11},
				p_{22}=-c_{21}, \\
				q_{10}=\frac{g_{10}}{b_{11}},
								c_{21},
				c_{11},
			\\	c_{10},
				g_{10}\in\mathbb{R}
			\end{array}\right.$
			&
			$\left\{\begin{array}{ll}
				{W}^1_{1,4}=    \text{Span}\{
				\sin(\sqrt{a_{10}}x_1),
				\cos(\sqrt{a_{10}}x_1),
				x_2\sin(\sqrt{a_{10}}x_1),
			\\	x_2\cos(\sqrt{a_{10}}x_1)\}
				\\
				{W}^1_{2,4}=    \text{Span}\{
				1,x_2, e^{-b_{11}x_1},
				x_2e^{-b_{11}x_1}\}
			\end{array}\right.
			$
						\\
			12. &
			$\left\{\begin{array}{ll}
				p_{12} = -c_{11},
				p_{22} = -c_{21} ,
				\\
				q_{10}=\frac{b_{11}g_{10}+b_{20}c_{10}}{b_{11}^2},
				\\
				g_{10},
				c_{21},
				c_{11},
				c_{10}
				\in\mathbb{R}
			\end{array}\right.$
			&
			$\left\{\begin{array}{ll}
				{W}^1_{1,4}=   \text{Span}\{
				\sin(\sqrt{b_{20}}x_2),
				\cos(\sqrt{b_{20}}x_2),
				x_1\sin(\sqrt{b_{20}}x_2),
			\\	x_1\cos(\sqrt{a_{20}}x_2)\}        \\
				{W}^1_{2,4}= \text{ Span}\{
				\sin(\sqrt{b_{20}}x_2),
						e^{-b_{11}x_1}\sin(\sqrt{b_{20}}x_2),
			\\	\cos(\sqrt{b_{20}}x_2),	e^{-b_{11}x_1}\cos(\sqrt{b_{20}}x_2)\}
			\end{array}\right.$
			\\
			
			13.	&
			$\left\{\begin{array}{ll}
				p_{12} = -c_{11},
				p_{22} = -c_{21} ,
			\\	q_{10}=\frac{g_{10}}{b_{11}},
					g_{10},
				c_{21},
				c_{11}\in\mathbb{R}
			\end{array}\right.$
			&
			$\left\{\begin{array}{ll}
				{W}^1_{1,4}=  \text{Span}\{
				\sin(\sqrt{b_{20}}x_2)\sin(\sqrt{a_{10}}x_1),
			\\	\cos(\sqrt{b_{20}}x_2)\sin(\sqrt{a_{10}}x_1),
					\sin(\sqrt{b_{20}}x_2)\cos(\sqrt{a_{10}}x_1),
			\\	\cos(\sqrt{a_{10}}x_1)\cos(\sqrt{b_{20}}x_2)\} \\
				{W}^1_{2,4}=
				\text{Span}\{	e^{-b_{11}x_1}\sin(\sqrt{b_{20}}x_2),
				\sin(\sqrt{b_{20}}x_2),
			\\	\cos(\sqrt{b_{20}}x_2),
				e^{-b_{11}x_1}\cos(\sqrt{b_{20}}x_2)\}
			\end{array}\right.
			$
		\\
			14.	&
			$\left\{\begin{array}{ll}
				p_{12} = -c_{11},
				p_{22} = -c_{21} ,	\\q_{10}=\frac{g_{10}}{b_{11}},
					c_{10}=-\frac{b_{10}g_{10}}{b_{11}b_{20}},
			\\	g_{10},
				c_{21},
				c_{11}\in\mathbb{R}
			\end{array}\right.$
			&
			$\left\{\begin{array}{ll}
				{W}^1_{1,4}= 
				\text{Span}\{
				\sin(\sqrt{b_{20}}x_2)\sin(\theta x_1),
			\\	\cos(\sqrt{b_{20}}x_2)\sin(\theta x_1),
					\sin(\sqrt{b_{20}}x_2)\cos(\theta x_1),	
			\\	\cos(\theta x_1)\cos(\sqrt{b_{20}}x_2)\},\theta= \sqrt{b_{10}-b_{11}^2}
				\\	{W}^1_{2,4}=\
				\text{Span}\{
				\sin(\sqrt{b_{20}}x_2)\sin(\sqrt{b_{10}}x_1),
			\\	\cos(\sqrt{b_{20}}x_2)\sin(\sqrt{b_{10}}x_1),
				\cos(\sqrt{b_{10}}x_1)\sin(\sqrt{b_{20}}x_2),
				\\				
				\cos(\sqrt{b_{10}}x_1)\cos(\sqrt{b_{20}}x_2)\}
			\end{array}\right.
			$
		\\
			15.	&
			$\left\{\begin{array}{ll}
				k_{12}=-\frac{c_{11}a_{20}}{a_{10}},
			\\
				q_{11}=\frac{c_{11}a_{20}}{a_{10}},	\\
				p_{12} = -c_{11},	c_{21},
				c_{11},
			\\	q_{10},
				c_{10}\in\mathbb{R},
			\end{array}\right.$
			&
			$\left\{\begin{array}{ll}
				{W}^1_{1,4}=  \text{Span}\{\sin(\sqrt{a_{10}}x_1)\sin(\sqrt{a_{20}}x_2),
				\\
				\cos(\sqrt{a_{10}}x_1)\sin(\sqrt{a_{20}}x_2),
				\cos(\sqrt{a_{20}}x_2)\sin(\sqrt{a_{10}}x_1),
				\\\cos(\sqrt{a_{20}}x_2)\cos(\sqrt{a_{10}}x_1)\}        	
				\\	{W}^1_{2,4}= 
				\text{Span}\{1,x_1,x_2,x_1x_2\}
			\end{array}\right.
			$\\
			16. &
			$\left\{\begin{array}{ll}
				p_{12} = -c_{11},
				q_{11}=\frac{c_{11}k_{20}}{p_{20}},\\
				k_{12}=-\frac{c_{11}k_{20}}{p_{20}},	\\
				c_{10},
				p_{20},
				q_{10},
				k_{20}\in\mathbb{R}
			\end{array}\right.$
			&
			$\left\{\begin{array}{ll}
				{W}^1_{1,4}=  \text{Span}\{e^{\theta x_1}\sin(\sqrt{a_{20}}x_2),
	e^{\theta x_1}		\cos(\sqrt{a_{20}}x_2),\\
		e^{-\theta x_1}\sin(\sqrt{a_{20}}x_2),	
	e^{-\theta x_1}	\cos(\sqrt{a_{20}}x_2)\},
	\\
	\theta=\sqrt{\frac{p_{20}a_{20}}{k_{20}}}
				\\
				{W}^1_{2,4}=\text{Span}\{1,x_1,x_2,x_1x_2\}
			\end{array}\right.
			$
			\\
			17.&
			$\left\{\begin{array}{ll}
				k_{12} = -\frac{c_{11}b_{20}}{a_{10}},
							\\	q_{11} = \frac{c_{11}b_{20}}{a_{10}},
										p_{12} = -c_{11},
										\\
										c_{11}, p_{10},q_{10}\in\mathbb{R}
			\end{array}\right.$
			&
			$\left\{\begin{array}{ll}
				{W}^1_{1,4}= \text{Span}\{\sin(\sqrt{a_{10}}x_1),
				\cos(\sqrt{a_{10}}x_1),
				x_2\sin(\sqrt{a_{10}}x_1),
			\\	x_2\cos(\sqrt{a_{10}}x_1)\}
				\\{W}^1_{2,4}= 	\text{Span}\{\sin(\sqrt{b_{20}}x_2),x_1\sin(\sqrt{b_{20}}x_2),\cos(\sqrt{b_{20}}x_2),
				\\x_1\cos(\sqrt{b_{20}}x_2)\} \end{array}\right.$
						\\
						18.	&
				$\left\{\begin{array}{ll}
				k_{12} =\frac{c_{11}(a_{20}-b_{20})}{a_{10}},
\\		q_{11} = \frac{(b_{20}-a_{20})c_{11}}{a_{10}},\\					p_{12} = -c_{11},										
			c_{11},	q_{10}\in\mathbb{R}
				\end{array}\right.$
						&
					$\left\{\begin{array}{ll}
			{W}^1_{1,4}= 	\text{Span}\{
		\sin(\sqrt{a_{10}}x_1)\sin(\sqrt{a_{20}}x_2),
\\											\cos(\sqrt{a_{10}}x_1)\sin(\sqrt{a_{20}}x_2),
	\cos(\sqrt{a_{20}}x_2)\sin(\sqrt{a_{10}}x_1),	
\\											\cos(\sqrt{a_{20}}x_2)\cos(\sqrt{a_{10}}x_1)\} \\
								{W}^1_{2,4}= 
								\text{Span}\{
						\sin(\sqrt{b_{20}}x_2),
						\cos(\sqrt{b_{20}}x_2),
						x_1\sin(\sqrt{b_{20}}x_2),
\\							x_1\cos(\sqrt{b_{20}}x_2)\}
										\end{array}\right.$
									\\
						19.&
				$\left\{\begin{array}{ll}
				k_{12} = -\frac{q_{21}c_{11}}{c_{21}},
\\							k_{22} = -q_{21},
									p_{12} = -c_{11},
									\\
							p_{22} = -c_{21},	
					q_{11}=\frac{q_{21}c_{11}}{c_{21}},
\\											
											c_{11},c_{21},q_{21},q_{10}\in\mathbb{R}
										\end{array}\right.$
										&
										$\left\{\begin{array}{ll}
											{W}^1_{1,4}=  \text{Span}\{
											\sin(\sqrt{\theta}x_1)\sin(\sqrt{b_{20}}x_2),
										\\	\cos(\sqrt{\theta}x_1)\sin(\sqrt{b_{20}}x_2),
		\cos(\sqrt{b_{20}}x_2)\sin(\sqrt{\theta}x_1),	\\
		\cos(\sqrt{b_{20}}x_2)\cos(\sqrt{\theta}x_1)\},\theta=\frac{(b_{20}-a_{20})c_{21}}{q_{21}}>0
	\\{W}^1_{2,4}=\text{Span}\{
			\sin(\sqrt{b_{20}}x_2),
											\cos(\sqrt{b_{20}}x_2),
											x_1\sin(\sqrt{b_{20}}x_2),
										\\	x_1\cos(\sqrt{b_{20}}x_2)\}
										\end{array}\right.$
										\\
										 \bottomrule
		\end{tabular} }\label{table2}
	\end{table}
%
%
%
%
%
%
%
%
\begin{table}[pt]
\small \caption{Type-1  product linear spaces $ {{	\mathbf{\hat{W}}^{1}_k}} $ and  their corresponding parameter restrictions of  $  \mathbf{F}\left(\mathbf{U}\right)  $ with $ a_{11}=b_{21}=0$ and $d_{10},\eta_{10},k_{10},p_{10},g_{20},h_{20},q_{20},c_{20}\in\mathbb{R}.$}
{\begin{tabular}{@{}llllc@{}} \toprule
Cases &	Some non-zero parameter 
&  Type-1 invariant product linear space 
\\	&	
restrictions of
$  \mathbf{F}\left(\mathbf{U}\right)  $ &  
$ {{	\mathbf{\hat{W}}^{1}_8}}={{W^1_{1,4}}}\times {{W^1_{2,4}}} $ \\ \midrule
			20. &
			$\left\{\begin{array}{ll}
				p_{12} = -c_{11},	p_{22}=-c_{21},
			\\	k_{12}=-\frac{c_{11}q_{21}}{c_{21}},		\\
					q_{11}=\frac{c_{11}q_{21}}{c_{21}},		k_{22}=-q_{21},\\	k_{20}=\frac{p_{20}q_{21}}{c_{21}},
							
					c_{11},
				c_{10},\\
				p_{20},	
				c_{21},
				q_{10},
				q_{21}\in\mathbb{R}
			\end{array}\right.$
			&
			$\left\{\begin{array}{ll}
				{W}^1_{1,4}= \text{Span}\{
				e^{ \theta x_1}\sin(\sqrt{a_{20}}x_2),\\
				e^{-\theta x_1}\sin(\sqrt{a_{20}}x_2),
				e^{ \theta x_1}\cos(\sqrt{a_{20}}x_2),	\\
				\cos(\sqrt{a_{20}}x_2)e^{- \theta x_1}\},
				\theta=\sqrt{\frac{c_{21}a_{20}}{q_{21}}}
				\\	{W}^1_{2,4}= 
				\text{Span}\{1,x_1,x_2,x_1x_2\}
			\end{array}\right.
			$
		\\
			21.&
			$\left\{\begin{array}{ll}
				k_{12} = -\frac{q_{21}c_{11}}{c_{21}} ,
			\\	k_{22} = -q_{21}, p_{12} = -c_{11}, \\ p_{22} = -c_{21},  q_{11} = \frac{q_{21}c_{11}}{c_{21}},\\
			p_{20},q_{10},c_{11},c_{21}\in\mathbb{R}
			\end{array}\right.$
			&
			$\left\{\begin{array}{ll}
				{W}^1_{1,4}=  \text{Span}\{\sin(\sqrt{\theta}x_1),
				\cos(\sqrt{\theta}x_1),
				x_2\sin(\sqrt{\theta}x_1),
			\\	x_2\cos(\sqrt{\theta}x_1)\}, \theta={\frac{c_{21}b_{20}}{q_{21}}}>0
				\\	{W}^1_{2,4}= \text{Span}\{\sin(\sqrt{b_{20}}x_2),\cos(\sqrt{b_{20}}x_2),
				x_1\sin(\sqrt{b_{20}}x_2),
				\\
				x_1\cos(\sqrt{b_{20}}x_2)\} \end{array}\right.$
		\\
			22.&
			$\left\{\begin{array}{ll}
				k_{12} = \frac{(b_{20}-a_{20})c_{11}k_{20}}{(a_{20}p_{20})},
			\\				
				q_{11} =\frac{c_{11}(a_{20}-b_{20})k_{20}}{(a_{20}p_{20})}
				,
				\\
				p_{12} = -c_{11},p_{20},k_{20},\\
				q_{10},c_{11}\in\mathbb{R}
			\end{array}\right.$
			&
			$\left\{\begin{array}{ll}
				{W}^1_{1,4}=  \text{Span}\{e^{\sqrt{\theta}x_1}\sin(\sqrt{a_{20}}x_2),
				e^{-\sqrt{\theta}x_1}\sin(\sqrt{a_{20}}x_2),
			\\	\cos(\sqrt{a_{20}}x_2)e^{\sqrt{\theta}x_1},	
					\cos(\sqrt{a_{20}}x_2)e^{-\sqrt{\theta}x_1}\},
					\\
					{\theta}={\frac{p_{20}a_{20}}{k_{20}}}>0
				\\ {W}^1_{2,4}= \text{Span}\{\sin(\sqrt{b_{20}}x_2),\cos(\sqrt{b_{20}}x_2),
				x_1\sin(\sqrt{b_{20}}x_2),
				\\
				x_1\cos(\sqrt{b_{20}}x_2)\}
			\end{array}\right.$
			\\	23.&
			$\left\{\begin{array}{ll}
						k_{22} = -q_{21}, p_{12} = -c_{11},\\
				p_{22} = -c_{21},
							q_{11} = \frac{q_{21}c_{11}}{c_{21}}
			\\	k_{12} = -\frac{q_{21}c_{11}}{c_{21}},	p_{20}, \\
			k_{20},
				q_{10},q_{21},c_{21},c_{11}\in\mathbb{R}
			\end{array}
			\right.$
			&
			$\left\{\begin{array}{ll}
				{W}^1_{1,4}=  \text{Span}\big\{e^{\sqrt{\theta}(x_1+x_2)},
				e^{-\sqrt{\theta}(x_1+x_2)},
				e^{\sqrt{\theta}(x_1-x_2)},	
			\\	e^{-\sqrt{\theta}(x_1-x_2)}\big\},\theta={\frac{c_{21}b_{20}p_{20}}{k_{20}c_{21}-p_{20}q_{21}}}>0
				\\ 	{W}^1_{2,4}= \text{Span}\{\sin(\sqrt{b_{20}}x_2),\cos(\sqrt{b_{20}}x_2),
				x_1\sin(\sqrt{b_{20}}x_2),
				\\
				x_1\cos(\sqrt{b_{20}}x_2)\}
			\end{array}\right.$
			\\
			24.& $ \left\{\begin{array}{ll}
				\left\{\begin{array}{ll} k_{12} = -q_{11}, k_{22} = -q_{21},  
					\\ p_{22} = -c_{21},
					
					d_{20},h_{10},\\
					p_{20},c_{21},c_{10},	q_{11},\\
					q_{21},q_{10},k_{20}\in\mathbb{R}
				\end{array} \right.\\
				\left\{\begin{array}{ll} k_{12} = -q_{11}, k_{22} = -q_{21}, \\  p_{12} = -c_{11}, p_{22} = -c_{21},
					\\
					d_{20},h_{10},p_{20},c_{21},c_{11},
					c_{10},\\	q_{11},q_{21},q_{10},k_{20}\in\mathbb{R}
				\end{array} \right.
						\end{array}
			\right.$
			&
			$\left\{\begin{array}{ll}
				{W}^1_{1,4}=  \text{Span}\{\sin(\sqrt{b_{10}}x_1)\sin(\sqrt{b_{20}}x_2),
				\\
				\cos(\sqrt{b_{10}}x_1)\sin(\sqrt{b_{20}}x_2),
				\cos(\sqrt{b_{20}}x_2)\sin(\sqrt{b_{10}}x_1),
				\\
				 \cos(\sqrt{b_{20}}x_2)\cos(\sqrt{b_{10}}x_1)\}
				
				\\
				{W}^1_{2,4}= \text{Span}\{\sin(\sqrt{b_{10}}x_1)\sin(\sqrt{b_{20}}x_2),
				\\
				\cos(\sqrt{b_{10}}x_1)\sin(\sqrt{b_{20}}x_2),
								\cos(\sqrt{b_{20}}x_2)\sin(\sqrt{b_{10}}x_1),
								\\\cos(\sqrt{b_{20}}x_2)\cos(\sqrt{b_{10}}x_1)\}
			\end{array}\right.$
		\\
			25.&
		$\left\{\begin{array}{ll}
			  k_{12} = -q_{11},
			k_{22} = -q_{21},\\
			q_{10} = -\frac{b_{20}c_{10}}{b_{10}},
			\\	c_{10},	q_{11},q_{21},q_{10}\in\mathbb{R}
		\end{array}\right.
		$
		&
		$\left\{\begin{array}{ll}
			{W}^1_{1,4}=  \text{Span}\{\sin(\sqrt{b_{10}}x_1)\sin(\sqrt{a_{20}}x_2),\\
			\cos(\sqrt{b_{10}}x_1)\sin(\sqrt{a_{20}}x_2),
			\cos(\sqrt{a_{20}}x_2)\sin(\sqrt{b_{10}}x_1),
			\\\cos(\sqrt{a_{20}}x_2)\cos(\sqrt{b_{10}}x_1)\}
			\\ {W}^1_{2,4}=  \text{Span}\{\sin(\sqrt{b_{10}}x_1)\sin(\sqrt{b_{20}}x_2),\\
			\cos(\sqrt{b_{10}}x_1)\sin(\sqrt{b_{20}}x_2),
					\cos(\sqrt{b_{20}}x_2)\sin(\sqrt{b_{10}}x_1),
					\\
					\cos(\sqrt{b_{20}}x_2)\cos(\sqrt{b_{10}}x_1)\}
		\end{array}\right.$
		\\
		26.&
		$\left\{\begin{array}{ll}
			k_{12} = -\frac{q_{21}c_{11}}{c_{21}},
			k_{22} = -q_{21},\\
			q_{11} = \frac{q_{21}c_{11}}{c_{21}},	p_{12} = -c_{11},\\ 
			q_{10} = -\frac{b_{20}c_{10}}{b_{10}},	p_{22} = -c_{21},\\
			p_{20},q_{21},c_{11},c_{21},c_{10}\in\mathbb{R}
		\end{array}\right.$
		&
		$\left\{\begin{array}{ll}
			{W}^1_{1,4}= 
			\text{Span}\{\sin(\sqrt{\theta}x_1),
			\cos(\sqrt{\theta}x_1),
			x_2\sin(\sqrt{\theta}x_1),\\
			x_2\cos(\sqrt{\theta}x_1)\},
			{\theta}={\frac{ (b_{10}q_{21}+b_{20}c_{21})}{q_{21}}}>0
			\\	{W}^1_{2,4}=  \text{Span}\{\sin(\sqrt{b_{10}}x_1)\sin(\sqrt{b_{20}}x_2),
			\\
			\cos(\sqrt{b_{10}}x_1)\sin(\sqrt{b_{20}}x_2),
			
			\cos(\sqrt{b_{20}}x_2)\sin(\sqrt{b_{10}}x_1), \\
			\cos(\sqrt{b_{20}}x_2)\cos(\sqrt{b_{10}}x_1)\}
		\end{array}\right.$
		\\
		\bottomrule
		\end{tabular}}\label{table3}
	\end{table}
%
%
%
%
%
%
%
%
	\begin{table}[htp!]
	\small\caption{Type-1  product linear spaces $ {{	\mathbf{\hat{W}}^{1}_k}} $ and  their corresponding parameter restrictions of  $  \mathbf{F}\left(\mathbf{U}\right)  $ with $ a_{11}=b_{21}=0$ and $d_{10},\eta_{10},k_{10},p_{10},g_{20},h_{20},q_{20},c_{20}\in\mathbb{R}.$}
	{\begin{tabular}{@{}llllc@{}} \toprule
			Cases &	Some non-zero parameter 
			&  Type-1 invariant product linear space 
			\\	&	
			restrictions of
			$  \mathbf{F}\left(\mathbf{U}\right)  $ &  
			$ {{	\mathbf{\hat{W}}^{1}_8}}={{W^1_{1,4}}}\times {{W^1_{2,4}}} $ \\ 
			\midrule
			27.&
			$\left\{\begin{array}{ll}
				k_{12} = \frac{b_{20}c_{11}}{b_{10}-a_{10}},\\
							q_{10} = -\frac{b_{20}c_{10}}{b_{10}}, \\
				p_{12} = -c_{11},p_{20},c_{11},\\
				c_{10}\in\mathbb{R},	q_{11} = \frac{c_{11}b_{20}}{a_{10}-b_{10}}
			\end{array}\right.$
			&
			$\left\{\begin{array}{ll}
				{W}^1_{1,4}= \text{Span}\{\sin(\sqrt{a_{10}}x_1),x_2\sin(\sqrt{a_{10}}x_1),
				\\
				\cos(\sqrt{a_{10}}x_1),x_2\cos(\sqrt{a_{10}}x_1)\}
				\\	{W}^1_{2,4}= \text{Span}\{\sin(\sqrt{b_{10}}x_1)\sin(\sqrt{b_{20}}x_2),
				\\
				\cos(\sqrt{b_{10}}x_1)\sin(\sqrt{b_{20}}x_2),
				\cos(\sqrt{b_{20}}x_2)\sin(\sqrt{b_{10}}x_1),
						\\
						\cos(\sqrt{b_{20}}x_2)\cos(\sqrt{b_{10}}x_1)\}
			\end{array}\right.$
			\\
					
			28.&
			$\left\{\begin{array}{ll}
				k_{12} = -q_{11}, k_{22} = -q_{21},\\
				q_{10} = -\frac{b_{20}c_{10}}{b_{10}},
				c_{10},	q_{11},\\
				q_{21},q_{10},p_{20},k_{20}\in\mathbb{R}
			\end{array}\right.$
			&
			$\left\{\begin{array}{ll}
				{W}^1_{1,4}= \text{Span}\{\sin(\sqrt{b_{10}}x_1) e^{\sqrt{\theta}x_2}, \\\cos(\sqrt{b_{10}}x_1) e^{\sqrt{\theta}x_2},
				e^{-\sqrt{\theta}x_2}
				\sin(\sqrt{b_{10}}x_1),
				\\	e^{-\sqrt{\theta}x_2}
				\cos(\sqrt{b_{10}}x_1)\},
				{\theta}={\frac{k_{20}b_{10}}{p_{20}}}>0
				
				\\
				{W}^1_{2,4}= \text{Span}\{\sin(\sqrt{b_{10}}x_1)\sin(\sqrt{b_{20}}x_2),\\
				\cos(\sqrt{b_{10}}x_1)\sin(\sqrt{b_{20}}x_2),	\cos(\sqrt{b_{20}}x_2)\sin(\sqrt{b_{10}}x_1),
							\\
							\cos(\sqrt{b_{20}}x_2)\cos(\sqrt{b_{10}}x_1)\}
			\end{array}\right.$
			\\
			29.
			&
			$\left\{\begin{array}{ll}
				k_{12} = \frac{c_{11}(a_{20}-b_{20})}{(a_{10}-b_{10})},\\
				p_{12} = -c_{11},
				q_{10} = -\frac{b_{20}c_{10}}{b_{10}}, \\
				q_{11} =\frac{ (b_{20}-a_{20})c_{11}}{(a_{10}-b_{10})},\\
				c_{11},c_{10}\in\mathbb{R}
			\end{array}\right.$
			&
			$\left\{\begin{array}{ll}
				{W}^1_{1,4}= \text{Span}\{\sin(\sqrt{a_{10}}x_1)\sin(\sqrt{a_{20}}x_2),\\
				\cos(\sqrt{a_{10}}x_1)\sin(\sqrt{a_{20}}x_2),
				\cos(\sqrt{a_{20}}x_2)\sin(\sqrt{a_{10}}x_1),	\\
				\cos(\sqrt{a_{20}}x_2)\cos(\sqrt{a_{10}}x_1)\}
				
				\\
				{W}^1_{2,4}= \text{Span}\{\sin(\sqrt{b_{10}}x_1)\sin(\sqrt{b_{20}}x_2),\\
				\cos(\sqrt{b_{10}}x_1)\sin(\sqrt{b_{20}}x_2),
				\cos(\sqrt{b_{20}}x_2)\sin(\sqrt{b_{10}}x_1),
				\\
				\cos(\sqrt{b_{20}}x_2)\cos(\sqrt{b_{10}}x_1)\}
			\end{array}\right.$
		\\
			30.
			&
			$\left\{\begin{array}{ll}
				k_{12} = -\frac{q_{21}c_{11}}{c_{21}},
				k_{22} = -q_{21}, \\p_{12} = -c_{11},
					p_{22} = -c_{21},\\
				q_{10} = -\frac{b_{20}c_{10}}{b_{10}},
				\\	q_{11} = \frac{q_{21}c_{11}}{c_{21}},
				q_{21},\\	c_{21},
				c_{11},c_{10}\in\mathbb{R}
			\end{array}
			\right.$
			&
			$\left\{\begin{array}{ll}
				{W}^1_{1,4}= \text{Span}\{\sin(\sqrt{\theta}x_1)
				\sin(\sqrt{a_{20}}x_2),
			\\	\cos(\sqrt{\theta}x_1)	\sin(\sqrt{a_{20}}x_2),
								\cos(\sqrt{a_{20}}x_2)\sin(\sqrt{\theta}x_1),
								\\
				\cos(\sqrt{a_{20}}x_2)\cos(\sqrt{\theta}x_1)\},
				{\theta}={\frac{ (b_{10}q_{21}-c_{21}(a_{20}-b_{20}))}{q_{21}}}>0
				\\
				{W}^1_{2,4}= \text{Span}\{\sin(\sqrt{b_{10}}x_1)\sin(\sqrt{b_{20}}x_2),\\
				\cos(\sqrt{b_{10}}x_1)\sin(\sqrt{b_{20}}x_2),
				\cos(\sqrt{b_{20}}x_2)\sin(\sqrt{b_{10}}x_1),	
				\\
				\cos(\sqrt{b_{20}}x_2)\cos(\sqrt{b_{10}}x_1)\}
			\end{array}
			\right.$
			\\
			31.
			&
			$\left\{\begin{array}{ll}
				k_{12}= -q_{11},  k_{22}= -q_{21},\\ p_{12}= \frac{d_{12}p_{20}}{d_{20}},\\  c_{11}= -\frac{d_{12}p_{20}}{d_{20}},
								d_{12},
				d_{20},
				h_{10},\\
				k_{20},
				p_{20},
				c_{10},	
				q_{10} ,
				q_{11},
				q_{21}\in\mathbb{R}
			\end{array}\right.$
			&
			$\left\{\begin{array}{ll}
				{W}^1_{1,4}= &\text{Span}\{1,x_1,e^{-\frac{d_{20}}{p_{20}}x_2},x_1 e^{-\frac{d_{20}}{p_{20}}x_2}
				\}
				\\	{W}^1_{2,4}= &\text{Span}\{1,x_1,x_2,x_1x_2
				\}
			\end{array}\right.$
		\\
			32.
			&
			$\left\{\begin{array}{ll}
				k_{12}= -  \frac{g_{11}}{b_{11}},
				k_{22}= -  \frac{g_{21}}{b_{11}},\\
				q_{10}=   \frac{g_{10}}{b_{11}},
				q_{11}=   \frac{g_{11}}{b_{11}},
				\\	q_{21}=   \frac{g_{21}}{b_{11}},
				g_{11},
				g_{21},\\
				g_{10},
				c_{10},
				k_{20}
				\in\mathbb{R}
			\end{array}\right.
			$
			&
			$\left\{\begin{array}{ll}
				{W}^1_{1,4}= \text{Span}\{
				 	e^{-\theta_1x_2},	
				x_1e^{-\theta_1x_2},e^{-\theta_2x_2},
			\\	x_1e^{-\theta_2x_2}
				\}, \theta_1=\frac{a_{21}+\sqrt{a_{21}^2-4a_{20}}}{2},
				\\
				\theta_2=\frac{a_{21}-\sqrt{a_{21}^2-4a_{20}}}{2}
				\\	{W}^1_{2,4}= 	\text{Span}\{1,e^{-b_{11}x_1},x_2,x_2e^{-b_{11}x_1}\}
			\end{array}\right.$
			\\
			\\
			33.
			&
			$\left\{\begin{array}{ll}
				k_{12}= -q_{11},  k_{22}= -q_{21},\\
				q_{11},
				q_{21},q_{10},
				c_{10},
				k_{20},
				\in\mathbb{R}
			\end{array}\right.
			$
			&
			$\left\{\begin{array}{ll}
					{W}^1_{1,4}= \text{Span}\{
				e^{-\theta_1x_2},	
				x_1e^{-\theta_1x_2},e^{-\theta_2x_2},
				\\	x_1e^{-\theta_2x_2}
				\}, \theta_1=\frac{a_{21}+\sqrt{a_{21}^2-4a_{20}}}{2},
				\\
				\theta_2=\frac{a_{21}-\sqrt{a_{21}^2-4a_{20}}}{2}
				\\	{W}^1_{2,4}=  	\text{Span}\{1,x_1,x_2,x_1x_2\}
				
			\end{array}
			\right.$
			\\
					\bottomrule
	\end{tabular}}\label{table4}
	\end{table}
%
%
%
%
%
%
%
%
	\begin{table}[htp!]
		\small \caption{Type-1  product linear spaces $ {{	\mathbf{\hat{W}}^{1}_k}} $ and  their corresponding parameter restrictions of  $  \mathbf{F}\left(\mathbf{U}\right)  $ with $ a_{11}=b_{21}=0$ and $d_{10},\eta_{10},k_{10},p_{10},g_{20},h_{20},q_{20},c_{20}\in\mathbb{R}.$}
		{\begin{tabular}{@{}llllc@{}} \toprule
				Cases &	Some non-zero parameter 
				&  Type-1 invariant product linear space 
				\\	&	
				restrictions of
				$  \mathbf{F}\left(\mathbf{U}\right)  $ &  
				$ {{	\mathbf{\hat{W}}^{1}_8}}={{W^1_{1,4}}}\times {{W^1_{2,4}}} $ \\ \midrule
			34.			&
			$\left\{\begin{array}{ll}
				k_{i2}= -\frac{g_{i1}}{b_{11}},  p_{i2}= -c_{i1},
				\\
				q_{i1}= \frac{g_{i1}}{b_{11}}, i=1,2,
				\\ q_{10}= \frac{g_{10}}{b_{11}},
				\eta_{20},
				g_{10},
								g_{21},g_{11},\\
				p_{20},
				k_{20},
				c_{10},c_{11}, c_{21}\in\mathbb{R}
			\end{array}\right.$
			&
			$\left\{\begin{array}{ll}
				{W}^1_{1,4}= \text{Span}\{1,x_1,x_2,x_1x_2
				\}
				\\{W}^1_{2,4}=
				\text{Span}\{1,e^{-b_{11}x_1},x_2,x_2e^{-b_{11}x_1}\}
			\end{array}\right.$
		\\
		\\
			35.
			&
			$\left\{\begin{array}{ll}
			
				 p_{12}= -\frac{d_{12}c_{21}}{d_{22}},k_{12}= -q_{11},\\
				p_{20}= -\frac{d_{20}c_{21}}{d_{22}},  	k_{22}= -q_{21},\\ c_{11}= \frac{d_{12}c_{21}}{d_{22}},
			p_{22}= -c_{21}, \\	d_{12},
				d_{22},
				d_{20},
				h_{10},
			k_{20},
		c_{10},\\	c_{21},	
				q_{10} ,
			q_{11},
				q_{21}\in\mathbb{R}
			\end{array}\right.$
			&
			$\left\{\begin{array}{ll}
				{W}^1_{1,4}= \text{Span}\{1,x_1,e^{\frac{d_{22}}{c_{21}}x_2},x_1 e^{\frac{d_{22}}{c_{21}}x_2}
				\\ {W}^1_{2,4}=\text{Span}\{1,x_1,x_2,x_1x_2		\}
			\end{array}\right.$
			\\
			\\
			36.
			&
			$\left\{\begin{array}{ll}
				k_{12}= -\frac{g_{11}}{b_{11}},
				k_{22}= -\frac{g_{21}}{b_{11}},\\
				p_{12}= \frac{d_{12}}{a_{21}},
				q_{10}= \frac{g_{10}}{b_{11}},\\
				q_{11}=\frac{g_{11}}{b_{11}},
				q_{21}=\frac{g_{21}}{b_{11}},
				\\	c_{11}= - \frac{d_{12}}{a_{21}}  ,
				d_{12},	k_{20}, g_{10},\\
				g_{21},g_{11},
				q_{10},q_{11}, q_{21},c_{10}\in\mathbb{R}
			\end{array}\right.$
			&
			$\left\{\begin{array}{ll}
				{W}^{1}_{1,4}= = \text{Span}\{1,x_1,e^{-a_{21}x_2},x_1e^{-a_{21}x_2}
				\}
				\\	{W}^{1}_{1,4}= 
				\text{Span}\{1,e^{-b_{11}x_1},x_2,e^{-b_{11}x_1}x_2\}
			\end{array}\right.$
			\\
			\\
			37.
			&
			$\left\{\begin{array}{ll}
				k_{12}= -\frac{g_{11}}{b_{11}},
				k_{22}= -\frac{g_{21}}{b_{11}},\\
				p_{12}= \frac{d_{12}p_{20}}{d_{20}},
				q_{10}= \frac{g_{10}}{b_{11}},\\
				q_{11}=\frac{g_{11}}{b_{11}},
				q_{21}= \frac{g_{21}}{b_{11}},
				\\ 	c_{11}= -\frac{d_{12}p_{20}}{d_{20}},
								d_{12},d_{20},	k_{20},p_{20},\\ g_{10},
				g_{21},g_{11},
				q_{10},q_{11}, q_{21},c_{10}\in\mathbb{R}
			\end{array}\right.$
			&
			$\left\{\begin{array}{ll}
				{W}^{1}_{1,4}=
				&\text{Span}\{1,x_1,e^{-\frac{d_{20}}{p_{20}}x_2},x_1 e^{-\frac{d_{20}}{p_{20}}x_2}
				\}
				\\	{W}^{1}_{2,4}=  &
				\text{Span}\{1,e^{-b_{11}x_1},x_2,e^{-b_{11}x_1}x_2\}
			\end{array}\right.$
			\\
			\\
			38.
			&
			$\left\{\begin{array}{ll}
				k_{12}= -\frac{g_{11}}{b_{11}},
				k_{22}= -\frac{g_{21}}{b_{11}},\\
				p_{12}= -\frac{d_{12}c_{21}}{d_{22}},
				p_{20}= -\frac{d_{20}c_{21}}{d_{22}},	\\
							q_{10}= \frac{g_{10}}{b_{11}},
				q_{11}= \frac{g_{11}}{b_{11}},\\
					p_{22}= -c_{21},c_{11}= \frac{d_{12}c_{21}}{d_{22}}, \\	d_{22},d_{12},d_{20},
				g_{10},g_{21},	
				g_{11},	k_{20},\\
				q_{21},
				c_{10},c_{21}\in\mathbb{R}
			\end{array}\right.
			$
			&
			$\left\{\begin{array}{ll}
				{W}^{1}_{1,4}=  &\text{Span}\{1,x_1,e^{\frac{d_{22}}{c_{21}}x_2},x_1 e^{\frac{d_{22}}{c_{21}}x_2}
				\}
				\\
				{W}^{1}_{2,4}= &
				\text{Span}\{1,e^{-b_{11}x_1},x_2,e^{-b_{11}x_1}x_2\}
			\end{array}
			\right.$
			\\	\\	39.
			&
			$\left\{\begin{array}{ll}
				k_{12}= -q_{11},k_{22}= -q_{21}, \\
				p_{12}= -c_{11},
				p_{22}= -c_{21},\\	d_{20},\eta_{20},
				k_{20},
				p_{20},
				g_{10} ,\\
				h_{10} ,
				q_{10} ,
				q_{11}, 	
				q_{21},	c_{21},
				c_{11},c_{10}\in\mathbb{R}
			\end{array}\right.$
			&
			$\left\{\begin{array}{ll}
				{W}^{1}_{1,4}= &\text{Span}\{1,x_1,x_2,x_1x_2
				\}
				\\
				{W}^{1}_{2,4}= &\text{ Span}\{1,x_1,x_2,x_1x_2
				\}
			\end{array}\right.$
			\\
			\\	40.
			&
			$\left\{\begin{array}{ll}
				k_{12}= -q_{11}, k_{22}= -q_{21}, p_{12}= \frac{d_{12}}{a_{21}},  c_{11}= -\frac{d_{12}}{a_{21}},
				\\
				d_{12},k_{20},
				h_{10} ,	q_{10} ,
				q_{11},
				q_{21},	c_{10}\in\mathbb{R}
			\end{array}\right.$
			&
			$\left\{\begin{array}{ll}
				{W}^{1}_{1,4}= &\text{Span}\{1,x_1,e^{-a_{21}x_2},x_1e^{-a_{21}x_2}
				\}
				\\
				{W}^{1}_{2,4}=&\text{Span}\{1,x_1,x_2,x_1x_2
				\}
			\end{array}\right.$
			\\
		
		\bottomrule
	\end{tabular}}
		\label{table5}
	\end{table}
	\subsection{Generalized separable {analytical} solutions for the given  coupled equations \eqref{cdsystemeqn}}
	This subsection provides a systematic study for deriving the generalized separable {analytical} solutions for the given coupled equations \eqref{cdsystemeqn}  using the {type-1 and type-2 invariant product linear spaces} discussed in the previous section.
	\subsubsection{Generalized separable {analytical} solutions for the initial value problems of the given coupled equations  \eqref{cdsystemeqn}}
	First, we discuss the systematic procedure for finding generalized separable {analytical} solutions for the initial value problems of the considered coupled equations \eqref{cdsystemeqn} using the obtained {type-1 and type-2 invariant product linear spaces}.
	\begin{example}\label{eg1}
		Consider the following  {multi-component} $(2+1)$-dimensional coupled   nonlinear {TFDCWEs}  as
		{	\begin{eqnarray}
				\begin{aligned}\label{poly-polyeq}
					{	\dfrac{\partial ^{\alpha_s} u_s}{\partial t^{\alpha_s}}}
					=
					F_s(u_1,u_2)
					\equiv&
					{\frac {\partial }{\partial x_1}}
					\left[	
					\left(
					k_{{s0}}-q_{{s1}}u_{{2}}
					\right)
					{\frac {\partial u_{{1}} }{	\partial x_1}}
					+
					\left(
					q_{{s1}}u_{{1}}
					+q_{{s0}}
					\right)
					{\frac {\partial u_{{2}} }{\partial x_1}}
					\right]
			\\&		+
					{\frac {\partial }{	\partial x_2}}
					\left[
					\left(
					p_{{s0}}-c_{{s1}}u_{{2}}
					\right)
					{\frac {\partial u_{{1}}}{\partial x_2}}
					+
					\left(
					c_{{s1}}u_{{1}}
					+c_{{s0}}
					\right)
					{\frac {\partial u_{{2}}}{\partial x_2}}
					\right]
					+
					\eta_{{s0}} {\frac {\partial u_{{1}}}{\partial x_1}}
				\\&	+
					g_{{s0}}
					{\frac {\partial u_{{2}} }{\partial x_1}}
					+
					d_{{s0}}
					{\frac {\partial u_{{1}}  }{\partial x_2}}
					+
					h_{{s0}}
					{\frac {\partial u_{{2}} }{\partial x_2}},{\alpha_s}\in(0,2],s=1,2,
				\end{aligned}
		\end{eqnarray}}
		with appropriate initial conditions
		{\begin{eqnarray}
				\begin{aligned}\label{icpoly-polyeq}
					&(i) \, u_s(x_1,x_2,0)= \omega_s(x_1,x_2)
					\,\text{ if }\  {\alpha_s}\in(0,1],s=1,2,\\
					&(ii)\,	u_s(x_1,x_2,0)= \omega_s(x_1,x_2),\, \,
					{\dfrac{\partial u_s}{\partial t}}\big{|}_{t=0}=\vartheta_s(x_1,x_2)
					\,	\text{if}\  {\alpha_s}\in(1,2],s=1,2.
				\end{aligned}
		\end{eqnarray}}
	The above-coupled equations   \eqref{poly-polyeq} admit the invariant product linear space
		$$
		{{	\mathbf{\hat{W}}^{1}_8}}=\text{Span}\left\{1,x_1,x_2,x_1x_2
		\right\} \times\text{Span}\left\{1,x_1,x_2,x_1x_2\right\},
		$$
		which is discussed in  case 39 of Table 5.
			Then, there exists a generalized separable  {analytical} solution for the coupled equations \eqref{poly-polyeq} in the form	
		\begin{equation}\label{solutionformeg1}
			u_s(x_1,x_2,t)=
			\delta_{s1}(t)+\delta_{s2}(t)x_1+\delta_{s3}(t)x_2+\delta_{s4}(t)x_1x_2,s=1,2.
		\end{equation}
		Now, we want to find the unknown functions
		$ \delta_{si}(t),i=1,2,3,4, $
		of the above {analytical} solution
		$ u_s=	u_s(x_1,x_2,t),s=1,2. $
		Thus, we substitute \eqref{solutionformeg1} into the above-coupled equations  \eqref{poly-polyeq}.
		Hence, we obtain the following system of FODEs:
		{	\begin{eqnarray}
				\begin{aligned}
					\label{poly-polyfodes}
					&{	\mathbf{D}^{\alpha_s}_t \delta_{s1}(t)}=
					g_{s0}\delta_{22}(t)+h_{s0}\delta_{23}(t)+\eta_{s0}\delta_{12}(t)+d_{s0}\delta_{13}(t),
					\\&
					{	\mathbf{D}^{\alpha_s}_t \delta_{s2}(t)}=h_{s0}\delta_{24}(t)+d_{s0}
					\delta_{14}(t),
					\\
					&{	\mathbf{D}^{\alpha_s}_t \delta_{s3}(t)}=
					g_{s0}\delta_{24}(t)+\eta_{s0}\delta_{14}(t),	\text{ and } &
				\\&	{\mathbf{D}^{\alpha_s}_t \delta_{s4}(t)}=0,s=1,2 ,
				\end{aligned}
		\end{eqnarray}}
		where
		{$ 	\mathbf{D}^{\alpha_s}_t(\cdot)=\dfrac{d^{\alpha_s}}{dt^{\alpha_s}}  (\cdot)$}
		represents the Caputo fractional-order derivative \eqref{caputo} of order {${\alpha_s}\in(0,2],s=1,2,$}  with respect to $t.$
		
		Next, we use the Laplace transformation method to solve the above-reduced system of FODEs \eqref{poly-polyfodes}.
	{In	\cite{Podlubny1999}, the Laplace transformation of the Caputo  derivative with fractional-order $\alpha>0$ is obtained as follows:}	$$\mathit{L}(\mathbf{D}^\alpha_t\delta(t);s)= s^\alpha\mathit{L}(\delta(t);s)-\sum\limits_{q=0}^{n-1}s^{\alpha-(q+1)} \dfrac{d^q\delta(t)}{dt^q}\Big{|}_{t=0},\,\text{Re}(s)>0,$$	
		where  $\alpha\in(n-1,n],n\in\mathbb{N}.$
		Here, we apply the Laplace transformation on the last equation of the obtained system of FODEs \eqref{poly-polyfodes}, and then taking the inverse Laplace transformation of the resultant equation, we have
		{	\begin{eqnarray}
				\begin{aligned}\label{deltai4}
					\delta_{s4}(t)=& \left\{
					\begin{array}{ll}
						A_{s4},\text{ if } {\alpha_s}\in(0,1],
						\\
						A_{s4}+B_{s4}t,\text{ if } {\alpha_s}\in(1,2],
					\end{array}
					\right.
				\end{aligned}
		\end{eqnarray}}
		where $A_{s4},B_{s4}\in\mathbb{R},s=1,2.$
		
	Now, we substitute the above-obtained function $\delta_{s4},s=1,2,$ in the second  and third equations of the given system of FODEs \eqref{poly-polyfodes}.  Thus,  {using the Laplace  and the inverse Laplace  transformations, } we get
		\begin{eqnarray}\label{delta12}
			\delta_{12}(t)=& \left\{
			\begin{array}{ll}
				A_{12}+\dfrac{\mu_{1}t^{\alpha_1}}{\Gamma(\alpha_1+1)},\text{ if } {\alpha_s}\in(0,1],s=1,2,
				\\
				A_{12}+B_{12}t+\Big( \dfrac{\mu_{1}}{\Gamma(\alpha_1+1)}+\dfrac{d_{10}B_{14}t}{\Gamma(\alpha_1+2)}\Big) t^{\alpha_1},\text{ if } {\alpha_1}\in(1,2]\,\&\,{\alpha_2}\in(0,1],
				\\
				A_{12}+\Big( \dfrac{\mu_{1}}{\Gamma(\alpha_1+1)}+\dfrac{h_{10}B_{24}t}{\Gamma(\alpha_1+2)}\Big) t^{\alpha_1},\text{ if } {\alpha_1}\in(0,1]\,\&\,{\alpha_2}\in(1,2],
				\\
				A_{12}+B_{12}t+\Big( \dfrac{\mu_{1}}{\Gamma(\alpha_1+1)}+
				\dfrac{\nu_{1}t}{\Gamma(\alpha_1+2)}
				\Big) t^{\alpha_1},\text{ if } {\alpha_s}\in(1,2],s=1,2,
			\end{array}\right.
	\\
		\label{delta22}
			\delta_{22}(t)=& \left\{
			\begin{array}{ll}
				A_{22}+\dfrac{\mu_{2}t^{\alpha_2}}{\Gamma(\alpha_2+1)},\text{ if } {\alpha_s}\in(0,1],s=1,2,
				\\
				A_{22}+\Big( \dfrac{\mu_{2}}{\Gamma(\alpha_2+1)}+\dfrac{d_{20}B_{14}t}{\Gamma(\alpha_2+2)}\Big) t^{\alpha_2},\text{ if } {\alpha_1}\in(1,2]\,\&\,{\alpha_2}\in(0,1],
				\\
				A_{22}+B_{22}t+\Big( \dfrac{\mu_{2}}{\Gamma(\alpha_2+1)}+\dfrac{h_{20}B_{24}t}{\Gamma(\alpha_2+2)}\Big) t^{\alpha_2},\text{ if } {\alpha_1}\in(0,1]\,\&\,{\alpha_2}\in(1,2],
				\\
				A_{22}+B_{22}t+\Big( \dfrac{\mu_{2}}{\Gamma(\alpha_2+1)}+
				\dfrac{\nu_{2}t}{\Gamma(\alpha_2+2)}
				\Big) t^{\alpha_2},\text{ if } \alpha_s\in(1,2],s=1,2,
			\end{array}\right.
		\\
			\label{delta13}	\delta_{13}(t)=& \left\{
			\begin{array}{ll}
				A_{13}+\dfrac{\gamma_{1}t^{\alpha_1}}{\Gamma({\alpha_1}+1)},\text{ if } {\alpha_s}\in(0,1],s=1,2,
				\\
				A_{13}+B_{13}t+\Big( \dfrac{\gamma_{1}}{\Gamma({\alpha_1}+1)} + \dfrac{\eta_{10}B_{14}t}{\Gamma({\alpha_1}+2)} \Big) t^{\alpha_1},\text{ if } {\alpha_1}\in(1,2]\,\&\,{\alpha_2}\in(0,1],	\\
				A_{13}+\Big( \dfrac{\gamma_{1}}{\Gamma({\alpha_1}+1)} + \dfrac{g_{10}B_{24}t}{\Gamma({\alpha_1}+2)} \Big) t^{\alpha_1},\text{ if } {\alpha_1}\in(0,1]\,\&\,{\alpha_2}\in(1,2],
				\\
				A_{13}+B_{13}t+\Big( \dfrac{\gamma_{1}}{\Gamma({\alpha_1}+1)}+
				\dfrac{\kappa_{1}t}{\Gamma({\alpha_1}+2)}
				\Big) t^{\alpha_1 },\text{ if } {\alpha_s}\in(1,2],s=1,2,\text{ and }
			\end{array}\right.
			\\
			\label{delta23}		\delta_{23}(t)=& \left\{
			\begin{array}{ll}
				A_{23}+\dfrac{\gamma_{2}t^{\alpha_2}}{\Gamma({\alpha_2}+1)},\text{ if } {\alpha_s}\in(0,1],s=1,2,
				\\
				A_{23}+\Big( \dfrac{\gamma_{2}}{\Gamma({\alpha_2}+1)} + \dfrac{\eta_{20}B_{14}t}{\Gamma({\alpha_2}+2)} \Big) t^{\alpha_2},\text{ if } {\alpha_1}\in(1,2]\,\&\,{\alpha_2}\in(0,1],	\\
				A_{23}+B_{23}t+\Big( \dfrac{\gamma_{2}}{\Gamma({\alpha_2}+1)} + \dfrac{g_{20}B_{24}t}{\Gamma({\alpha_2}+2)} \Big) t^{\alpha_2},\text{ if } {\alpha_1}\in(0,1]\,\&\,{\alpha_2}\in(1,2],
				\\
				A_{23}+B_{23}t+\Big( \dfrac{\gamma_{2}}{\Gamma({\alpha_2}+1)}+
				\dfrac{\kappa_{2}t}{\Gamma({\alpha_2}+2)}
				\Big) t^{\alpha_2},\text{ if } {\alpha_s}\in(1,2],s=1,2,
			\end{array}\right.\quad
		\end{eqnarray}
		{where }	$\gamma_{s}=g_{s0}A_{24}+\eta_{s0}A_{14},
		$
		$\kappa_{s}=g_{s0}B_{24}+\eta_{s0}B_{14},$
		$
		\mu_{s}=h_{s0}A_{24}+d_{s0}A_{14},
		$
		$\nu_{s}=h_{s0}B_{24}+d_{s0}B_{14},$   { and } $  A_{si},B_{si}\in\mathbb{R}, i=2,3,4,s=1,2.$
	
			Finally, we substitute the equations \eqref{deltai4}-\eqref{delta23} in the first equation of the obtained system \eqref{poly-polyfodes}, and then using the Laplace and inverse Laplace transformations, we have
		{\begin{eqnarray}
					\label{delta21}		\delta_{21}(t)=& \left\{
					\begin{array}{ll}
						A_{21}+\Big( \dfrac{l_{2}}{\Gamma({\alpha_2}+1)} + \sum\limits_{i=1}^2 \dfrac{\varrho_{2i}t^{\alpha_i}}{\Gamma({\alpha_2}+{\alpha_i}+1)} \Big) t^{\alpha_2},\text{ if } \alpha_s\in(0,1],s=1,2,
						\\
A_{21}+\Big( \dfrac{l_{2}}{\Gamma({\alpha_2}+1)}   +\dfrac{(\eta_{20}B_{12}+d_{20}B_{13})t}{\Gamma({{\alpha_2}+2)}}+\dfrac{(g_{20}d_{20} +h_{20}\eta_{20})B_{14}t^{\alpha_2+1}}{\Gamma(2{\alpha_2}+2)}
						\\+\sum\limits_{i=1}^2 \dfrac{\varrho_{2i}t^{\alpha_i}}{\Gamma({\alpha_2}+{\alpha_i}+1)}	+ \dfrac{(d_{20}\eta_{10}+d_{10}\eta_{20})B_{14}t^{\alpha_1+1}}{\Gamma({\alpha_1}+{\alpha_2}+2)}
						\Big) t^{\alpha_2},
					\\	\text{ if } \alpha_1\in(1,2]\ \&\ \alpha_2\in(0,1],
						\\
	A_{21}+B_{21}t+\Big( \dfrac{l_{2}}{\Gamma({\alpha_2}+1)} +  \sum\limits_{i=1}^2 \dfrac{\varrho_{2i}t^{\alpha_i}}{\Gamma({\alpha_2}+{\alpha_i}+1)}
						+ \dfrac{(g_{20}B_{22}+h_{20}B_{23})t}{\Gamma{({\alpha_2}+2)}}
					\\	+\dfrac{2g_{20}h_{20}B_{24}t^{\alpha_2+1}}{\Gamma(2{\alpha_2}+2)}
+\dfrac{B_{24}(\eta_{20}h_{10}+d_{20}g_{10})t^{\alpha_1+1}}{\Gamma(\alpha_1+{\alpha_2}+2)}						\Big) t^{\alpha_2},
			\\			\text{ if } \alpha_1\in(0,1]\,\&\,\alpha_2\in(1,2],
						\\
A_{21}+B_{21}t+	\dfrac{l_{2} t^{\alpha_2}}{\Gamma({\alpha_2}+1)}
+ \dfrac{\varsigma_{2} t^{\alpha_2+1}}{\Gamma({{\alpha_2}+2)}}
+\sum\limits_{i=1}^2 \Big(\dfrac{\varrho_{2i}t^{\alpha_i}}{\Gamma({\alpha_i}+{\alpha_2}+1)}
					\\	+
						\dfrac{\lambda_{2i}t^{\alpha_i+1}}{\Gamma({\alpha_i}+{\alpha_2}+2)}
						\Big) t^{\alpha_2},
							\text{ if } \alpha_s\in(1,2],s=1,2,
					\end{array}\right.
\\
						\label{delta11}	\delta_{11}(t)=& \left\{
					\begin{array}{ll}
						A_{11}+\Big( \dfrac{l_{1}}{\Gamma({\alpha_1}+1)} + \sum\limits_{i=1}^2\dfrac{\varrho_{1i}t^{\alpha_i}}{\Gamma({\alpha_1}+{\alpha_i}+1)}  \Big) t^{\alpha_1},\text{ if } \alpha_s\in(0,1],s=1,2,
						\\
	A_{11}+B_{11}t+\Big( \dfrac{l_{1}}{\Gamma({\alpha_1}+1)}  +\dfrac{(\eta_{10}B_{12}+d_{10}B_{13})t}{\Gamma({{\alpha_1}+2)}}
+\sum\limits_{i=1}^2\dfrac{\varrho_{1i}t^{\alpha_i}}{\Gamma({\alpha_1}+{\alpha_i}+1)}
\\ + \dfrac{2d_{10}\eta_{10}B_{14}t^{\alpha_1+1}}{\Gamma(2{\alpha_1}+2)}
+ \dfrac{(g_{10}d_{20}+\eta_{20}h_{10})B_{14}t^{\alpha_2+1}}{\Gamma({\alpha_1}+\alpha_2+2)}						\Big) t^{\alpha_1},
\\										\text{ if } \alpha_1\in(1,2]\,\&\,\alpha_2\in(0,1],
						\\
						A_{11}+ \Big( \dfrac{B_{24}(g_{20}h_{10}+g_{10}h_{20})t^{\alpha_2+1}}{\Gamma({\alpha_1+{\alpha_2}+2)}}	+ \dfrac{B_{24}(\eta_{10}h_{10}+d_{10}g_{10})t^{\alpha_1+1}}{\Gamma({2{\alpha_1}+2)}}
						\\ + \dfrac{l_{1}}{\Gamma({\alpha_1}+1)}
						+ \sum\limits_{i=1}^2\dfrac{\varrho_{1i}t^{\alpha_i}}{\Gamma({\alpha_1}+{\alpha_i}+1)}	+			 \dfrac{(h_{10}B_{23}+g_{10}B_{22})t}{\Gamma({\alpha_1}+2)}\Big) t^{\alpha_1}
						,
					\\	\text{if} \ \alpha_1\in(0,1]\ \&\ \alpha_2\in(1,2],
						\\
						A_{11}+B_{11}t+ \dfrac{l_{1}t^{\alpha_1}}{\Gamma({\alpha_2}+1)}	+ \dfrac{\varsigma_{1}t^{\alpha_1+1}}{\Gamma({{\alpha_1}+2)}} +\sum\limits_{i=1}^2\Big(  \dfrac{\varrho_{1i}t^{\alpha_i}}{\Gamma({\alpha_1}+{\alpha_i}+1)}
					\qquad\qquad\quad	\\+
						\dfrac{\lambda_{1i}t^{\alpha_i+1}}{\Gamma({\alpha_1}+{\alpha_i}+2)}
						\Big) t^{\alpha_1},
											\text{ if } \alpha_s\in(1,2],s=1,2,
					\end{array}\right.
				\end{eqnarray}
				where $l_{s}=g_{s0}A_{22}+h_{s0}A_{23}+d_{s0}A_{13}+\eta_{s0}A_{12},$
				$
				\varrho_{s1}=d_{s0}\gamma_{1}+\eta_{s0}\mu_{1},	$
				$	\varrho_{s2}=g_{s0}\mu_{2}+h_{s0}\gamma_{2},$
				$	\varsigma_{s}=g_{s0}B_{22}+h_{s0}B_{23}+d_{s0}B_{13}+\eta_{s0}B_{12},
				$
				$\lambda_{s1}=d_{s0}\kappa_{1}+\eta_{s0}\nu_{1},\, \lambda_{s2}=g_{s0}\nu_{2}+h_{s0}\kappa_{2},
				$ and
				$	A_{si},B_{si}\in{\mathbb{R}},s=1,2,i=1,2,3,4.
				$}
			
			Thus, the obtained  generalized separable {analytical} solution of the given system \eqref{poly-polyeq}  is $\mathbf{U}=(u_1,u_2)$ with components $u_s=u_s(x_1,x_2,t),s=1,2,$  as follows;
			{	\begin{eqnarray}\label{exactsolution1}
					\begin{aligned}
						&	u_1(x_1,x_2,t)=\delta_{11}(t)+\delta_{12}(t)x_1+ [\delta_{13}(t)+\delta_{14}(t)x_1]x_2,
						\text{ and}
						\\&
						u_2(x_1,x_2,t)=\delta_{21}(t)+\delta_{22}(t)x_1+ [\delta_{23}(t)+\delta_{24}(t)x_1]x_2,
					\end{aligned}
			\end{eqnarray}}
			where the functions $\delta_{si},s=1,2,i=1,2,3,4,$ are obtained in \eqref{deltai4}-\eqref{delta11}.
			We observe that the above {analytical} solution \eqref{exactsolution1} of the given coupled equations \eqref{poly-polyeq} satisfies the initial conditions \eqref{icpoly-polyeq} with
			$$\begin{aligned}
				&\omega_s(x_1,x_2)=A_{s1}+A_{s2}x_1+A_{s3}x_2+A_{s4}x_1x_2,
				\text{ and}
			\\	&
				\vartheta_s(x_1,x_2)=B_{s1}+B_{s2}x_1+B_{s3}x_2+B_{s4}x_1x_2,s=1,2.
			\end{aligned}
			$$
		{
		Here, note that whenever $\alpha_s=1$ and $\alpha_s=2$, the obtained solutions \eqref{exactsolution1} coincide with integer-order cases of the discussed coupled system \eqref{poly-polyeq}. 
							Also, we wish to point out that the given  equations \eqref{2+1diffzhang-reduced} are a particular case of the considered coupled equations   \eqref{poly-polyeq}  
				when the parameters are considered as $ q_{s1}=c_{s1}=\eta_{s0}=g_{s0}=d_{s0}=h_{s0}=0,s=1,2, $ $\kappa_{10}=\kappa_{1}\tau_{0},q_{10}=\kappa_{1}\tau_{1},\kappa_{20}=\kappa_{2}\tau_{2}, $ and $q_{20}=\kappa_{2}\tau_{0}.$ For the above-coupled equations \eqref{poly-polyeq}, we have obtained a generalized separable {analytical}  solution   \eqref{exactsolution1}  systematically using the developed invariant subspace method for the first time. However, the approximate solution of the above system \eqref{2+1diffzhang} can be expressed in terms of Mittag-Leffler functions, which are studied through  the numerical inversion of the Laplace transformation  technique \cite{zhang2019}.} We note that the
				obtained {analytical} solution  \eqref{exactsolution1}   of the {considered} coupled equations \eqref{poly-polyeq} is valid for all $\alpha_s, \alpha_s\in(0,2], s=1,2.$
				%
				%
				%
				%
				%
				%
				%
				%
				%
				%
				%
				%
				%
				%
				The two and three-dimensional  graphical representations of the obtained  {analytical} solution \eqref{exactsolution1} of the given  coupled equations \eqref{poly-polyeq} for arbitrary values of $\alpha_s,$ where $\alpha_s\in(0,2],s=1,2,$ are shown from Figure \ref{figa1} to Figure \ref{figc2}.
				In these Figures \ref{figa1}-\ref{figc2},  we considered the constants in the obtained {analytical} solution \eqref{exactsolution1} as $A_{11}=3,
				A_{21}=A_{14}=-h_{10}=d_{10}=A_{22}=A_{24}=B_{12}=B_{22}=2, A_{13}=\eta_{10}=-A_{21}=A_{23}=-d_{20}=
				-B_{13}=B_{14}=-B_{21}=B_{23}=-B_{24}=1,
				g_{10}=50,$
				and $B_{11}=5.$
\begin{figure}
	\begin{subfigure}{0.48\textwidth}
		\includegraphics[width=6.75cm, height=5cm]{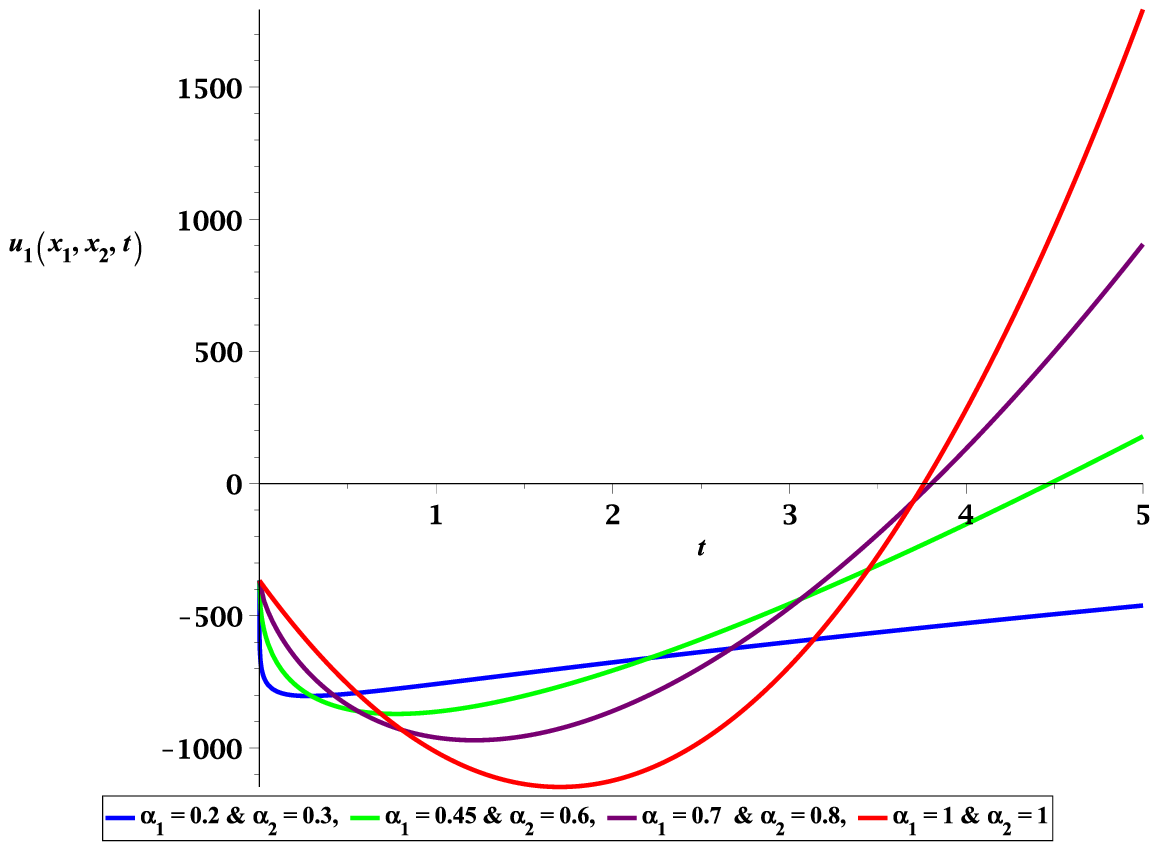}
		\caption{}
	\end{subfigure}
	\begin{subfigure}{0.47\textwidth}
		\includegraphics[width=6.75cm, height=5cm]{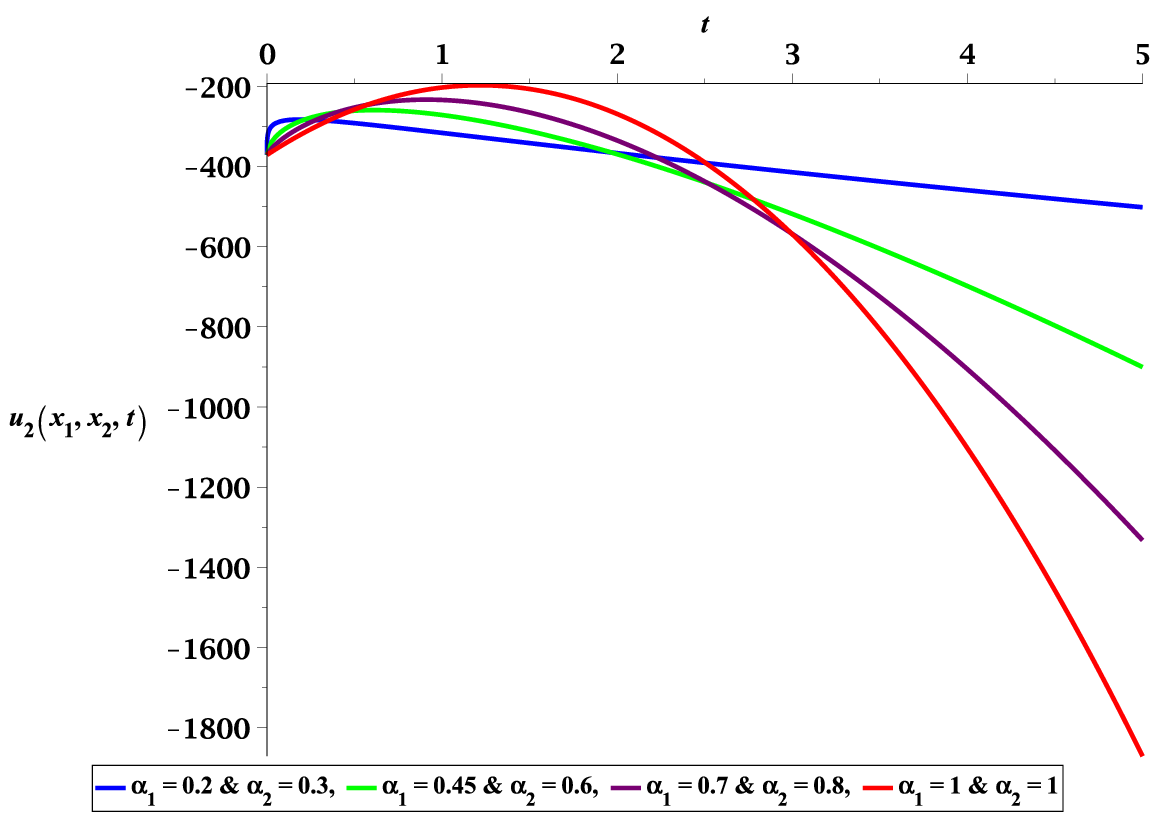}
		\caption{}
	\end{subfigure}
	\caption{{ 2D graphical representations of the {analytical} solution $ u_1(x_1,x_2,t)$ and $ u_2(x_1,x_2,t) $ given in \eqref{exactsolution1}  for various values of  $\alpha_s\in(0,1],s=1,2,$ $\eta_{20}=1,h_{20}=4,x_1=20,$ and $g_{20}=x_2=-10.$}}
	\label{figa1}
\end{figure}
\begin{figure}
	\begin{center}
		\begin{subfigure}{0.5\textwidth}
		\includegraphics[width=8.75cm, height=6.2cm]{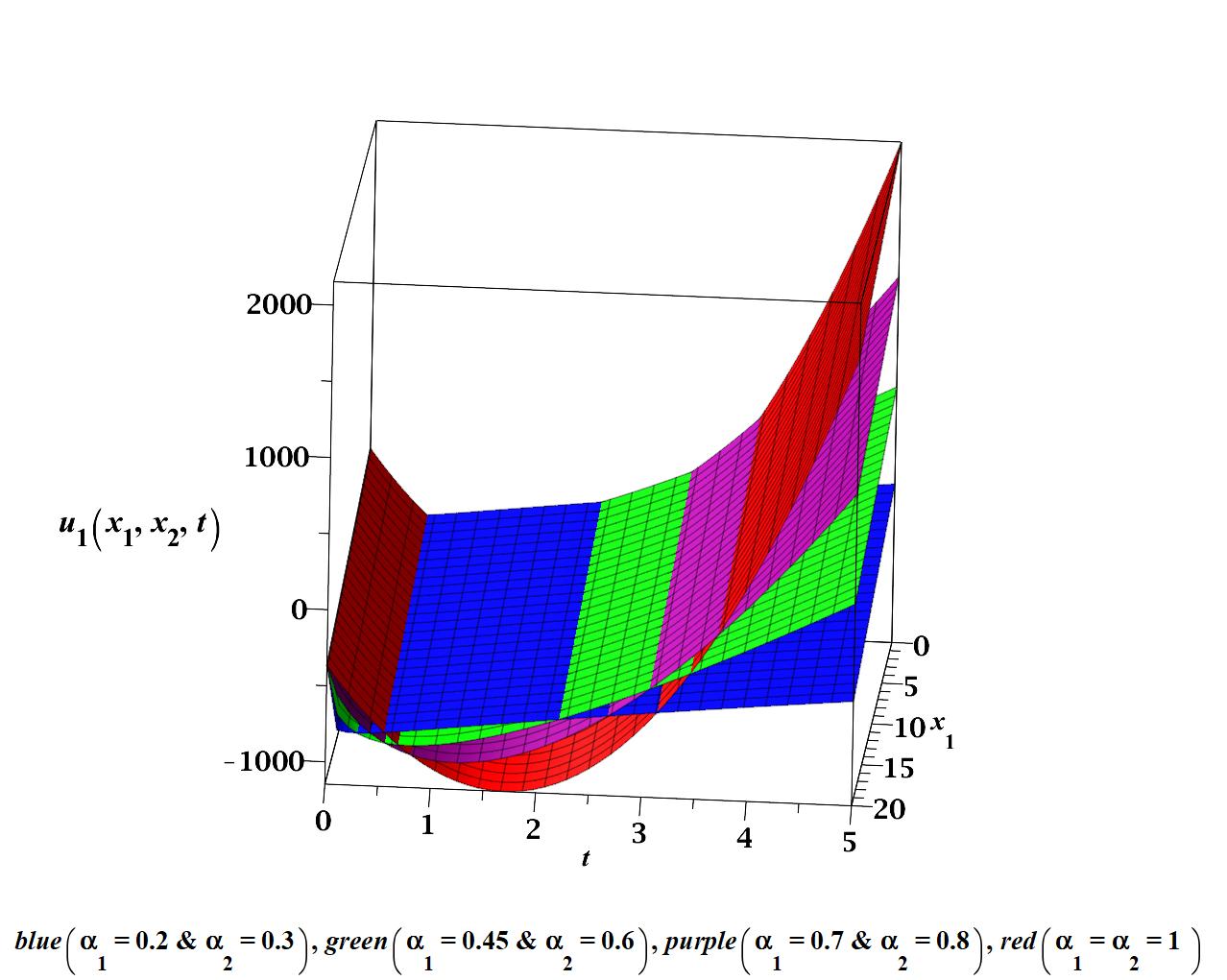}
		\caption{}
	\end{subfigure}
	\begin{subfigure}{0.5\textwidth}
		\includegraphics[width=9cm, height=6.2cm]{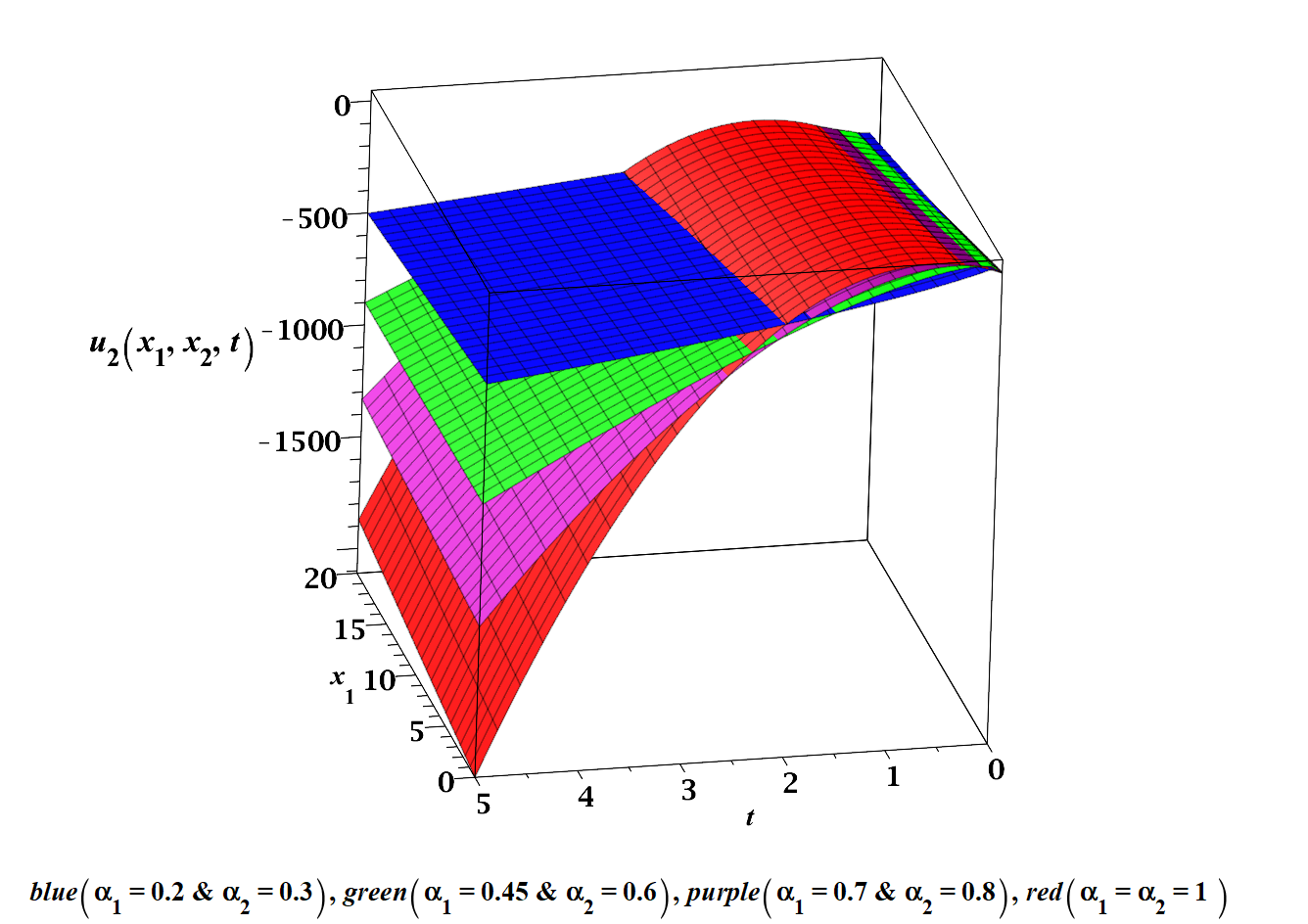}
		\caption{}
	\end{subfigure}
	\end{center}
	\caption{{3D graphical representations of the {analytical} solution $ u_1(x_1,x_2,t)$ and $ u_2(x_1,x_2,t) $ given in \eqref{exactsolution1} for various values of  $\alpha_s\in(0,1],s=1,2,$ $\eta_{20}=1,h_{20}=4,$ and $g_{20}=x_2=-10.$}}
	\label{figa2}
\end{figure}
\begin{figure}
	\begin{subfigure}{0.48\textwidth}\includegraphics[width=6.45cm, height=5.7cm]{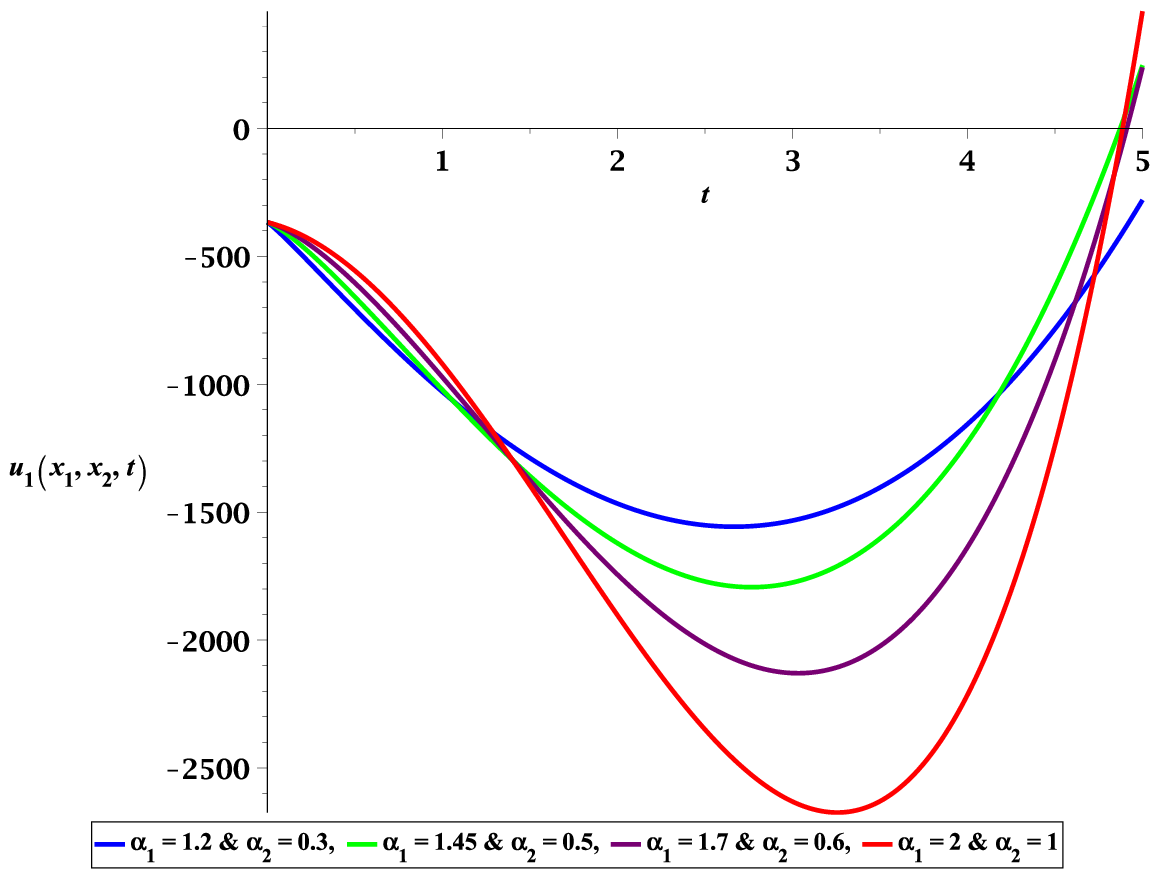}
		\caption{}
	\end{subfigure}
	\begin{subfigure}{0.48\textwidth}
		\includegraphics[width=6.45cm, height=5.7cm]{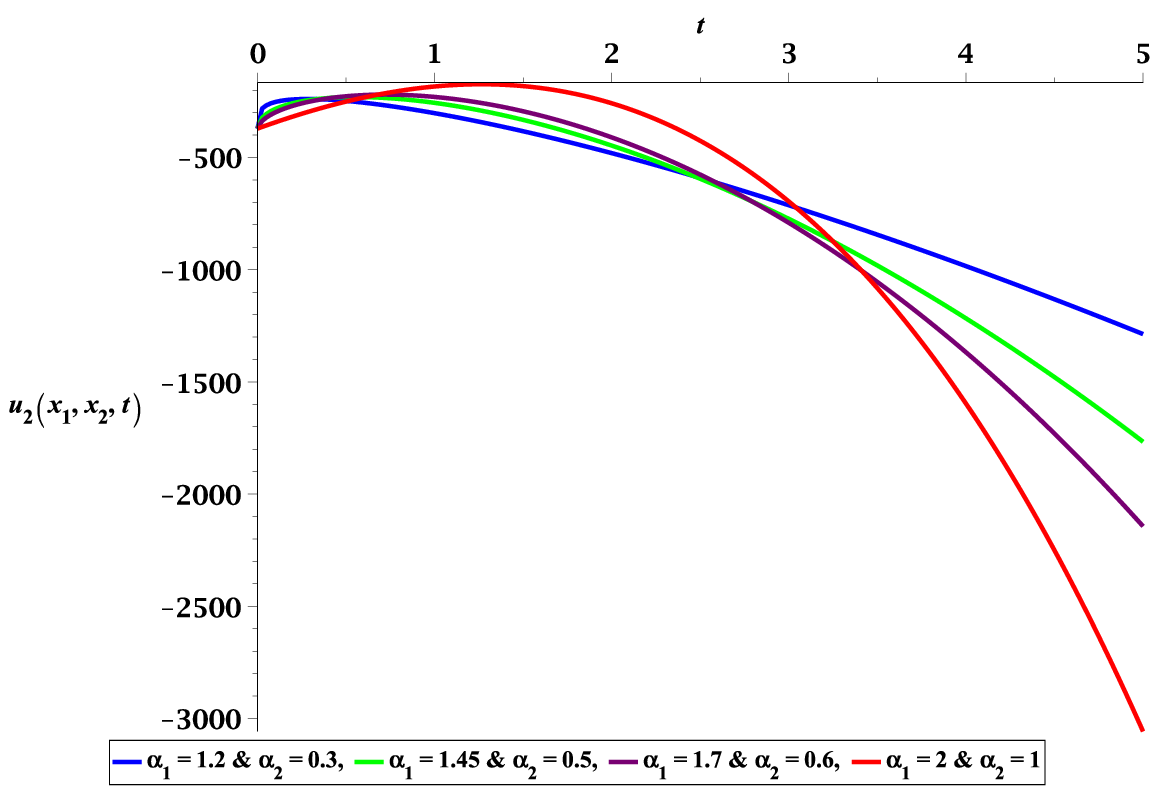}
		\caption{}
	\end{subfigure}
	\caption{{2D graphical representations of the {analytical} solution $ u_1(x_1,x_2,t)$ and $ u_2(x_1,x_2,t) $ given in \eqref{exactsolution1}  for various values of  $\alpha_1\in(1,2]$, $\alpha_2\in(0,1]$, $\eta_{20}=1,h_{20}=4,x_1=20,$ and $g_{20}=x_2=-10.$}}
	\label{figb1}
\end{figure}
\begin{figure}
\begin{center}
		\begin{subfigure}{0.6\textwidth}\includegraphics[width=8.8cm, height=6.2cm]{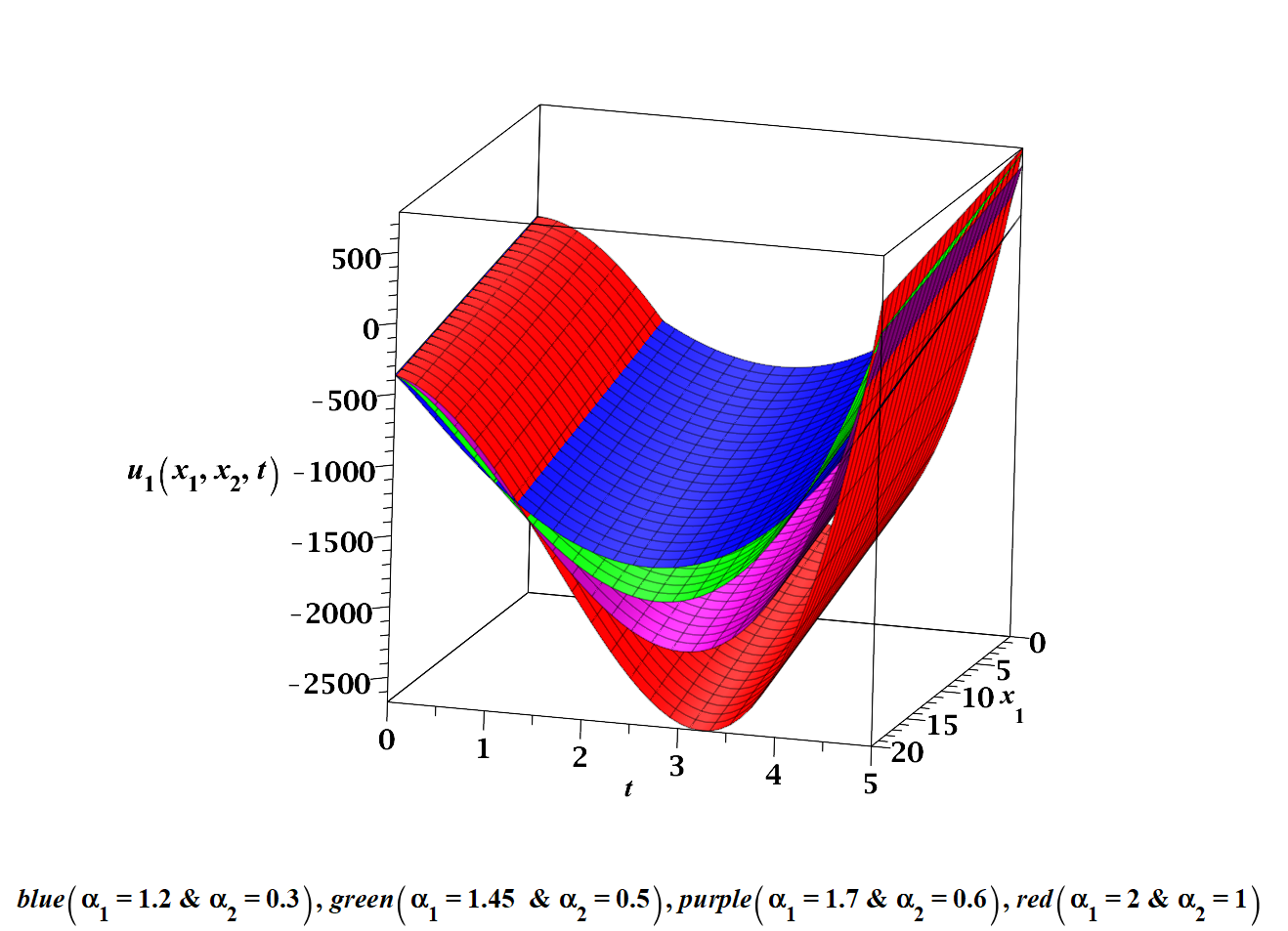}
		\caption{}
	\end{subfigure}
	\begin{subfigure}{0.6\textwidth}
		\includegraphics[width=8.8cm, height=6.2cm]{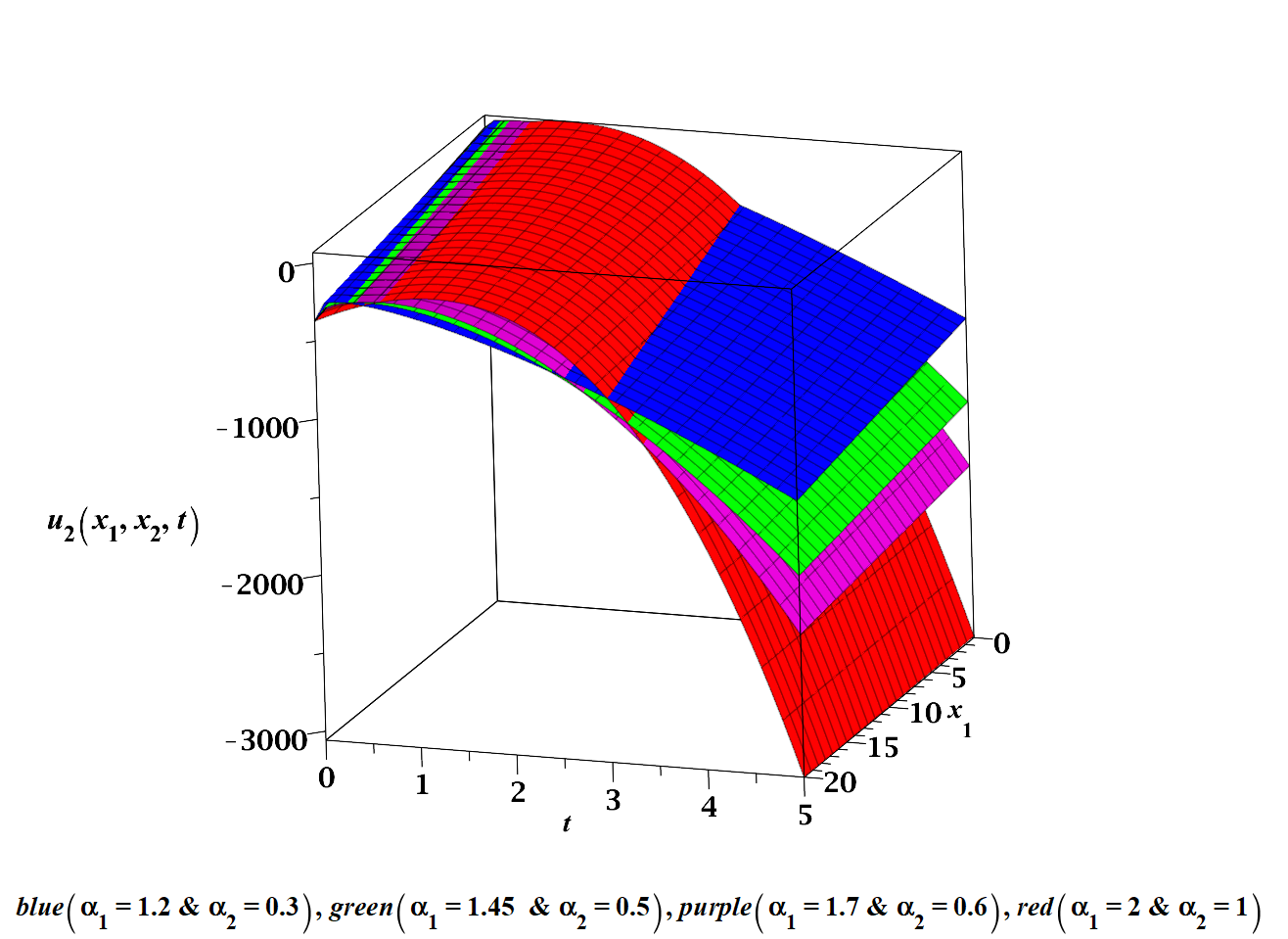}
		\caption{}
	\end{subfigure}
\end{center}
	\caption{{3D graphical representations of the {analytical} solution  $ u_1(x_1,x_2,t)$ and $ u_2(x_1,x_2,t) $ given in \eqref{exactsolution1} for various values of  $\alpha_1\in(1,2]$, $\alpha_2\in(0,1]$, $\eta_{20}=1,h_{20}=4,$ and $g_{20}=x_2=-10.$}}
	\label{figb2}
\end{figure}
\begin{figure}
	\begin{subfigure}{0.48\textwidth}\includegraphics[width=6.45cm, height=5cm]{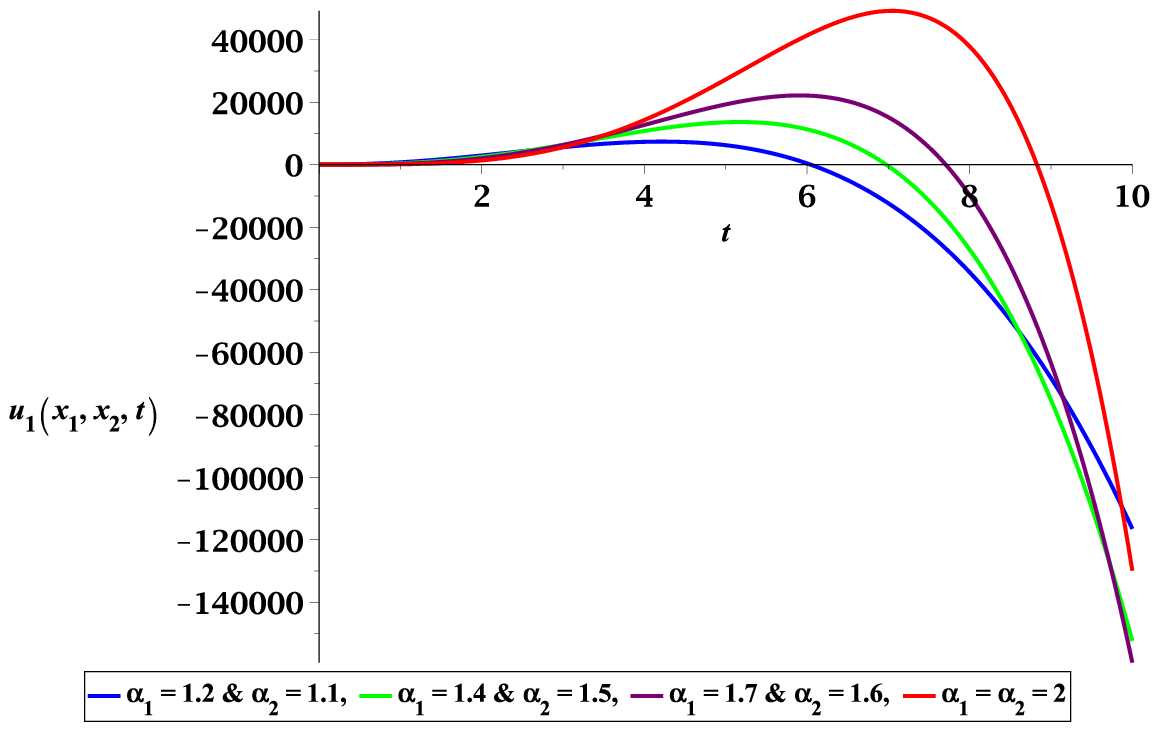}
		\caption{}
	\end{subfigure}
	\begin{subfigure}{0.48\textwidth}
		\includegraphics[width=6.45cm, height=5cm]{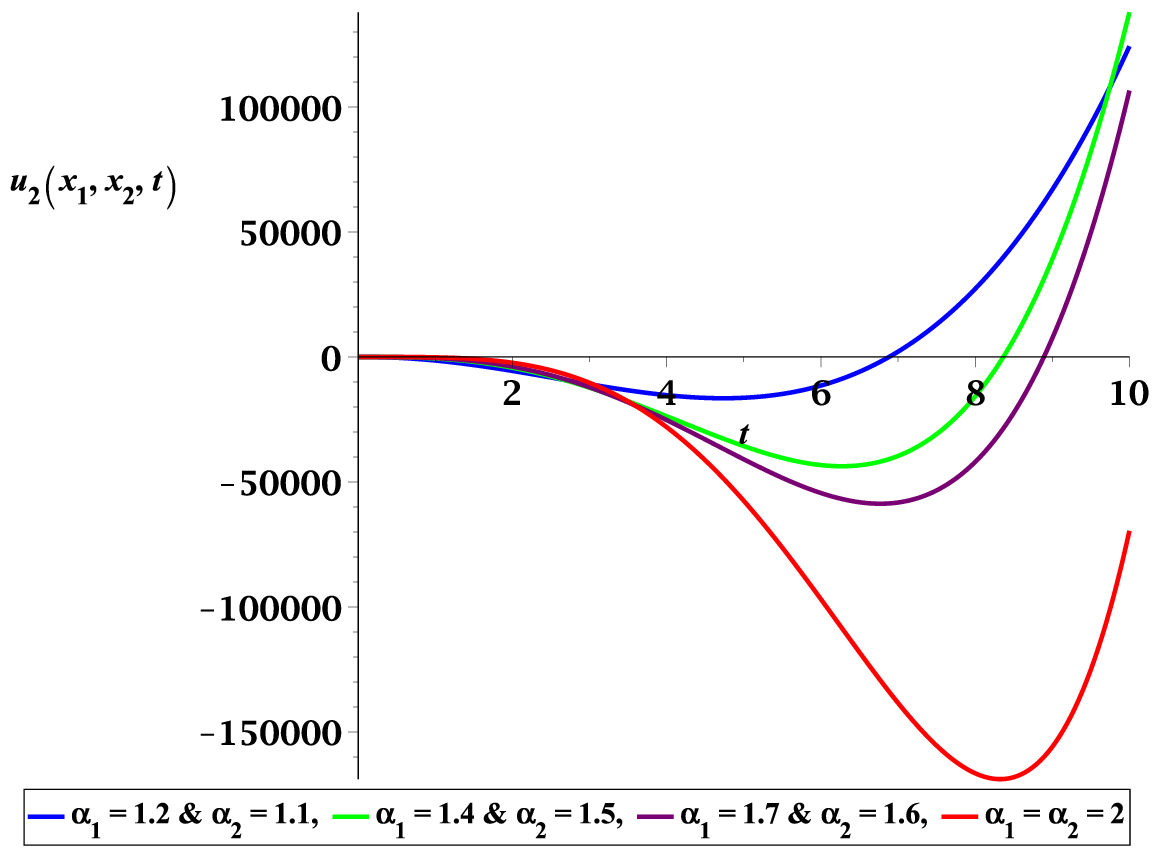}
		\caption{}
	\end{subfigure}
	\caption{{2D graphical representations of the {analytical} solution $ u_1(x_1,x_2,t)$ and $ u_2(x_1,x_2,t) $ given in \eqref{exactsolution1} for various values of  $\alpha_s\in(1,2]$, $s=1,2$, $\eta_{20}=-1,h_{20}=20,g_{20}=-50, x_1=20,$ and $x_2=-10.$}}
	\label{figc1}
\end{figure}
\begin{figure}
\begin{center}
		\begin{subfigure}{0.5\textwidth}\includegraphics[width=9.79cm, height=6.5cm]{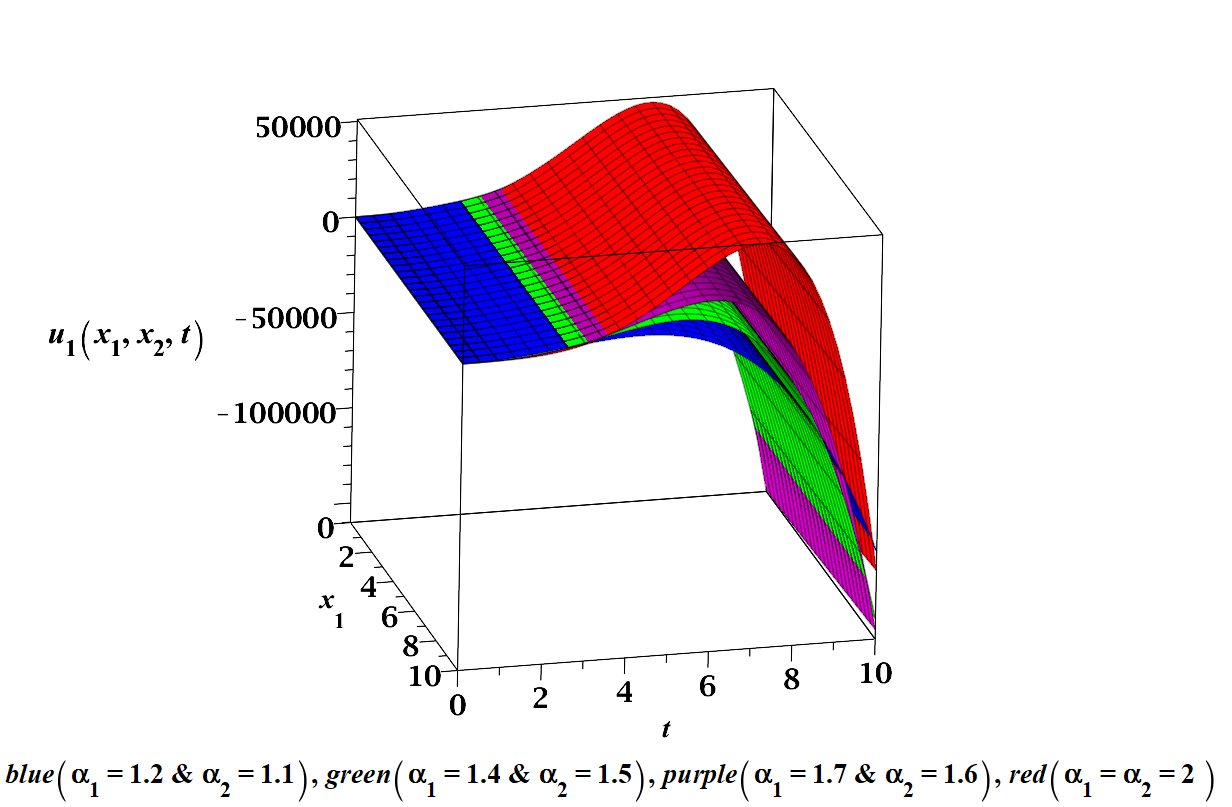}
		\caption{}
	\end{subfigure}
	\begin{subfigure}{0.5\textwidth}
		\includegraphics[width=9.5cm, height=6.75cm]{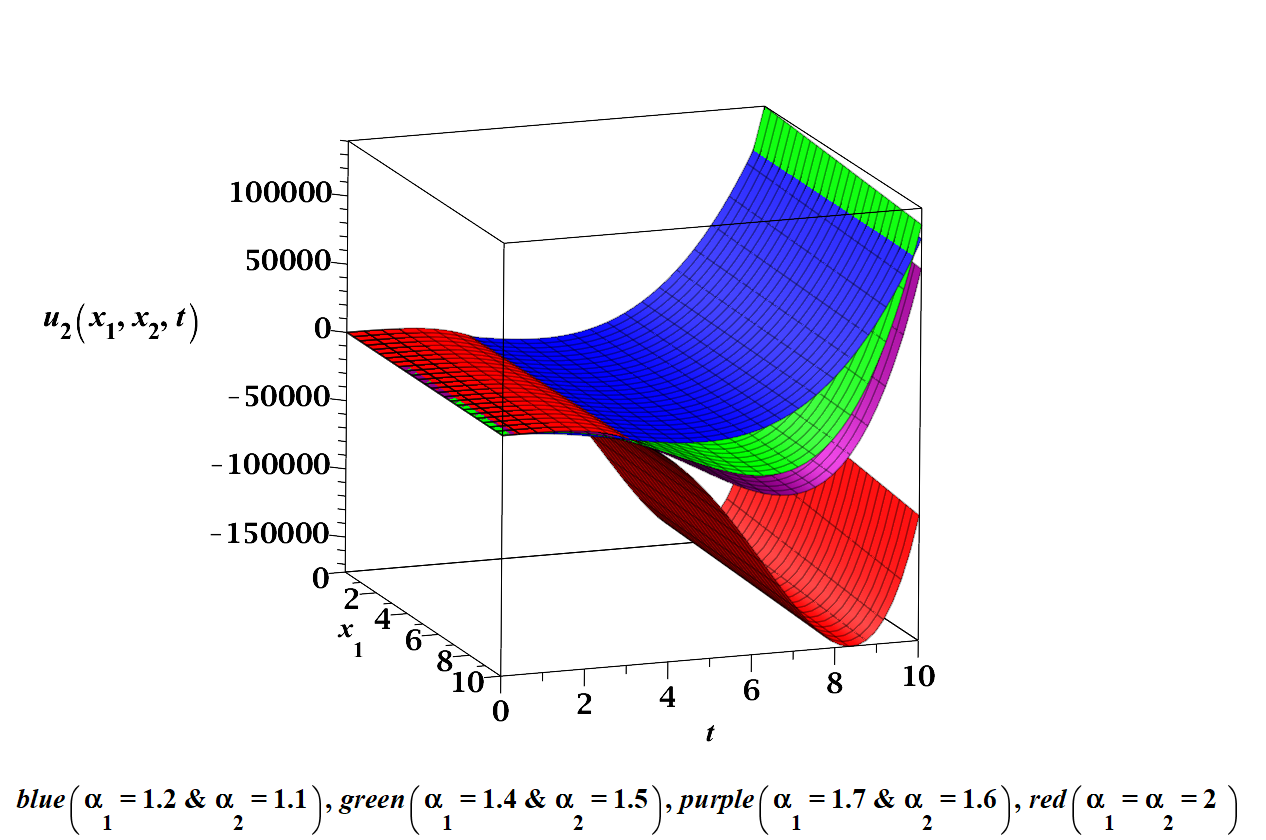}
		\caption{}
	\end{subfigure}
\end{center}
	\caption{{3D graphical representations of the {analytical} solution  $ u_1(x_1,x_2,t)$ and $ u_2(x_1,x_2,t) $ given in \eqref{exactsolution1} for various values of  $\alpha_s\in(1,2]$, $s=1,2$, $\eta_{20}=-1,h_{20}=20,g_{20}=-50,$ and $x_2=-10.$}}
	\label{figc2}
\end{figure}
				{
					\begin{note}
						From the graphical representations given in Figure \ref{figa1} to Figure \ref{figc2}, we observe that the behavior of  the corresponding physical process is highly dependent on the value  of the orders $\alpha_s,\alpha_s\in(0,2],s=1,2,$ of the Caputo fractional derivative that appears in the given coupled equations \eqref{poly-polyeq}. The nature of  the {analytical} solution \eqref{exactsolution1} of the {considered} coupled equations  \eqref{poly-polyeq} varies according to different values of $\alpha_s,\alpha_s\in(0,2],s=1,2.$ We observe that some important properties of the {analytical} solution \eqref{exactsolution1} of the considered coupled equations \eqref{poly-polyeq}  are given below.
						\begin{itemize}
							\item[(a)] From Figures \ref{figa1}-\ref{figa2}, for various values of $\alpha_s\in(0,1],s=1,2,$ the rate of change of  diffusion occurred  fast when $\alpha_1=\alpha_2=1 $ and slow if   $\alpha_1=0.2,\alpha_2=0.3.$   Note that for $\alpha_s<1,s=1,2,$ the slow-diffusion (sub-diffusion)  is visible for both the components $u_1$  and $u_2.$
							\item[(b) ] From the plotted Figures \ref{figb1}-\ref{figb2}, we {note} that the {analytical} solution \eqref{exactsolution1} of the considered coupled equations \eqref{poly-polyeq} follows the wave structure for the component  $u_1$ whereas the diffusion  behavior is prominent for the component  $u_2,$ for various values of $\alpha_1\in(1,2],$ and  $\alpha_2\in(0,1].$
							\item[(c) ] For the Caputo derivative of   order $\alpha_s\in(1,2],s=1,2,$ we observe a prominent wave behavior  for the components $u_1$ and $u_2$ as visible in Figures \ref{figc1}-\ref{figc2}.
						\end{itemize}
					\end{note}
				}
			\begin{remark}
		{In \cite{wu2019},   to describe  unsteady magnetohydrodynamic  flows, Wu et al.  have studied the integer-order {multi-component} $(2+1)$-dimensional coupled  diffusion-convection equations having the form of  \eqref{poly-polyeq} with parametric restrictions $\alpha_1=\alpha_2=1,$ and $ q_{s1}=q_{10}=c_{s1}=c_{10}=\eta_{10}=d_{10}=\kappa_{20}=p_{20}=g_{20}=h_{20}=0,s=1,2.$ They have obtained an {analytical} solution of the considered coupled equations  in terms of exponential and trigonometric functions in $t ,x_1,$ and $x_2.$ We note that  the given coupled  equations \eqref{poly-polyeq} with
					$ q_{s1}=q_{10}=c_{s1}=c_{10}=\eta_{10}=d_{10}=\kappa_{20}=p_{20}=g_{20}=h_{20}=0,s=1,2,$ {preserve  the invariant product subspace}
					\begin{eqnarray*}
						\begin{aligned}\label{invariant-comparison} {W}_8&=\text{Span}\{sin(a_1x_1)cos(a_2x_2),sin(a_1x_1)sin(a_2x_2),cos(a_1x_1)cos(a_2x_2),
							\\& cos(a_1x_1)sin(a_2x_2)\}
							 \times \text{Span}\{sin(a_1x_1)cos(a_2x_2),sin(a_1x_1)sin(a_2x_2),
							 \\&cos(a_1x_1)cos(a_2x_2),cos(a_1x_1)sin(a_2x_2)\},
						\end{aligned}
					\end{eqnarray*}
				where $a_1,a_2\in\mathbb{R}.$	Thus, using the above invariant product linear space ${W}_8,$   we can obtain a generalized separable {analytical} solution of the considered coupled equations  \eqref{poly-polyeq} with the above-mentioned  parametric restrictions
					as follows:
					\begin{eqnarray}\begin{aligned}\label{wu219soln}
							u_s(x_1,x_2,t)=&sin(a_1x_1)[\delta_{s1}(t)cos(a_2x_2)+\delta_{s2}(t)sin(a_2x_2)]
					\\&		+cos(a_1x_1)[\delta_{s3}(t)cos(a_2x_2)+\delta_{s4}(t)sin(a_2x_2)],s=1,2,
						\end{aligned}
					\end{eqnarray} where the unknown functions $\delta_{sj}(t),s=1,2,j=1,2,3,4,$ can be determined by following the same procedure of the  invariant subspace method as used in  Example \ref{eg1}.
					Note that the obtained {analytical} solution {\eqref{wu219soln}} is more general than the derived solution  in \cite{wu2019}.}
			\end{remark}
		\end{example}	
		%
		%
		%
		%
		%
		%
		%
		%
		%
		%
		%
		%
		%
		%
		\begin{example}\label{eg2}
	{	Now, we consider a coupled   nonlinear {TFDCWEs} in $(2+1)$-dimensions of the following  form}
			{	\begin{eqnarray}
					\begin{aligned}\label{poly-exp}
						\dfrac{\partial ^{\alpha_1} u_1}{\partial t^{\alpha_1}}
						=&
						{\frac {\partial }{\partial x_1}}
						\left[  \left(
						k_{{10}} -q_{{11}}u_{{2}}\right) {\frac {\partial u_{{1}} }{	\partial x_1}}
						+ \left(
						q_{{10}} +q_{{11}}u_{{1}}\right) {\frac {\partial u_{{2}} }{
								\partial x_1}}  \right]
						+{\frac {\partial }{	\partial x_2}}
						\left[\left( \dfrac{d_{12}}{a_{21}}u_{{2}}
						+p_{{10}} \right) {\frac {\partial u_{{1}}}{\partial x_2}}
							\right.\\& \left.
						+ \left(c_{{10}}- \dfrac{d_{12}}{a_{21}}u_{{1}} \right) {\frac {\partial u_{{2}}}{\partial x_2}}  \right]
							+ \eta_{{10}}  {\frac {\partial u_{{1}} }{\partial x_1}}
						+
						\left(d_{{12}}u_{{2}} +d_{{10}}
						\right){\frac {\partial u_{{1}} }{\partial x_2}}
						+
						h_{{10}}
						{\frac {\partial u_{{2}} }{\partial x_2}},
						\\
						%
						%
						%
						%
						%
						%
						\dfrac{\partial ^{\alpha_2} u_2}{\partial t^{\alpha_2}}
						=&
						{\frac {\partial }{\partial x_1}}
						\left[  \left(
						k_{{20}}-q_{{21}}u_{{2}} \right) {\frac {\partial u_{{1}} }{	\partial x_1}}
						+ \left( q_{{21}}u_{{1}}
						+q_{{20}} \right) {\frac {\partial u_{{2}} }{
								\partial x_1}}  \right]
						+c_{{20}}  {\frac {\partial^2 u_{{2}}}{
								\partial x_2^2}}
						+g_{{20}}{\frac {\partial u_{{2}}}{\partial x_1}}
						+h_{{20}}
						{\frac {\partial u_{{2}} }{\partial x_2}},
											\end{aligned}
			\end{eqnarray}}
		${\alpha_s}\in(0,2],s=1,2,$	with appropriate initial conditions
			{\begin{eqnarray}
					\begin{aligned}\label{icpoly-exp}
						&(i) \, u_s(x_1,x_2,0)= \omega_s(x_1,x_2)
						\,\text{ if }\  {\alpha_s}\in(0,1],s=1,2,\\
						&(ii)\,	u_s(x_1,x_2,0)= \omega_s(x_1,x_2),\, \,
						{\dfrac{\partial u_s}{\partial t}}\big{|}_{t=0}=\vartheta_s(x_1,x_2)
						\,	\text{if}\  {\alpha_s}\in(1,2],s=1,2.
					\end{aligned}
			\end{eqnarray}}
			The exponential-polynomial product linear space $${{	\mathbf{\hat{W}}^{1}_8}}=\text{Span}\left\{1,x_1,e^{-a_{21}x_2},x_1e^{-a_{21}x_2}\right\}\times\text{Span}\left\{1,x_1,x_2,x_1x_2\right\}$$ is an invariant product linear space for the given coupled equations \eqref{poly-exp}, which is discussed in  case 40 of Table 5.
			By proceeding with the similar technique employed in Example \ref{eg1}, we obtain the generalized separable {analytical}  solution for the {considered} coupled equations  \eqref{poly-exp} as follows:
			{	\begin{eqnarray}
					\begin{aligned}\label{solutionu1ofeg2}
						&	u_1(x_1,x_2,t)=\delta_{11}(t)+\delta_{12}(t)x_1+ [\delta_{13}(t)+\delta_{14}(t)x_1]e^{-a_{{21}}x_2},
						\text{ and}
						\\&
						u_2(x_1,x_2,t)=\delta_{21}(t)+\delta_{22}(t)x_1+ [\delta_{23}(t)+\delta_{24}(t)x_1]x_2,
					\end{aligned}
			\end{eqnarray}}
			where 	\begin{eqnarray*}				&	\begin{aligned} &	\delta_{11}(t)=	
			\left\{
			\begin{array}{ll}
				A_{11}+\left( \dfrac{l_{1}}{\Gamma({\alpha_1}+1)} +  \sum\limits_{i=1}^2\dfrac{\varrho_{1i}t^{\alpha_i}}{\Gamma({\alpha_1}+{\alpha_i}+1)} \right) t^{\alpha_1},\text{ if } \alpha_s\in(0,1],s=1,2,
				\\
				A_{11}+B_{11}t+\left( \dfrac{l_{1}}{\Gamma({\alpha_1}+1)} +
				\dfrac{\eta_{10}B_{12}t}{\Gamma({\alpha_1}+2)} +  \sum\limits_{i=1}^2\dfrac{\varrho_{1i}t^{\alpha_i}}{\Gamma({\alpha_1}+{\alpha_i}+1)} \right) t^{\alpha_1},
				\\	\text{ if } \alpha_1\in(1,2]\ \&\ \alpha_2\in(0,1],
				\\
				A_{11}+\Big[\dfrac{l_{1}}{\Gamma({\alpha_1}+1)}  +  \sum\limits_{i=1}^2\left( \dfrac{\varrho_{1i}}{\Gamma({\alpha_1}+{\alpha_i}+1)} +
				\dfrac{\lambda_{1i}t}{\Gamma({\alpha_1}+{\alpha_i}+2)}
				\right)t^{\alpha_i}
				\\	+
				\dfrac{h_{10}B_{23}t}{\Gamma({\alpha_1}+2)}	\Big]  t^{\alpha_1},
				\text{ if } \alpha_1\in(0,1]\,\&\,\alpha_2\in(1,2],
				\\
				A_{11}+B_{11}t+
				\Big[\dfrac{l_{1}}{\Gamma({\alpha_1}+1)} +
				\dfrac{l_2t}{\Gamma({\alpha_1}+2)} +  \sum\limits_{i=1}^2 \Big( \dfrac{\varrho_{1i}}{\Gamma({\alpha_1}+{\alpha_i}+1)} \\+
				\dfrac{\lambda_{1i}t}{\Gamma({\alpha_1}+{\alpha_i}+2)}
				\Big)t^{\alpha_i} \Big]  t^{\alpha_1},
				\text{ if }\alpha_s\in(1,2],s=1,2,
			\end{array}
			\right.	
\\								&
						\delta_{12}(t)=	
						\left\{
						\begin{array}{ll}
							A_{12}+ \dfrac{h_{10}A_{24}t^{\alpha_1}}{\Gamma({\alpha_1}+1)} ,\text{ if } \alpha_s\in(0,1],s=1,2,
							\\
							A_{12}+B_{12}t+\dfrac{h_{10}A_{24}t^{\alpha_1}}{\Gamma({\alpha_1}+1)},
							\text{ if } \alpha_1\in(1,2]\ \&\ \alpha_2\in(0,1],
							\\
							A_{12}+\dfrac{h_{10}A_{24}t^{\alpha_1}}{\Gamma({\alpha_1}+1)} +\dfrac{h_{10}B_{24}t^{\alpha_1+1}}{\Gamma({\alpha_1}+2)} ,	
							\text{ if } \alpha_1\in(0,1]\,\&\,\alpha_2\in(1,2],
							\\
							A_{12}+B_{12}t
							+\dfrac{h_{10}A_{24}t^{\alpha_1}}{\Gamma({\alpha_1}+1)} +\dfrac{h_{10}B_{24}t^{\alpha_1+1}}{\Gamma({\alpha_1}+2)},
							\text{ if }\alpha_s\in(1,2],s=1,2,
							\qquad\qquad\qquad\qquad\qquad\qquad
						\end{array}
						\right.
				\\
							& \delta_{13}(t)=	
															\left\{
															\begin{array}{ll}
																A_{13}E_{\alpha_1,1}(\gamma_0 t^{\alpha_1})+ \eta_{10}A_{14}t^{\alpha_1-1}E_{\alpha_1,\alpha_1}(\gamma_0 t^{\alpha_1})\ast E_{\alpha_1,1}(\gamma_0 t^{\alpha_1})  ,
																\\
																\text{ if } \alpha_1\in(0,1]\ \&\ \alpha_2\in(0,2],
																\\
																A_{13}E_{\alpha_1,1}(\gamma_0 t^{\alpha_1})+
																B_{13}tE_{\alpha_1,2}(\gamma_0 t^{\alpha_1})+ \eta_{10}t^{\alpha_1-1}E_{\alpha_1,\alpha_1}(\gamma_0 t^{\alpha_1})
																\\
																\ast\Big(  A_{14}E_{\alpha_1,1}(\gamma_0 t^{\alpha_1})
															+	B_{14}tE_{\alpha_1,2}(\gamma_0 t^{\alpha_1}) \Big) ,
																\text{ if } \alpha_1\in(1,2]\ \&\ \alpha_2\in(0,2],
								\end{array}
									\right.
				\\&		\delta_{14}(t)=	
						\left\{
						\begin{array}{ll}
							A_{14}E_{\alpha_1,1}(\gamma_0 t^{\alpha_1}),\text{ if } \alpha_1\in(0,1]\ \&\ \alpha_2\in(0,2],
							\\
							A_{14}E_{\alpha_1,1}(\gamma_0 t^{\alpha_1})+
							B_{14}tE_{\alpha_1,2}(\gamma_0 t^{\alpha_1}),
							\text{ if } \alpha_1\in(1,2]\ \&\ \alpha_2\in(0,2].
						\end{array}
						\right.
									\\
			&	\delta_{21}(t)=	
							\left\{
							\begin{array}{ll}
								A_{21}+\left( \dfrac{\varrho_{21}}{\Gamma({\alpha_2}+1)}+\dfrac{\varrho_{22}t^{\alpha_2}}{\Gamma({2\alpha_2}+1)}
								\right) t^{\alpha_2} ,\text{ if } \alpha_1\in(0,2]\ \&\ \alpha_2\in(0,1],
								\\
								A_{21}
								+	B_{21}t	+\left( \dfrac{\varrho_{21}}{\Gamma({\alpha_2}+1)}
								+\dfrac{\lambda_{21}t}{\Gamma({\alpha_2}+2)} +\dfrac{\varrho_{22}}{\Gamma({2\alpha_2}+1)}
								+\dfrac{\lambda_{22}t}{\Gamma({2\alpha_2}+2)}
								\right) t^{\alpha_2} ,
							\\	\text{ if } \alpha_1\in(0,2]\ \&\ \alpha_2\in(1,2],
							\end{array}
							\right.
					\\	&	\delta_{22}(t)=	
							\left\{
							\begin{array}{ll}
								A_{22}+ \dfrac{h_{20}A_{24}t^{\alpha_2}}{\Gamma({\alpha_2}+1)} ,\text{ if } \alpha_1\in(0,2]\,\&\,\alpha_2\in(0,1],
								\\
								A_{22}+B_{22}t+ \dfrac{h_{20}A_{24}t^{\alpha_2}}{\Gamma({\alpha_2}+1)}
								+\dfrac{h_{20}B_{24}t^{\alpha_2+1}}{\Gamma({\alpha_2}+2)} ,
								\text{ if } \alpha_1\in(0,2]\,\&\,\alpha_2\in(1,2],
							\end{array}
							\right.
							\\		&	\delta_{23}(t)=	
							\left\{
							\begin{array}{ll}
								A_{23}+ \dfrac{g_{20}A_{24}t^{\alpha_2}}{\Gamma({\alpha_2}+1)} ,\text{ if } \alpha_1\in(0,2]\,\&\,\alpha_2\in(0,1],
								\\
								A_{23}+B_{23}t+ \dfrac{g_{20}A_{24}t^{\alpha_2}}{\Gamma({\alpha_2}+1)}
								+\dfrac{g_{20}B_{24}t^{\alpha_2+1}}{\Gamma({\alpha_2}+2)} ,
								\text{ if } \alpha_1\in(0,2]\,\&\,\alpha_2\in(1,2],
							\end{array}
							\right.
							\\	&		\text{ and }
							\delta_{24}(t)=	
							\left\{
							\begin{array}{ll}
								A_{24} ,\text{ if } \alpha_1\in(0,2]\,\&\,\alpha_2\in(0,1],
								\\
								A_{24}+B_{24}t ,\text{ if } \alpha_1\in(0,2]\,\&\,\alpha_2\in(1,2],
							\end{array}
							\right.
					\end{aligned}	\end{eqnarray*}	
				Here the two-parameter Mittag-Leffler function \cite{Podlubny1999}  is defined  as  $E_{a,b}(z)=\sum\limits_{k=0}^\infty\dfrac{z^k}{\Gamma(ak+b)},$  $\mathfrak{R}(a),\mathfrak{R}(b)>0,$
				$\ast$ is the convolution of  two functions defined as	$	 w_1(t)\ast w_2(t)=\int_{0}^tw_1(y)w_{2}(t-y)dy, $
				$$ \begin{array}{llll}
					l_1=\eta_{10}A_{12}+h_{10}A_{23},
					&\, \varrho_{11}=h_{10}\eta_{10}A_{24},
					&\, \lambda_{11}=h_{10}\eta_{10}B_{24},
					&\, \gamma_0=p_{10}a_{21}^2-d_{10}a_{21},
					\\
					l_2=\eta_{10}B_{12}+h_{10}B_{23},
					&\,  \varrho_{12}=h_{10}g_{20}A_{24},
					& \, \lambda_{12}=h_{10}g_{20}B_{24},
					&\,	\lambda_{21}=h_{20}B_{23}+g_{20}B_{22},
					
					\\
					\varrho_{21}=h_{20}A_{23}+g_{20}A_{22} ,&	\, \varrho_{22}=2h_{20}g_{20}A_{24},
					&\,
					\lambda_{22}=2h_{20}g_{20}B_{24},&	 \, \&\, A_{si},B_{si}\in\mathbb{R},
				\end{array}
				$$
		for	$s=1,2,i=1,2,3,4.$	Finally, we observe that the above {analytical} solution \eqref{solutionu1ofeg2} 
				for the given coupled equations \eqref{poly-exp} satisfies the initial conditions \eqref{icpoly-exp} with
				$$\begin{aligned}
					&\omega_1(x_1,x_2)=A_{11}+A_{12}x_1+(A_{13}+A_{14}x_1)e^{-a_{21}x_2}, \\	& \omega_2(x_1,x_2)=A_{21}+A_{22}x_1+(A_{23}+A_{24}x_1)x_2,
					\\ &
					\vartheta_1(x_1,x_2)=B_{11}+B_{12}x_1+(B_{13}+B_{14}x_1)e^{-a_{21}x_2},
				\\	\&\,&
					\vartheta_{2}(x_1,x_2)=B_{21}+B_{22}x_1+(B_{23}+B_{24}x_1){x_2}.
				\end{aligned}
				$$
Furthermore, we note that the {analytical} solution \eqref{solutionu1ofeg2} of the given coupled equations  \eqref{poly-exp} is valid for all $ \alpha_s, \alpha_s\in(0,2], s=1,2.
	  $ {Also, we observe that when  $\alpha_s=1$ and $\alpha_s=2$, the obtained solutions  \eqref{solutionu1ofeg2}  coincide with integer-order cases of the given coupled system \eqref{poly-exp}.}	
			\end{example}
			%
			
			%
			%
			%
			%
			%
 	\begin{example}\label{eg4}
				Consider the   {multi-component} $(2+1)$-dimensional coupled  nonlinear {TFDCWEs}  for {${\alpha_s}\in(0,2],s=1,2,$} as in the form
				{	\begin{eqnarray}\label{tri-poly-expeq}
						\begin{aligned}
							\dfrac{\partial ^{\alpha_1} u_1}{\partial t^{\alpha_1}}
							=
							&	k_{{10}}{\frac {\partial^2 u_{{1}} }{\partial x_1^2}}
							+
							{\frac {g_{10} }{b_{11}} }{\frac {\partial^2 u_{{2}} }{\partial x_1^2}}
							+
							{\frac {\partial }{	\partial x_2}}
							\left[
							\left(
							-c_{{11}}u_{{2}}
							+p_{{10}} \right)
							{\frac {\partial u_{{1}}}{\partial x_2}}
							+
							\left(
							c_{{11}}u_{{1}}
							+c_{{10}}
							\right)
							{\frac {\partial u_{{2}}}{
									\partial x_2}}
							\right]
						\\&	+
							\eta_{{10}}
							{\frac {\partial u_{{1}} }{\partial x_1}}
							+
							g_{{10}}
							{\frac {\partial u_{{2}} }{\partial x_1}}
							+
							d_{{10}}
							{\frac {\partial u_{{1}}  }{\partial x_2}}
							,
							\\
							\dfrac{\partial ^{\alpha_2} u_2}{\partial t^{\alpha_2}}
							=
							&
							q_{{20}}
							{\frac {\partial^2 u_{{2}} }{\partial x_1^2}}
							+
							{\frac {\partial }{	\partial x_2}}
							\left[
							\left(
							-c_{{21}}u_{{2}}
							+p_{{20}} \right)
							{\frac {\partial u_{{1}}}{\partial x_2}}
							+
							\left(
							c_{{21}}u_{{1}}
							+c_{{20}}
							\right)
							{\frac {\partial u_{{2}}}{\partial x_2}}
							\right]
					\\&		+
							g_{{20}}
							{\frac {\partial u_{{2}} }{\partial x_1}}
							+
							h_{{20}}
							{\frac {\partial u_{{2}} }{\partial x_2}},
						\end{aligned}
				\end{eqnarray}	}
				with appropriate initial conditions
				{\begin{eqnarray}
						\begin{aligned}\label{ictri-poly-expeq}
							&(i) \, u_s(x_1,x_2,0)= \omega_s(x_1,x_2)
							\,\text{ if }\  {\alpha_s}\in(0,1],s=1,2,\\
							&(ii)\,	u_s(x_1,x_2,0)= \omega_s(x_1,x_2),\, \,
							{\dfrac{\partial u_s}{\partial t}}\big{|}_{t=0}=\vartheta_s(x_1,x_2)
							\,	\text{if}\  {\alpha_s}\in(1,2],s=1,2.
						\end{aligned}
				\end{eqnarray}}
				The product linear space
				$$\begin{aligned}
					{{	\mathbf{\hat{W}}^{1}_8}}= & \text{Span}\left\{\sin(\sqrt{a_{10}}x_1),x_2\sin(\sqrt{a_{10}}x_1),\cos(\sqrt{a_{10}}x_1),x_2\cos(\sqrt{a_{10}}x_1)\right\}
				\\&\times\text{Span}\left\{1,x_2,e^{b_{11}x_1},x_2e^{b_{11}x_1}\right\},a_{10}>0,
				\end{aligned}
				$$ is  invariant for the coupled equations \eqref{tri-poly-expeq}, which is discussed in  case 11 of Table 2.
				{Now, we apply the above  procedure here and obtain} the generalized separable  {analytical} solution for the given coupled equations  \eqref{tri-poly-expeq} {along with $\eta_{10}=d_{10}=h_{20}=0$}
				as follows:
				\begin{eqnarray}
					\begin{aligned}
						\label{exactsolution4} u_1(x_1,x_2,t)=&(\delta_{11}(t)+\delta_{12}(t)x_2)\sin(\sqrt{a_{10}}x_1)+(\delta_{13}(t)+\delta_{14}(t)x_2)\cos(\sqrt{a_{10}}x_1),
						\\
						u_2(x_1,x_2,t)=&\delta_{21}(t)+\delta_{22}(t)x_2+\delta_{23}(t)e^{b_{11}x_1}+\delta_{24}(t)x_2e^{b_{11}x_1},
				\end{aligned}\end{eqnarray}
				where
				{	 $$\begin{aligned}
					&	\delta_{1i}(t)=\left\{\begin{array}{ll}
							A_{1i}E_{{\alpha_1},1}(-k_{10}a_{10}t^{\alpha_1})
							,\text{ if } \, {\alpha_1}\in(0,1]\,\&\,{\alpha_2}\in(0,2],
							\\
							A_{1i}E_{{\alpha_1},1}(-k_{10}a_{10}t^{\alpha_1})+	B_{1i}tE_{{\alpha_1},2}(-k_{10}a_{10}t^{\alpha_1}), 
							\\
							\text{ if } \, {\alpha_1}\in(1,2]\,\&\,{\alpha_2}\in(0,2],i=1,2,3,4,
						\end{array}\right.
				\\&
						\delta_{2i}(t)=\left\{\begin{array}{ll}
							A_{2i}	,\text{ if }  {\alpha_1}\in(0,2]\,\&\, {\alpha_2}\in(0,1],
							\\
							A_{2i}+B_{2i}t	,\text{ if } \, {\alpha_1}\in(0,2]\,\&\, {\alpha_2}\in(1,2],i=1,2,
							\text{ and }	
						\end{array}\right.
						\\
						&	\delta_{2i}(t)=\left\{\begin{array}{ll}
							A_{2i}E_{ {\alpha_2},1}(\gamma_0t^ {\alpha_2})	,\text{ if } \, {\alpha_1}\in(0,2]\,\&\, {\alpha_2}\in(0,1],
							\\
							A_{2i}E_{ {\alpha_2},1}(\gamma_0t^ {\alpha_2})+B_{2i}tE_{ {\alpha_2},2}(\gamma_0t^ {\alpha_2}),\text{ if } \, {\alpha_1}\in(0,2]\,\&\, {\alpha_2}\in(1,2],i=3,4.\qquad\qquad\qquad\qquad\qquad
						\end{array}\right.
					\end{aligned}
					$$}
				Here 	$
				\gamma_0=q_{20}b_{11}^2-g_{20}b_{11},
				$ and $A_{si},B_{si}\in\mathbb{R},s=1,2,i=1,2,3,4.$
				{Additionally}, we note that the obtained generalized separable {analytical} solution \eqref{exactsolution4} of the given coupled equations \eqref{tri-poly-expeq} satisfies the initial conditions \eqref{ictri-poly-expeq} with
				$$ \begin{aligned}
					&	\omega_1(x_1,x_2)=(A_{11}+A_{12}x_2) \sin(\sqrt{a_{10}}x_1)+(A_{13}+A_{14}x_2)\cos(\sqrt{a_{10}}x_1),
					\\	&\vartheta_1(x_1,x_2)=(B_{11}+B_{12}x_2)\sin(\sqrt{a_{10}}x_1)+(B_{13}+B_{14}x_2)\cos(\sqrt{a_{10}}x_1),
					\\& \omega_2(x_1,x_2)=A_{21}+A_{22}x_2+ (A_{23}+A_{24}x_2)e^{b_{11}x_1},	\\	\& \,&	\vartheta_2(x_1,x_2)=B_{21}+B_{22}x_2+(B_{23}+B_{24}x_2)e^{b_{11}x_1}.
				\end{aligned}$$
						{Also, we would like to mention here that the obtained {analytical} solution \eqref{exactsolution4} of the {considered} coupled equations \eqref{tri-poly-expeq} is valid for all $\alpha_s, \alpha_s\in(0,2], s=1,2.$
				}{Additionally, we note that when  $\alpha_s=1$ and $\alpha_s=2$, the obtained solutions \eqref{exactsolution4} coincide with integer-order cases of the coupled system \eqref{tri-poly-expeq}. }
			\end{example}
		
\subsubsection{Generalized separable {analytical} solutions for the {initial-boundary} value problems of the given  coupled equations  \eqref{cdsystemeqn}}
			Next, we provide {the  systematic approach of the developed method to derive } {analytical} solutions for the initial and the Dirichlet boundary value problems of the considered coupled equations \eqref{cdsystemeqn}  {based on the determined invariant product subpaces}.
			\begin{example}\label{linear}
				Now, we consider the following    {multi-component} $(2+1)$-dimensional coupled  linear {TFDCWEs}  as
				{	\begin{eqnarray}\label{linearcdsystem}
						\begin{aligned}
							\dfrac{\partial ^{\alpha_1} u_1}{\partial t^{\alpha_1}}
							=
							&		
							k_{{10}}
							{\frac{\partial^2 u_{{1}} }{\partial x_1^2}}
							+	q_{{10}}
							{\frac {\partial^2 u_{{2}} }{\partial x_1^2}}
							+p_{{10}}
							{\frac {\partial^2 u_{{1}}}{\partial x_2^2}} 					
							+q_{{10}}b_{11}
							{\frac {\partial u_{{2}} }{\partial x_1}}
							,
							\\
							\dfrac{\partial ^{\alpha_2} u_2}{\partial t^{\alpha_2}}
							=
							&q_{{20}}
							{\frac{\partial^2 u_{{2}} }{\partial x_1^2}}
							+	c_{{20}}
							{\frac {\partial^2 u_{{2}} }{\partial x_2^2}}
							+g_{{20}}
							{\frac {\partial u_{{2}} }{\partial x_1}},{\alpha_s}\in(0,2],s=1,2,
						\end{aligned}
				\end{eqnarray}	}
				with appropriate initial conditions
				{	\begin{eqnarray}
						\begin{aligned}\label{iclinear}
							&(i) \, u_s(x_1,x_2,0)= \omega_s(x_1,x_2)
							\,\text{ if }\  {\alpha_s}\in(0,1],s=1,2,\\
							&(ii)\,	u_s(x_1,x_2,0)= \omega_s(x_1,x_2),\, \,
							{\dfrac{\partial u_s}{\partial t}}\big{|}_{t=0}=\vartheta_s(x_1,x_2)
							\,	\text{ if }\  {\alpha_s}\in(1,2],s=1,2,
						\end{aligned}
				\end{eqnarray}}
				and the Dirichlet boundary conditions
				\begin{eqnarray}
			\begin{aligned}\label{bclinear}
						(\text{B-I}) &
						\left\{\begin{array}{ll}
							u_1(0,x_2,t)=\xi_{11}(t)\sin(\sqrt{b_{20}}x_2)+\xi_{12}(t)\cos(\sqrt{b_{20}}x_2),\text{if}\ (0,x_2,t)\in \partial \Omega\times[0,\infty),
							\\
							u_1(\sigma_1,x_2,t)=\xi_{13}(t)\sin(\sqrt{b_{20}}x_2)+\xi_{14}(t)\cos(\sqrt{b_{20}}x_2),
							\\
							\text{if}\ (\sigma_1,x_2,t)\in \partial \Omega\times[0,\infty),
							\\
							u_1(x_1,0,t)=\xi_{21}(t)\sin(\sqrt{a_{10}}x_1)+\xi_{22}(t)\cos(\sqrt{a_{10}}x_1),\text{if}\ (x_1,0,t)\in \partial \Omega\times[0,\infty),
							\\
							u_1(x_1,\sigma_2,t)= \xi_{23}(t)\sin(\sqrt{a_{10}}x_1) +\xi_{24}(t)\cos(\sqrt{a_{10}}x_1),
							\\ \text{if}\ (x_1,\sigma_2,t)\in \partial \Omega\times[0,\infty),
						\end{array}	
						\right.		
						\\
						(\text{B-II})&	\left\{\begin{array}{ll}
							
							u_2(0,x_2,t)=f_{11}(t)\sin(\sqrt{b_{20}}x_2)+
							f_{12}(t)\cos(\sqrt{b_{20}}x_2),\text{if}\ (0,x_2,t)\in \partial \Omega\times[0,\infty),
							\\
							u_2(\sigma_1,x_2,t)=f_{13}(t)\sin(\sqrt{b_{20}}x_2)+
							f_{14}(t)\cos(\sqrt{b_{20}}x_2),
							\\
							\text{if}\ (\sigma_1,x_2,t)\in \partial \Omega\times[0,\infty),
							\\
							u_2(x_1,0,t)=f_{21}(t)+
							f_{22}(t)e^{-b_{11}x_1},\text{if}\ (x_1,0,t)\in \partial \Omega\times[0,\infty),
							\\
							u_2(x_1,\sigma_2,t)=f_{23}(t)+
							f_{24}(t)e^{-b_{11}x_1},\text{if}\ (x_1,\sigma_2,t)\in \partial \Omega\times[0,\infty).
						\end{array}	
						\right.		
					\end{aligned}
				\end{eqnarray}
				The given coupled equations \eqref{linearcdsystem} are defined on a domain $\Omega\times[0,\infty),$ where $\Omega=\{(x_1,x_2)\in\mathbb{R}^2|\, 0\leq x_i\leq \sigma_i,\sigma_i\in\mathbb{R},i=1,2\}$ is a subset of $\mathbb{R}^2$ whose boundary $\partial\Omega$ is defined by the lines $x_1=0,x_2=0,x_1=\sigma_1,$ and $x_2=\sigma_2.$
				It is easy to observe that the coupled equations \eqref{linearcdsystem} admit the following invariant  product linear space
				$$\begin{aligned}
					{{	\mathbf{\hat{W}}^{1}_8}}= & \text{Span}\left\{
					\sin(\sqrt{a_{10}}x_1)\sin(\sqrt{b_{20}}x_2), \cos(\sqrt{a_{10}}x_1)\sin(\sqrt{b_{20}}x_2), \sin(\sqrt{a_{10}}x_1)\cos(\sqrt{b_{20}}x_2),  \right. \\&\left. \cos(\sqrt{a_{10}}x_1)\cos(\sqrt{b_{20}}x_2)\right\}  \\&  		     \times
					\text{Span}\left\{\sin(\sqrt{b_{20}}x_2), \cos(\sqrt{b_{20}}x_2),
					e^{-b_{11}x_1}\sin(\sqrt{b_{20}}x_2),  e^{-b_{11}x_1}\cos(\sqrt{b_{20}}x_2)\right\},
				\end{aligned}
				$$
				where $a_{10},b_{20}\geq0.$
				Applying the similar procedure as explained in  Example \ref{eg1}, the obtained generalized separable  {analytical} solution for the given coupled equations \eqref{linearcdsystem}  is $\mathbf{U}=(u_1,u_2)$ with $u_s=u_s(x_1,x_2,t),s=1,2,$
				as
				\begin{eqnarray}
					\begin{aligned}\label{exactsolutionlinearcdsystem}
						u_1=& [\delta_{11}(t)\sin(\sqrt{a_{10}}x_1) +\delta_{13}(t)\cos(\sqrt{a_{10}}x_1)] \sin(\sqrt{b_{20}}x_2)
					\\&	+ [\delta_{12}(t)\sin(\sqrt{a_{10}}x_1) +\delta_{14}(t)\cos(\sqrt{a_{10}}x_1)] \cos(\sqrt{b_{20}}x_2),
						\\
						u_2=& [\delta_{21}(t)+\delta_{23}(t)e^{-b_{11}x_1}] \sin(\sqrt{b_{20}}x_2) +[\delta_{22}(t)+\delta_{24}(t)e^{-b_{11}x_1}] \cos(\sqrt{b_{20}}x_2),
				\end{aligned}\end{eqnarray}
				where
				{	\begin{eqnarray}
						\begin{aligned}\label{linearCDsystemsolutionofFODEs}
							\delta_{1i}(t)=&\left\{\begin{array}{ll}
								A_{1i}E_{{\alpha_1},1}(\gamma_0t^{\alpha_1}),
								\text{ if } \, {\alpha_1}\in(0,1]\,\&\,{\alpha_2}\in(0,2],
								\\
								A_{1i}E_{{\alpha_1},1}(\gamma_0t^{\alpha_1})+	B_{1i}tE_{{\alpha_1},2}(\gamma_0t^{\alpha_1}), \text{ if } \,{\alpha_1}\in(1,2] \\
								\&\,{\alpha_2}\in(0,2],i=1,2,3,4,
							\end{array}\right.
							\\
							\delta_{2i}(t)=&\left\{\begin{array}{ll}
								A_{2i}E_{{\alpha_2},1}(-c_{20}b_{20}t^{\alpha_2}),
								\text{ if} \,{\alpha_1}\in(0,2]\,\&\,{\alpha_2}\in(0,1],
								\\
								A_{2i}E_{{\alpha_2},1}(-c_{10}b_{20}t^{\alpha_2})+	B_{2i}tE_{{\alpha_2},2}(-c_{10}b_{20}t^{\alpha_2}), 
								\\
								\text{ if } \, {\alpha_1}\in(0,2]\,\&\,{\alpha_2}\in(1,2],i=1,2,
							\end{array}\right.
							\\
							\&\,			\delta_{2i}(t)=&\left\{\begin{array}{ll}
								A_{2i}E_{{\alpha_2},1}(\gamma_1t^{\alpha_2}),
								\text{ if } \, {\alpha_1}\in(0,2]\,\&\,{\alpha_2}\in(0,1],
								\\
								A_{2i}E_{{\alpha_2},1}(\gamma_1t^{\alpha_2})+	B_{2i}tE_{{\alpha_2},2}(\gamma_1t^{\alpha_2}), 
								\text{ if } \, {\alpha_1}\in(0,2]\,\&\,{\alpha_2}\in(1,2],i=3,4.
							\end{array}\right.
						\end{aligned}
				\end{eqnarray}}
				Here  $\gamma_0=-k_{10}a_{10}-p_{10}b_{20}, \gamma_1=q_{20}b_{11}^2-c_{20}b_{20}-g_{20}b_{11},$ and $A_{si},B_{si}\in\mathbb{R},s=1,2,i=1,2,3,4.$
				We also observe that the above {analytical} solution \eqref{exactsolutionlinearcdsystem} for the given coupled equations \eqref{linearcdsystem} satisfies the {initial-boundary conditions} \eqref{iclinear}-\eqref{bclinear} with
				$$\begin{aligned}
					\omega_1(x_1,x_2)=&(A_{11}\sin(\sqrt{a_{10}}x_1) +A_{13}\cos(\sqrt{a_{10}}x_1)) \sin(\sqrt{b_{20}}x_2)
				\\&	+ (A_{12}\sin(\sqrt{a_{10}}x_1) +A_{14}\cos(\sqrt{a_{10}}x_1)) \cos(\sqrt{b_{20}}x_2),
					\\
					\vartheta_1(x_1,x_2)=&(B_{11}\sin(\sqrt{a_{10}}x_1) +B_{13}\cos(\sqrt{a_{10}}x_1)) \sin(\sqrt{b_{20}}x_2)
				\\&	+ (B_{12}\sin(\sqrt{a_{10}}x_1) +B_{14}\cos(\sqrt{a_{10}}x_1)) \cos(\sqrt{b_{20}}x_2)
					,
					\\
					\omega_2(x_1,x_2)=&(A_{21}+A_{23}e^{-b_{11}x_1})\sin(\sqrt{b_{20}}x_2)+(A_{22}+A_{24}e^{-b_{11}x_1})\cos(\sqrt{b_{20}}x_2),
					\\
					\vartheta_2(x_1,x_2)=&(B_{21}+B_{23}e^{-b_{11}x_1})\sin(\sqrt{b_{20}}x_2)+(B_{22}+B_{24}e^{-b_{11}x_1})\cos(\sqrt{b_{20}}x_2)
					,
				\end{aligned}
				$$
				$$
				\begin{array}{lllll}
					&
					\xi_{11}(t)=\delta_{12}(t),
										&	
					\xi_{13}(t)=C_1\delta_{11}(t)+ C_2\delta_{12}(t),
				\\	& 	\xi_{12}(t)=\delta_{14}(t),
					&	\xi_{14}(t)=C_1\delta_{13}(t)+ C_2\delta_{14}(t),
					\\	&	\xi_{21}(t)=\delta_{13}(t),
									&	
					\xi_{23}(t)=C_3\delta_{11}(t)+ C_4\delta_{13}(t),
				\\	& \xi_{22}(t)=\delta_{14}(t),	&	\xi_{24}(t)=C_3\delta_{12}(t)+ C_4\delta_{14}(t),
					\\	&	f_{11}(t)=\delta_{21}(t)+\delta_{23}(t),
					&f_{12}(t)=\delta_{22}(t)+\delta_{24}(t),
				\\	&		
					f_{13}(t)=\delta_{21}(t)+\delta_{23}(t)e^{-b_{11}\sigma_1},
					&	f_{14}(t)=\delta_{22}(t)+\delta_{24}(t)e^{-b_{11}\sigma_1},
					\\
					&
					f_{21}(t)=\delta_{22}(t),
									&f_{23}(t)=C_3\delta_{21}(t)+C_4\delta_{22}(t),
				\\	&	f_{22}(t)=\delta_{24}(t),	\&\,& f_{24}(t)=C_3\delta_{23}(t)+C_4\delta_{24}(t),
				\end{array}$$
				where  $C_1=\sin(\sigma_1\sqrt{a_{10}}), C_2=\cos(\sigma_1\sqrt{a_{10}}), C_3=\sin(\sigma_2\sqrt{b_{20}}), C_4=\cos(\sigma_2\sqrt{b_{20}}),$ and the functions $\delta_{si}$ as given in \eqref{linearCDsystemsolutionofFODEs} for $s=1,2,$ and $i=1,2,3,4.$
				{Additionally, we observe that the obtained {analytical} solution \eqref{exactsolutionlinearcdsystem} of the given coupled equations \eqref{linearcdsystem}  is valid for all $\alpha_s, \alpha_s\in(0,2], s=1,2.$ }
				{Also, we note that  the derived solutions \eqref{exactsolutionlinearcdsystem} coincide with integer-order cases of the discussed coupled system \eqref{linearcdsystem} if  $\alpha_s=1$ and $\alpha_s=2$. }
				\begin{note}
	It is interesting to observe that in the above Example \ref{linear}, the obtained {analytical} solution \eqref{exactsolutionlinearcdsystem} for the given initial-boundary value problem \eqref{linearcdsystem}-\eqref{bclinear}   is spatially periodic. From the physical point of view, linear analysis of possible patterns in an experiment is observed to be repeating if the experimental domain is subjected to periodic boundary conditions. While studying the stability of the underlying patterns, these periodic boundary conditions in a bounded domain are crucial in both fractional-order and integer-order PDEs. 
				\end{note}
				
			\end{example}
			%
			%
			%
			%
			%
			%
			%
		\begin{example}\label{eg3}
				Next, we consider the following  {multi-component} $(2+1)$-dimensional coupled nonlinear {TFDCWEs}  as
				{	\begin{eqnarray}\label{tri-polyeqa}
						\begin{aligned}
							\dfrac{\partial ^{\alpha_1} u_1}{\partial t^{\alpha_1}}
							=
							&
							{\frac {\partial }{\partial x_1}}
							\Big[\left(k_{{10}}-q_{{11}}u_{{2}}\right)
							{\frac{\partial u_{{1}} }{\partial x_1}}
							+	\left(q_{{11}}u_{{1}}+q_{{10}}\right)
							{\frac {\partial u_{{2}} }{\partial x_1}} \Big]
						\\&	+p_{{10}}
							{\frac {\partial^2 u_{{1}}}{\partial x_2^2}} 								+\eta_{{10}}
							{\frac {\partial u_{{1}} }{\partial x_1}}
						 +
							d_{{10}}{\frac {\partial u_{{1}}  }{\partial x_2}}
							,
							\\
							\dfrac{\partial ^{\alpha_2} u_2}{\partial t^{\alpha_2}}
							=
							&
							{\frac {\partial }{\partial x_1}}
							\left[	
							\left(
							k_{{20}}	-q_{{21}}u_{{2}}
							\right)
							{\frac {\partial u_{{1}} }{	\partial x_1}}
							+
							\left(
							q_{{21}}u_{{1}}
							+q_{{20}}
							\right)
							{\frac {\partial u_{{2}} }{\partial x_1}}
							\right]
						\\&	+
							c_{{20}} 							{\frac {\partial^2 u_{{2}}}{\partial x_2^2}}
							+
							g_{{20}}
							{\frac {\partial u_{{2}} }{\partial x_1}}
						 	+h_{{20}}
							{\frac {\partial u_{{2}} }{\partial x_2}},
						\end{aligned}
				\end{eqnarray}	}
				${\alpha_s}\in(0,2],s=1,2,$	with appropriate initial conditions
				{	\begin{eqnarray}
						\begin{aligned}\label{icpoly-tri}
							&(i) \, u_s(x_1,x_2,0)= \omega_s(x_1,x_2)
							\,\text{ if }\  {\alpha_s}\in(0,1],s=1,2,\\
							&(ii)\,	u_s(x_1,x_2,0)= \omega_s(x_1,x_2),\, \,
							{\dfrac{\partial u_s}{\partial t}}\big{|}_{t=0}=\vartheta_s(x_1,x_2)
							\,	\text{if}\  {\alpha_s}\in(1,2],s=1,2,
						\end{aligned}
				\end{eqnarray}}
				and the Dirichlet boundary conditions
				\begin{eqnarray}
					\label{bc-tri-polyeqa}
					(\text{B-I}) &
					\left\{\begin{array}{ll}
						u_1(0,x_2,t)=\xi_{11}(t)\sin(\sqrt{a_{20}}x_2)+\xi_{12}(t)\cos(\sqrt{a_{20}}x_2),\text{if}\ (0,x_2,t)\in \partial \Omega\times[0,\infty),
						\\
						u_1(\sigma_1,x_2,t)=\xi_{13}(t)\sin(\sqrt{a_{20}}x_2)+\xi_{14}(t)\cos(\sqrt{a_{20}}x_2),
						\\
						\text{if}\ (\sigma_1,x_2,t)\in \partial \Omega\times[0,\infty),
						\\
						u_1(x_1,0,t)=\xi_{21}(t)+\xi_{22}(t)x_1,\text{if}\ (x_1,0,t)\in \partial \Omega\times[0,\infty),
						\\
						u_1(x_1,\sigma_2,t)= \xi_{23}(t) +\xi_{24}(t)x_1, \text{if}\ (x_1,\sigma_2,t)\in \partial \Omega\times[0,\infty),
					\end{array}	
					\right.		
			\\
								\label{bc2-tri-polyeqa}
					(\text{B-II})&	\left\{\begin{array}{ll}
						u_2(0,x_2,t)=f_{11}(t)\sin(\sqrt{b_{20}}x_2)+f_{12}(t)\cos(\sqrt{b_{20}}x_2),\text{if}\ (0,x_2,t)\in \partial \Omega\times[0,\infty),
					\\
						u_2(\sigma_1,x_2,t)=f_{13}(t)\sin(\sqrt{b_{20}}x_2)+f_{14}(t)\cos(\sqrt{b_{20}}x_2),\\
						\text{if}\ (\sigma_1,x_2,t)\in \partial \Omega\times[0,\infty),
						\\
						u_2(x_1,0,t)=f_{21}(t)+f_{22}(t)x_1,\text{if}\ (x_1,0,t)\in \partial \Omega\times[0,\infty),
						\\
						u_2(x_1,\sigma_2,t)= f_{23}(t) +f_{24}(t)x_1, \text{if}\ (x_1,\sigma_2,t)\in \partial \Omega\times[0,\infty).
					\end{array}	
					\right.		
				\end{eqnarray}
				The given coupled equations  \eqref{tri-polyeqa} are defined  on a domain $\Omega\times[0,\infty),$ where $\Omega=\{(x_1,x_2)\in\mathbb{R}^2|\, 0\leq x_i\leq \sigma_i,\sigma_i\in\mathbb{R},i=1,2\}$ is a subset of $\mathbb{R}^2$ whose boundary $\partial\Omega$ is defined by the lines $x_1=0,x_2=0,x_1=\sigma_1,$ and $x_2=\sigma_2.$
				It is given in case 7 of Table 1 that the coupled equations \eqref{tri-polyeqa} admit the following invariant  product linear space
				$$\begin{aligned}
					{{	\mathbf{\hat{W}}^{1}_8}}= & \text{Span}\{\sin(\sqrt{a_{20}}x_2),\cos(\sqrt{a_{20}}x_2),x_1\sin(\sqrt{a_{20}}x_2),x_1\cos(\sqrt{a_{20}}x_2)\}        \\
					& \times\text{Span}\{\sin(\sqrt{b_{20}}x_2),\cos(\sqrt{b_{20}}x_2),x_1\sin(\sqrt{b_{20}}x_2),x_1\cos(\sqrt{b_{20}}x_2)\},
				\end{aligned}
				$$
				where $a_{20},b_{20}\geq0.$
				Applying a similar procedure as explained in  Example \ref{eg1}, the obtained generalized separable  {analytical} solution for the given coupled equations \eqref{tri-polyeqa} {along with $\eta_{10}=d_{10}=h_{20}=g_{20}=0$} is $\mathbf{U}=(u_1,u_2)$ with $u_s=u_s(x_1,x_2,t),s=1,2,$
				as \begin{eqnarray}
					\begin{aligned}\label{exactsolution3}
						u_1(x_1,x_2,t)=& (\delta_{11}(t)+\delta_{13}(t)x_1)\sin(\sqrt{a_{20}}x_2)+(\delta_{12}(t)+\delta_{14}(t)x_1)\cos(\sqrt{a_{20}}x_2),
						\\
						u_2(x_1,x_2,t)=&(\delta_{21}(t)+\delta_{23}(t)x_1)\sin(\sqrt{b_{20}}x_2)+(\delta_{22}(t)+\delta_{24}(t)x_1)\cos(\sqrt{b_{20}}x_2),
				\end{aligned}\end{eqnarray}
				where
				{	\begin{eqnarray}
						\begin{aligned}\label{eg3solnFODEs}
							\delta_{1i}(t)=&\left\{\begin{array}{ll}
								A_{1i}E_{{\alpha_1},1}(-p_{10}a_{20}t^{\alpha_1}),
								\text{ if } \,{\alpha_1}\in(0,1]\,\&\,{\alpha_2}\in(0,2],
								\\
								A_{1i}E_{{\alpha_1},1}(-p_{10}a_{20}t^{\alpha_1})+	B_{1i}tE_{{\alpha_1},2}(-p_{10}a_{20}t^{\alpha_1}), \\
								\text{if } \, {\alpha_1}\in(1,2]\,\&\,{\alpha_2}\in(0,2],i=1,2,3,4,
							\end{array}\right.
							\\
							\text{and }		\delta_{2i}(t)=&\left\{\begin{array}{ll}
								A_{2i}E_{{\alpha_2},1}(-c_{20}b_{20}t^{\alpha_2}),
								\text{ if } \, {\alpha_1}\in(0,2]\,\&\,{\alpha_2}\in(0,1],
								\\
								A_{2i}E_{{\alpha_2},1}(-c_{10}b_{20}t^{\alpha_2})+	B_{2i}tE_{{\alpha_2},2}(-c_{10}b_{20}t^{\alpha_2}),
								\\ \text{ if } \, {\alpha_1}\in(0,2]\,\&\,{\alpha_2}\in(1,2],i=1,2,3,4.
							\end{array}\right.
						\end{aligned}
				\end{eqnarray}}
				Here  $A_{si},B_{si}\in\mathbb{R},s=1,2,i=1,2,3,4.$
				We also observe that the above {analytical} solution \eqref{exactsolution3} of the given coupled equations \eqref{tri-polyeqa} satisfies the {initial-boundary conditions} \eqref{icpoly-tri}-\eqref{bc2-tri-polyeqa} with
				$$\begin{aligned}
				&\begin{array}{lllll}	&\omega_1(x_1,x_2)=(A_{11}+A_{13}x_1)\sin(\sqrt{a_{20}}x_2)+(A_{12}+A_{14}x_1)\cos(\sqrt{a_{20}}x_2)
					,
					\\ &
					\vartheta_1(x_1,x_2)=(B_{11}+B_{13}x_1)\sin(\sqrt{a_{20}}x_2)+(B_{12}+B_{14}x_1)\cos(\sqrt{a_{20}}x_2)
					, \\&
					\omega_2(x_1,x_2)=(A_{21}+A_{23}x_1)\sin(\sqrt{b_{20}}x_2)+(A_{22}+A_{24}x_1)\cos(\sqrt{b_{20}}x_2),
					\\
					&
					\vartheta_2(x_1,x_2)=(B_{21}+B_{23}x_1)\sin(\sqrt{b_{20}}x_2)+(B_{22}+B_{24}x_1)\cos(\sqrt{b_{20}}x_2)  ,
					\end{array}	\\ &\begin{array}{lllll}
					&	\xi_{11}(t)=\delta_{11}(t),
							\quad	&	\xi_{12}(t)=	\xi_{21}(t)=\delta_{12}(t), \quad	&	\xi_{22}(t)=\delta_{14}(t),	
				\\&	f_{11}(t)=\delta_{21}(t), 
				&	f_{12}(t)=	f_{21}(t)=\delta_{22}(t), 
				& 
				f_{22}(t)=\delta_{24}(t),
			\end{array}
								\\	&\begin{array}{lllll}
				& 
			\xi_{23}(t)=C_1\delta_{11}(t)+C_2 \delta_{12}(t),		
			&\xi_{13}(t)=\delta_{11}(t)+ \sigma_{1}\delta_{13}(t),	&\xi_{14}(t)=\delta_{12}(t)+ \sigma_{1}\delta_{14}(t), 
				\\	&	\xi_{24}(t)=C_1\delta_{13}(t)+ C_2\delta_{14}(t),
					 & 	f_{13}(t)=\delta_{21}(t)+{\sigma_1}\delta_{23}(t), & 	f_{14}(t)=\delta_{22}(t)+{\sigma_1}\delta_{24}(t),
					\\
					& 
					f_{23}(t)=C_3\delta_{21}(t)+C_4\delta_{22}(t),\, \& \,	 &  	f_{24}(t)=C_3\delta_{23}(t)+C_4\delta_{24}(t),
				\end{array}					 				\end{aligned} $$
				where $C_1=\sin(\sigma_2\sqrt{a_{20}}), C_2=\cos(\sigma_2\sqrt{a_{20}}),
				C_3=\sin(\sigma_2\sqrt{b_{20}}), C_4=\cos(\sigma_2\sqrt{b_{20}}),$ and the functions $\delta_{si}$ as given in \eqref{eg3solnFODEs} for $s=1,2,$ and $i=1,2,3,4.$
				
				{Additionally, we wish to mention here that  when $I_s=\tau_1=0,s=1,2,$  the  given system  \eqref{2+1diffyang} can be considered as a particular case of the above-discussed coupled equations  \eqref{tri-polyeqa}  with $\alpha_s=\alpha,\alpha\in(0,1],$   $ q_{s1}=h_{20}=d_{10}=g_{20}=\eta_{10}=0,$ and $\kappa_{s0}=q_{s0}=p_{s0}=c_{s0}=\tau_{0},s=1,2,$  for which we have obtained the {analytical} solution in the generalized separable form  \eqref{exactsolution3} {using the developed invariant subspace method.} Also, Yang et al. \cite{yang2019} have found an approximate solution for the given coupled equations \eqref{2+1diffyang} using the Jacobi special collocation method along with the reduction $r=\sqrt{x_1^2+x_2^2}.$
				}
				{Also, we note that the obtained  {analytical} solution  \eqref{exactsolution3} for the considered} coupled equations \eqref{tri-polyeqa} is valid for all $\alpha_s, \alpha_s\in(0,2], s=1,2. $
		{Additionally, we observe that   the obtained solutions \eqref{exactsolution3} coincide with integer-order cases of the given coupled system \eqref{tri-polyeqa} when $\alpha_s=1$ and $\alpha_s=2$. }
						\end{example}
			%
			%
			%
			%
			%
			%
			\section{Discussions and Concluding remarks}\label{conclusion}
			\subsection{Discussions}
			Diffusion is one of the most fundamental phenomena that occur naturally in physical, biological, chemical and ecology systems. With the development of science and technology, the need to study more complicated processes, usually involving the long-time interactions of more than one particle, chemical agent or species has become more significant. In particular, the flow of fluids in a porous media, dispersive sedimentation on the surfaces, and catalytic chemical reactors are studied extensively using the class of coupled  diffusion-convection equations
			\cite{espinov,rossi2012,fletcher1983,pandey2019}.
			One of the revolutionary investigations was made by Turing \cite{turing1952} to formulate the various pattern formations in organisms using coupled diffusion systems.
			Later, the interest in complex processes governed by more than one particle described mathematically using the class of coupled  diffusion equations increased significantly \cite{langlands2008,datsko2012,datsko2018,gafiychuk2009}.
			The convective flux in  such a complex process has an intrinsic effect whenever an external force is present. For example, considering the effect of gravitational force in  {the porous medium fluids dynamics}, the coupled system of diffusion-convection equations plays a vital role in  studying  the transport of contamination in groundwater \cite{rossi2012}.
			
			One of the most celebrated coupled  diffusion-convection equations is the coupled Burger's equation, which has many applications in the physical context. The mathematical formulation of the  $(1+1)$ dimensional coupled  Burger's equation \cite{espinov}  having the form \begin{eqnarray}
				\dfrac{\partial u_i}{\partial t}= v_i\dfrac{\partial}{\partial x}[(1-k_{i1}u_1+k_{i2}u_2)u_i] +d_i\dfrac{\partial^2 u_i}{\partial x^2},
			\end{eqnarray}
			which  Esipov \cite{espinov} proposed initially to describe the poly-dispersive sedimentation on rocks under gravitational force using the following governing equations based on Fick's law:
			\begin{eqnarray*}
				\begin{aligned}(i)&
					\text{ the continuum equation: }
					&
					\dfrac{\partial\mathbf{U}}{\partial t} +\dfrac{\partial\mathbf{q}}{\partial x}=0, \text{ and }
					\\(ii)
					& \text{ the constitutive equation:}
					&{q}_i= V_i(\mathbf{U})-D_{ij}(\mathbf{U}) \dfrac{\partial u_j}{\partial x},
				\end{aligned}
			\end{eqnarray*}
			where $\mathbf{q}=(q_1,q_2), q_i$ is the $i^{th}$ particle flux and   $ \mathbf{U}=(u_1,u_2),u_i=u_i(x,t) $ is the $i^{th}$ particle concentration $ i=1,2. $
			Here, $V_i(\mathbf{U})$ denotes the concentration-dependent velocity, and  $D_{ij}(\mathbf{U})$ is the diffusion coefficient for each $i,j=1,2,$ which were taken in particular as $ -v_i(1-k_{i1}u_1+k_{i2}u_2)u_i$ and $d_i,$  respectively for $v_i,d_i,k_{i1},k_{i2}\in\mathbb{R},i=1,2.$
			Due to its broad applicability, the continuous interest in the class of coupled systems of diffusion-convection equations led to the study of {analytical} solutions for more complex systems of this kind, even in higher dimensions.
			For instance, Fletcher \cite{fletcher1983} {studied} the {analytical} solutions for the $(2+1)$-dimensional coupled  Burger's equations  of the form
			\begin{eqnarray}
				\begin{aligned}
					\dfrac{\partial u_1}{\partial t}=& \dfrac{1}{Re}\sum\limits_{r=1}^2\dfrac{\partial^2 u_1}{\partial x_r^2} -u_1\dfrac{\partial u_1}{\partial x_1}-u_2\dfrac{\partial u_1}{\partial x_2}
					\\
					\dfrac{\partial u_2}{\partial t}=& \dfrac{1}{Re}\sum\limits_{r=1}^2\dfrac{\partial^2 u_2}{\partial x_r^2} -u_1\dfrac{\partial u_2}{\partial x_1}-u_2\dfrac{\partial u_2}{\partial x_2},
				\end{aligned}\label{2+1Burgerseqn}
			\end{eqnarray}
			where $Re$ represents the Reynold's number, and $u_s=u_s(x_1,x_2,t),x_r\in\mathbb{R},t>0.$ The  above-coupled  equations \eqref{2+1Burgerseqn} coincides with the well-known Navier-Stokes equations for incompressible fluids without pressure gradient terms. These are used to study the properties of various  hydro-dynamical models and chemical reactors, to name only a few.
			Using the local discontinuous Galerkin method, an extensive study on different properties of the most generalized coupled  diffusion-convection equations in the form of  \eqref{DGM} was discussed by Cockburn and Shu \cite{cockburn1998}.
			
			However, the anisotropic or heterogeneous media diffusion processes challenged the scientists with their unpredictable nature, commonly called anomalous diffusion \cite{lenzi2016,evangelista2018}. Such anomalous diffusion has the characteristic property of nonlinear mean square displacement in time. More precisely, $\left\langle (x-< x >)^2\right\rangle \propto t^\alpha,\alpha>0,$ where $x$ is the displacement and $t$ is the time \cite{lenzi2016,evangelista2018}. Povstenko \cite{Povstenko2015} considered the
			nonlocal time-dependence related to the matter flux $\mathbf{Q}$ as follows:
			$$	\mathbf{Q}=\left\{\begin{array}{ll}
				\dfrac{1}{\Gamma (\alpha)} \dfrac{\partial}{\partial t} \int_{0}^t(t-y)^{\alpha-1}(-\kappa\bigtriangledown u )dy,  \text{ if }\alpha\in(0,1],
				\\
				\dfrac{1}{\Gamma (\alpha-1)}  {\int_{0}^t}(t-y)^{\alpha-2} (-\kappa\bigtriangledown u) dy, \text{ if } \alpha\in(1,2].
			\end{array}\right.
			$$
			The above-equation gives the time-fractional diffusion-wave equation of the form
			\begin{equation}
				\label{2007povstenko}
				\dfrac{\partial^\alpha u}{ \partial t^\alpha}= \kappa\bigtriangleup u,\alpha\in(0,2],
			\end{equation}
			where $\kappa$ is the diffusive constant, $u=u(x_1,\dots,x_N,t),x_r\in\mathbb{R},$  $t>0,$ and $	\dfrac{\partial^\alpha }{ \partial t^\alpha}(\cdot)$ represents the Caputo fractional derivative \eqref{caputo}.
			We  observed that the above equation  \eqref{2007povstenko} is used to describe many real-world processes in semiconductors, chemical reactors, biological systems, thermoelastic materials, viscous-elastic materials, hydrothermal medium, pattern formations, and so on \cite{Mainardi1997,Tarasov2013a,Tarasov2011,ionescu2017,bagley1984,Tarasov2013b,Povstenko2015,evangelista2018}.
			Moreover, {to investigate the heat conduction phenomena in an infinite dimensional composite medium, Povstenko
			\cite{povstenko2013} used  the time-fractional {multi-component} diffusion-wave (heat) equations of the form}
			\begin{equation}
				\label{systempovstenko}
				\dfrac{\partial^{\alpha_s}
					u_s}{ \partial t^{\alpha_s}}= \kappa_s\dfrac{\partial^2u_s}{\partial x^2} ,\alpha_s\in(0,2],s=1,2,
			\end{equation}
			where $\kappa_s$ is the heat conduction constant, and $u_s=u_s(x,t)$ is the temperature of the $s$-th surface, $x\in\mathbb{R},$  $t>0,$ and $	\dfrac{\partial^{\alpha_s} }{ \partial t^{\alpha_s}}(\cdot)$ represents $\alpha_s$-th order of the Caputo {time-fractional derivative} \eqref{caputo}. {Lenzi \cite{lenzi2016}  used the {multi-component} equations \eqref{systempovstenko} to describe the transport  and anomalous diffusion separated by a
				membrane in heterogeneous
				systems.}
			In many complex physical phenomena, the convection terms appear naturally along with the anomalous diffusion processes.
			Moreover, long-time interactions in systems with more than one species can be described using the class of system of time-fractional diffusion equations \cite{lenzi2016,povstenko2013,langlands2008,datsko2012,datsko2018}.
			The class of coupled  {TFDCWEs}, with or without including source/sink terms, have  {many potential applications in engineering,  biology, chemistry, and  physics, primarily in the fluid flow problems in porous media.} The approximate analytical solution for  the $(1+1)$-dimensional coupled nonlinear time-fractional  diffusion-convection equations with source term {by imposing the Dirichlet boundary conditions} were investigated in \cite{pandey2019}.

			{In general, finding {analytical} solutions for  coupled nonlinear fractional diffusion equations is a challenging and difficult task.
			However,  analytical methods like the invariant subspace method \cite{sahadevan2017a,choudary2019a,choudary2019b,prakash2020a,prakash2024} and   Lie symmetry method \cite{sahadevan2017b,cherniha2005,bindu2001,bindu2004} were proven to be helpful in finding  {analytical} solutions of  {coupled and scalar}   nonlinear diffusion equations and   coupled  nonlinear diffusion-convection equations in $(1+1)$-dimensions in both integer-order and fractional-order cases.
		A few analytical studies only exist in the literature for studying    {analytical} solutions of   coupled  NTFPDEs, particularly in higher dimensions \cite{choudary2019a}.}
				{
				Additionally, we wish to point out that we can expect the following kinds of {analytical} solutions 	
				$\mathbf{U}=\big( u_1(x_1,x_2,t),u_2(x_1,x_2,t)\big) $  for the {multi-component} $(2+1)$-dimensional  coupled  nonlinear integer-order diffusion-convection equations  as
	\begin{eqnarray}
		\begin{aligned}
	(a)	\,& \text{Functional separable analytical solution: }\mathbf{U}=\left( \Phi_1(z),\Phi_2(z)\right), 
	\\& \text{where }z=\sum\limits^{n_2}_{i_2=1}\sum\limits^{n_1}_{i_1=1}\varphi_{i_1,i_2}(t)\prod_{j=1}^{2}\psi_{i_j}(x_j).
\\(b)\,& \text{Travelling-wave analytical solution: }
\mathbf{U}=\left( \Psi_1(z),\Psi_2(z)\right), \\& \text{where }z=k_1t+k_2x_1+k_3x_2 ,k_r\in\mathbb{R},r=1,2,3.
 		\\
(c)\, &
\text{Generalized separable analytical solutions: }	\\& \begin{array}{ll}
\text{(I)}\,&	\mathbf{U}=\Big(  \sum\limits^{n_1}_{i=1}\varphi^1_{i}(t) \psi^1_{i }(z),
\sum\limits^{n_2}_{i=1}\varphi^1_{i}(t) \psi^1_{i }(z)
\Big) ,z= k_1x_1+k_2x_2 ,k_r\in\mathbb{R},r=1,2,
	\\	\text{(II)}\,& 	\mathbf{U}=\Big(  \sum\limits^{n_{ 2}}_{i_2=1}\sum\limits^{n_{ 1}}_{i_1=1}\varphi^1_{i_1,i_2}(t)\prod\limits_{j=1}^{2}\psi^1_{i_j}(x_j),
	\sum\limits^{m_{2}}_{i_2=1}\sum\limits^{m_{ 1}}_{i_1=1}\varphi^2_{i_1,i_2}(t)\prod\limits_{j=1}^{2}\psi^2_{i_j}(x_j)
	\Big) .
\end{array}
					\end{aligned}
				\end{eqnarray}
				Here $n_j,m_{j}\in\mathbb{N},$ the unknown functions $\varphi_{i_1,i_2}(t),\varphi^s_{i_1,i_2}(t)$ and $\psi_{i_j}(x_j),\psi^s_{i_j}(x_j),s,j=1,2,$   are to be determined, which depend on the considered coupled system.
				We can obtain {solutions of the form (a), (b) and (c)} for the multi-component $(2 + 1)$-dimensional coupled nonlinear integer-order diffusion-convection equations through  the  functional and  generalized separable methods, respectively. However, for fractional-order case in the sense of the Riemann-Liouville and Caputo fractional-order derivatives, the functional separable solutions (a) and travelling-wave analytical solutions (b) are not possible for    the    
				the {multi-component} $(2+1)$-dimensional coupled  {TFDCWEs} \eqref{cdsystemeqn}   with singular kernel   fractional-order derivatives.   Hence, we can expect only generalized separable solutions   (c) for the {multi-component} $(2+1)$-dimensional coupled  {TFDCWEs}  \eqref{cdsystemeqn} under singular kernel fractional derivatives. So, in this work, we discussed the generalized separable solutions   (c)-(II) for the given {multi-component} $(2+1)$-dimensional coupled  nonlinear {TFDCWEs} \eqref{cdsystemeqn} under the Caputo fractional derivative through the invariant subspace method.} {We also  note that only  a few papers} studied the {analytical} solutions for the coupled  NTFPDEs in higher dimensions.
				For the first time, we developed the theory of invariant subspace method to find the {analytical} solutions of the {multi-component} $(N+1)$-dimensional coupled  NTFPDEs. In this work, we also studied how to find {analytical} solutions of a more general class of {multi-component} $(2+1)$-dimensional coupled  nonlinear {TFDCWEs} {with suitable initial-boundary conditions}. Here, we used the Dirichlet boundary conditions for the {multi-component} $(2+1)$-dimensional coupled  {TFDCWEs} defined in a bounded domain.
{Furthermore, we would like to emphasize that the invariant subspace method discussed in this article is applicable to all kinds of  fractional derivatives, which are available in the literature, not only the Caputo fractional derivative. Also, we note that some interesting investigations have been made in the existing literature for the invariant subspace method to find generalized separable analytical solutions for fractional PDEs in $(1+1)$-dimension under different fractional derivatives, such as the Riemann-Liouville fractional derivative \cite{gazizov2013, choudary2019b,choudary2019a,prakash2020a}, and the Regularized Prabhakar fractional derivative  \cite{chu2022}. In   future, we will study the details of applications for the developed invariant subspace method to obtain generalized separable analytical solutions for higher-dimensional coupled systems of time-fractional PDEs under different types of fractional derivatives.}

			In the following subsection, we summarize the importance of the developed method and  the results obtained from this work.
\subsection{Concluding remarks}
			In the present work, we have systematically developed the generalization of the invariant subspace method for finding the invariant  product linear space and generalized separable {analytical} solutions of {multi-component} $(N+1)$-dimensional coupled  NTFPDEs \eqref{Nd1.1}  for the first time.
			Explicitly, we have provided a detailed description for constructing {the various kinds of type-$p$} {invariant product  linear spaces} for the coupled NTFPDEs  given in \eqref{Nd1.1} for {$p=1,2,\dots,N.$}
			Also, we wish to mention here that one can expect $ N $ different types of generalized separable {analytical} solutions for the coupled NTFPDEs \eqref{Nd1.1} {based on the obtained  type-$p$  invariant product  linear spaces, {$p=1,2,\dots,N.$}}
			Additionally, {we have explained how to derive}  {analytical} solutions for the given coupled NTFPDEs \eqref{Nd1.1} by reducing it into a system of FODEs {using the invariant product linear spaces determined}.
			More precisely, we successfully illustrated the efficacy of the developed method through the well-known  {multi-component} $(2+1)$-dimensional coupled  nonlinear {TFDCWEs} \eqref{cdsystemeqn}  and constructed their invariant  product linear spaces.
			Using the {type-1 and type-2 invariant product linear spaces determined for coupled equations  \eqref{cdsystemeqn}}, we have explained {how to obtain} the {analytical} solutions of the coupled equations  \eqref{cdsystemeqn} with the appropriate initial conditions.
			Finally, we wish to point out that one can expect two types of {analytical} solutions for the coupled equations \eqref{cdsystemeqn} based on the {type-1 and type-2 invariant product linear spaces}.
		
{Additionally, we wish to point out that the discussed invariant subspace method is applicable for finding the analytical solutions of  scalar and coupled systems of nonlinear PDEs involving both time and space fractional derivatives.
	 So, we can use this method to find analytical solutions for coupled system of time and space fractional diffusion-convection-wave equations under different fractional-order derivatives. 
		In addition, we wish to point out that analytical solutions of scalar and coupled time-space fractional nonlinear PDEs have been derived through the invariant subspace method in  \cite{choudary2019b,prakash2020a}. However, they have not developed the systematic algorithmic approach of the invariant subspace method to obtain analytical solutions for scalar and coupled time-space fractional nonlinear PDEs in the literature. So, in the future, we plan to develop systematically the algorithmic approach of the invariant subspace method for finding the analytical solutions of scalar and coupled time-space fractional nonlinear PDEs.
}	
We hope that the discussed systematic investigation will help to study the exact descriptions of many mutually connected real-world complex phenomena more accurately.
			This study will help to prove the effectiveness of the invariant subspace method, which is a {powerful systematic and analytical tool} for finding the generalized separable {analytical} solutions for complex  {coupled  nonlinear fractional PDEs.}
			We strongly believe that due to the  {nonstandard properties} of singular kernel fractional-order derivatives,  the present development of the invariant subspace method will be very useful in constructing various kinds of generalized separable {analytical} solutions for  {coupled and scalar  nonlinear fractional PDEs arising in various areas of science and engineering. }
		
		\section*{Acknowledgments}
			The first author is thankful for the financial assistance received  in the form of International Mathematical Union Breakout Graduate fellowship-2023 (IMU-BGF-2023-06) provided by the  IMU, Germany. Another author (M.L.) receives support from the Department of Science and Technology, India-SERB National Science Chair position (NSC/2020/000029).
			\section*{Declaration of competing interest}
		The authors declare hat they have no known competing interest.

\end{document}